\documentclass[11pt,a4paper]{book}
\usepackage{amsmath,amssymb,amsthm,amsfonts}
\usepackage{latexsym}
\usepackage{makeidx}

\usepackage{lscape}
\usepackage{graphics}

\renewcommand{\chaptermark}[1]%
        {\markboth{#1}{}}
\renewcommand{\sectionmark}[1]%
        {\markright{#1}}

\newtheorem{thm}{Theorem}[section]

\newtheorem*{thma}{Theorem~\ref{thm:no2}}
\newtheorem*{thmcA}{Theorem A}
\newtheorem*{thmcB}{Theorem B}
\newtheorem*{thmcC}{Theorem C}
\newtheorem*{thmcD}{Theorem D}
\newtheorem*{thmcE}{Theorem E}

\newtheorem{lem}[thm]{Lemma}
\newtheorem{cor}[thm]{Corollary}
\newtheorem{const}[thm]{Construction}
\newtheorem{prop}[thm]{Proposition}
\newtheorem{deff}[thm]{Definition}
\newtheorem{rk}[thm]{Remark}
\newtheorem{examp}[thm]{Example}
\newtheorem{quest}[thm]{Question}
\newtheorem{claim}{Claim}
\newtheorem{prob}[thm]{Problem}

\newcommand{\brac}[1]{\langle #1\rangle}
\newcommand{\qedd}{\hfill \qedsymbol \vspace{1.5mm}}
\newcommand{\vs}{\vspace{0.3cm}}
\newcommand{\vsl}{\vspace{0.1cm}}

\newcommand{\0}{\mathbf{0}}

\newcommand{\Or}{\mathcal{O}}

\newcommand{\normleq}{\trianglelefteq}
\newcommand{\normgeq}{\trianglerighteq}
\newcommand{\normlt}{\vartriangleleft}
\newcommand{\normgt}{\vartriangleright}
\newcommand{\Z}{{\ensuremath{\mathbb{Z}}}}
\newcommand{\Fp}{\mathbb{F}_p}
\newcommand{\F}[1]{\mathbb{F}_{#1}}
\newcommand{\M}[1]{\mathcal{#1}}
\newcommand{\al}{\alpha}
\newcommand{\oline}[1]{\overline{#1}}

\newcommand{\ifff}{\Longleftrightarrow}
\newcommand{\imply}{\Longrightarrow}
\newcommand{\mapp}{\longrightarrow}
\newcommand{\om}{\omega}
\newcommand{\omb}{\widehat{\omega}}

\newcommand{\pone}{\widehat{1}}
\newcommand{\mone}{\widehat{-1}}
\newcommand{\FF}{\mathbb{F}}
\newcommand{\EE}{\mathbb{E}}
\newcommand{\ppd}{p.p.d.}
\newcommand{\GP}{GPaley(p^R, \frac{p^R-1}{k})}
\newcommand{\GPq}{GPaley(q, \frac{q-1}{k})}
\newcommand{\GPP}[2]{GPaley(#1, \frac{#1 - 1}{#2})}

\newcommand{\Magm}{\textsc{Magma  }}

\newcommand{\Atlas}{\textsc{Atlas  }}

\newcommand{\bpf}{\begin{proof}}
\newcommand{\epf}{\end{proof}}

\newcommand{\bpro}{\begin{prob}}
\newcommand{\epro}{\end{prob}}

\newcommand{\bt}{\begin{thm}}
\newcommand{\et}{\end{thm}}

\newcommand{\btaa}{\begin{thma}}
\newcommand{\etaa}{\end{thma}}

\newcommand{\btcA}{\begin{thmcA}}
\newcommand{\etcA}{\end{thmcA}}

\newcommand{\btcB}{\begin{thmcB}}
\newcommand{\etcB}{\end{thmcB}}

\newcommand{\btcC}{\begin{thmcC}}
\newcommand{\etcC}{\end{thmcC}}

\newcommand{\btcD}{\begin{thmcD}}
\newcommand{\etcD}{\end{thmcD}}

\newcommand{\btcE}{\begin{thmcE}}
\newcommand{\etcE}{\end{thmcE}}

\newcommand{\bl}{\begin{lem}}
\newcommand{\el}{\end{lem}}

\newcommand{\bp}{\begin{prop}}
\newcommand{\ep}{\end{prop}}

\newcommand{\bc}{\begin{cor}}
\newcommand{\ec}{\end{cor}}

\newcommand{\bct}{\begin{const}}
\newcommand{\ect}{\end{const}}

\newcommand{\bdeff}{\begin{deff}}
\newcommand{\edeff}{\end{deff}}

\newcommand{\bve}{\begin{verbatim}}
\newcommand{\eve}{\end{verbatim}}

\newcommand{\brk}{\begin{rk}}
\newcommand{\erk}{\end{rk}}

\newcommand{\bexp}{\begin{examp}}
\newcommand{\eexp}{\end{examp}}

\newcommand{\bq}{\begin{quest}}
\newcommand{\eq}{\end{quest}}

\newcommand{\bcl}{\begin{claim}}
\newcommand{\ecl}{\end{claim}}

\newcommand{\be}{\begin{eqnarray*}}
\newcommand{\ee}{\end{eqnarray*}}

\newcommand{\ben}{\begin{eqnarray}}
\newcommand{\een}{\end{eqnarray}}

\newcommand{\bi}{\begin{itemize}}
\newcommand{\ei}{\end{itemize}}

\newcommand{\bnum}{\begin{enumerate}}
\newcommand{\enum}{\end{enumerate}}

\begin{document}

\pagenumbering{roman}

\title {\bf Edge-Transitive Homogeneous Factorisations of Complete Graphs}


\author
{\textbf{Tian Khoon Lim}  \\
B. Math (Hons) Newcastle, NSW \\[15ex]
{\small This thesis is presented for} \\
{\small the Degree of Doctor of Philosophy of} \\
{\small The University of Western Australia} \\[15ex]
{\small School of Mathematics and Statistics (M019)} \\
{\small The University of Western Australia} \\
{\small 35 Stirling Highway, Crawley, WA 6009} \\
{\small Australia} }

\date{April 28, 2004}
\maketitle

\chapter*{Abstract}

\addcontentsline{toc}{chapter}{\textbf{Abstract}}

This thesis concerns the study of homogeneous factorisations of complete graphs with
edge-transitive factors. A factorisation of a complete graph $K_n$ is a partition of its
edges into disjoint classes. Each class of edges in a factorisation of $K_n$ corresponds
to a spanning subgraph called a factor. If all the factors are isomorphic to one
another, then a factorisation of $K_n$ is called an isomorphic factorisation. A
homogeneous factorisation of a complete graph is an isomorphic factorisation where there
exists a group $G$ which permutes the factors transitively, and a normal subgroup $M$ of
$G$ such that each factor is $M$-vertex-transitive. If $M$ also acts edge-transitively
on each factor, then a homogeneous factorisation of $K_n$ is called an edge-transitive
homogeneous factorisation. The aim of this thesis is to study edge-transitive homogeneous
factorisations of $K_n$. We achieve a nearly complete explicit classification except for
the case where $G$ is an affine $2$-homogeneous group of the form $Z_p^R \rtimes G_0$,
where $G_0 \leq \Gamma L(1,p^R)$. In this case, we obtain necessary and sufficient
arithmetic conditions on certain parameters for such factorisations to exist, and give a
generic construction that specifies the homogeneous factorisation completely, given that
the conditions on the parameters hold. Moreover, we give two constructions of infinite
families of examples where we specify the parameters explicitly. In the second infinite
family, the arc-transitive factors are generalisations of certain arc-transitive,
self-complementary graphs constructed by Peisert in 2001.

\chapter*{Acknowledgements}

\addcontentsline{toc}{chapter}{\textbf{Acknowledgements}}

There are many people whom I wish to express my utmost gratitude for making this thesis
possible.

Foremost, I would like to thank my supervisor Professor Cheryl E. Praeger for her
guidance and encouragement throughout the work on this thesis. She has been an
inspiration to me. I am most grateful to her for introducing me to the area of algebraic
graph theory, for sharing with me her vast knowledge and for providing me with various
opportunities.

I also want to thank my co-supervisor Associate Professor Cai Heng Li for his support and many
invaluable suggestions during the course of my study. His willingness to answer my various
questions is greatly appreciated.

I will always be thankful to Dr. James A. MacDougall, my honours supervisor at the University of
Newcastle (NSW), who first introduced me to mathematical research and showed me how much fun and
satisfying the experience can be.

I am greatly indebted to my parents for their support. Although deprived of a decent
education at a very young age, they never cease to encourage me in pursuing a higher
education. I am grateful to my two sisters, Sok Ting and Sok Khim, who while I am away
from home, have provided me with the ease of mind that my parents are being well taken
care of, especially during the time when Singapore was hit with SARS. I also want to
thank my girlfriend Sandra Tan for her love, patience and constant support during the
course of my study in UWA.

I would like to thank Mr. Daryl Foster, principal of Currie Hall, under whom I have the
privilege to work as a residential adviser. His warmth and friendliness will be fondly
remembered. Also, to the many friends I have got to know throughout my three and a half
years in Currie Hall, I would like to say a big ``thank you" --- for all your friendship
and a great sense of humour.

I want to acknowledge the support of a PhD completion scholarship from the University of
Western Australia during the final semester of my research study.

Last but not least I would like to thank my Lord and Saviour, Jesus Christ, who has
blessed me with so much, and saved me from that which I could not save myself.

\chapter*{List of publications arising from this thesis}

\addcontentsline{toc}{chapter}{\textbf{List of publications arising from this thesis}}

\cite{LLP} C. H. Li, T. K. Lim and C. E. Praeger, Homogeneous factorisations of complete
graphs with edge-transitive factors, \textit{in preparation}.

\vspace{10mm}

\noindent \cite{Lim} T. K. Lim, Arc-transitive homogeneous factorisations, Hamming graphs and
affine planes, \textit{submitted} (also in Research Report No. 2003/21, U. Western Australia).

\vspace{10mm}

\noindent \cite{LP} T. K. Lim and C. E. Praeger, On generalised Paley graphs and their
automorphism groups, \textit{preprint}.

\tableofcontents


\pagenumbering{arabic}

\setcounter{page}{1}


\chapter{Introduction} \label{c1}


We begin by defining what we mean by a homogeneous factorisation of a complete graph and
pose the central question addressed in this thesis. In addition, we also briefly review
the study of factorisations of complete graphs, and provide some background to the
research development that motivates the study of homogeneous factorisations of graphs.
Finally, we give an overview of the results achieved in this thesis.

\section{Introduction} \label{c1:intro}
The study of graph factorisations has long been a very important topic in graph theory.
Roughly speaking, a factorisation of a graph $\Gamma$ is an expression of $\Gamma$ as an
edge-disjoint union of subgraphs called factors. Here, we are only interested in
factorisations of complete graphs. (We refer readers to Chapter~\ref{c2} of this thesis
and other standard texts \cite{Biggs_b93, Cameron_b99, DM_b96, GR} in permutation groups
and algebraic graph theory for definitions, notations and terminology not fully explained
in this chapter.)

Let $K_n = (V,E)$ be a complete graph on $n$ vertices with vertex set $V$ and edge set $E =
\{\{x,y\} \mid x\neq y \ \mbox{and} \ x ,y \in V\}$. A \textit{factorisation} of $K_n$ is a
partition of its edges into disjoint classes $\mathcal{E}:=\{E_1, \ldots, E_k\}$. Each class of
edges in a factorisation of $K_n$ corresponds to some spanning subgraph called a \textit{factor}
$\Gamma_i =(V,E_i)$, with vertex set $V$ and edge set $E_i$. If all the factors are isomorphic to
one another, then the factorisation of $K_n$ is called an \textit{isomorphic factorisation}. The
following is an example of an isomorphic factorisation of a complete graph.

\begin{figure}[ht]

\begin{center}
\includegraphics{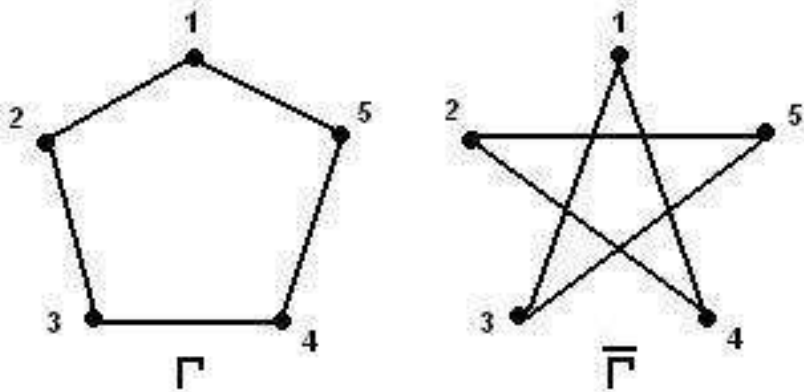}

\caption{Isomorphic factorisation of $K_5$ into $2$ factors.} \label{pic}
\end{center}


\end{figure}

\bexp \label{ex:hom} \rm{Figure~\ref{pic} shows an isomorphic factorisation of $K_5=(V,E)$ into 2
factors $\Gamma$ and $\oline{\Gamma}$ (where $\Gamma \cong \oline{\Gamma} =C_5$ and $V=\{1,2,3,4,5
\}$). Notice that there are two groups associated with this example. First, we have $M =D_{10}
=Aut(\Gamma)= Aut(\oline{\Gamma})$ ($D_{10}$ is the dihedral group of order $10$). Then there is
$G =F_{20}$ (the Frobenius group of order $20$), which contains $M=D_{10}$ as a normal subgroup,
and interchanges $\Gamma$ and $\oline{\Gamma}$. } \eexp

Inspired by the beautiful examples of self-complementary vertex-transitive graphs (a graph that is
isomorphic to its complement is said to be \textit{self-complementary}, see Section~\ref{c2:graphs}
for a detailed definition) such as the one shown in Figure~\ref{pic}, Li and Praeger in
\cite{LP2003} introduced the concept of homogeneous factorisations of complete graphs. We shall
now define what we mean by a homogeneous factorisation.

\bdeff \label{def:homfac} \rm{\textbf{(Homogeneous factorisation)} \ A
\textit{homogeneous factorisation} of a complete graph $K_n=(V,E)$ is defined to be an
isomorphic factorisation $\mathcal{E}=\{E_1, \ldots, E_k\}$ (where $k >1$) such that
there exists a transitive subgroup $G \leq$ Sym$(V)$ which leaves the partition
$\mathcal{E}$ invariant, permutes the parts $E_i$ transitively, and the kernel $M$ of
this action acts transitively on $V$. A homogeneous factorisation of a complete graph is
said to be of \textit{index} $k$ if the edge set is partitioned into $k$ parts. We shall
use the 4-tuple $(M,G,V,\mathcal{E})$ to denote a homogeneous factorisation of $K_n
=(V,E)$.} \edeff

This thesis is dedicated to studying homogeneous factorisations of complete graphs with
edge-transitive (undirected) factors. More precisely, we look at the situation where
there exists a group $G$ that is transitive on the factors of an isomorphic
factorisation $\mathcal{E}$ of $K_n=(V,E)$ with $M$, the kernel of this $G$-action on
$\mathcal{E}$, acting transitively on both vertices and edges of each factor. Thus we
have $M < G \leq$ Sym$(V)$ such that: \bi
\item For all $E_i \in \mathcal{E}=\{E_1, \ldots, E_k\}$
and $g \in G$, we have $E_i^g \in \mathcal{E}$.
\item For any $E_i$ and $E_j$, there exists $g \in G$ such that
$E_i^g =E_j$.
\item For each $i$, $M$ fixes $E_i$ setwise and is transitive on $V$ and on $E_i$.
\ei We shall call such a factorisation an \textit{edge-transitive homogeneous factorisation} of a
complete graph. We also call $(M,G,V,\M{E})$ an \textit{arc-transitive homogeneous factorisation}
if $M$ is arc-transitive on each factor. For convenience, we will often refer to $(M,G,V,\M{E})$
simply as an edge or arc-transitive homogeneous factorisation. Clearly $M$ is a normal subgroup of
$G$ (since it is the kernel of the $G$-action on $\mathcal{E}$) and $M = G \cap Aut(\Gamma_1) \cap
\ldots \cap Aut(\Gamma_k)$. Furthermore, since the group $G$ permutes the factors transitively, it
also induces isomorphisms between each pair of the graphs $\Gamma_i =(V,E_i)$ ($1 \leq i \leq k)$.
(Note that Example~\ref{ex:hom} is an arc-transitive homogeneous factorisation $(D_{10}, F_{20},
V, \M{E})$ of $K_5$ of index $2$, where $V=\{1,2,3,4,5\}$ and $\M{E} =\{E\Gamma_1, E\Gamma_2\}$
(see Figure~\ref{pic}).)

There are many interesting questions that we can ask pertaining to edge-transitive
homogeneous factorisations of complete graphs. However, our point of departure with
regards to this topic will be the following natural question.

\bq \label{quest:main} \rm{Given an edge-transitive homogeneous factorisation of a complete graph,
what can we say about its (undirected) factors? } \eq

The focus of this thesis will be centered around the above general question. In Chapter~\ref{c6},
we first study the structure of an infinite class of graphs called the generalised Paley graphs.
We will prove in Chapter~\ref{c3} that they arise as factors of arc-transitive homogeneous
factorisations of complete graphs. Also in Chapter~\ref{c3}, we will discuss the general approach
used to investigate the structure of the factors of edge-transitive homogeneous factorisations
$(M,G,V,\M{E})$ of $K_n=(V,E)$. In subsequent chapters (Chapters~\ref{c4} to \ref{c5}), we will
determine the possible groups $M$ and $G$ in such factorisations, characterise the edge-transitive
homogeneous factors that arise, and give a nearly complete and explicit classification. In the
case where we do not have an explicit classification (Chapter~\ref{c5}), we are also able to give
two infinite families of examples. A more detailed introduction to the main results of this thesis
will be given in Section~\ref{c1:m_results}.

We now give a brief review to the study of graph factorisations and discuss the
motivation behind the investigations of edge-transitive homogeneous factorisations of
complete graphs.

\section{Brief literature review} \label{c1:lit_rev}
Investigation into the factorisations of graphs have long been an integral part of graph theory.
More important is the study of isomorphic factorisations of graphs and in particular, isomorphic
factorisations of complete graphs, which have received considerable attention (see for example
\cite{HRW78} and \cite{HR85}). Isomorphic factorisations of complete graphs are studied for
various reasons, but one of the main foci of the research so far has been the study of
\textit{one-factorisations} and \textit{two-factorisations} of complete graphs (see for instance
\cite{BL2002} and \cite{Ober}).

A \textit{one-factor} in a graph is a set of edges no two of which are adjacent (that
is, they do not share a common vertex) and each vertex of the graph is incident to
precisely one of the edges in the set. A one-factorisation of a graph is a partition of
the edge set into one-factors. One-factorisations of complete graphs have been widely
studied and have implications in many areas, such as tournament scheduling (see \cite[p.
566] {CRC}) and in design theory. We refer readers to the surveys \cite{mend, WalA} for
a detailed account. Two-factorisations of $K_n$ coincide with the celebrated
\textit{Oberwolfach problem}, which asks for the existence of factorisations of complete
graphs into isomorphic 2-regular graphs (graphs of regular valency 2) on $n$ vertices,
whose components are simply cycles of length $c_1, \ldots, c_t$ such that $c_i \geq 3$
and $c_1 + \cdots + c_t =n$. Many results had been achieved and we refer readers to
\cite{Brian} for a summary.

Another focus in isomorphic factorisations of complete graphs has been the study of
isomorphic Ramsey numbers. For a graph $\Gamma$ containing no isolated vertices and any
number $k$, the \textit{isomorphic Ramsey number} $F(\Gamma ; k)$ is defined to be the
minimum $m$ such that for every $n \geq m$, $\Gamma$ is isomorphic to a subgraph of every
factor of an isomorphic factorisation of $K_n$ into $k$ parts. It follows that $F(\Gamma
; k) \leq R(\Gamma ; k)$ where $R(\Gamma ; k)$ denotes the $k$-colouring Ramsey number
of $\Gamma$ (which is, the minimum number $m$ such that if we colour the edges of $K_m$
with $k$ colours denoted as $\{C_1, \ldots, C_k\}$, then there exists some $i$ ($1 \leq
i \leq k$) such that $\Gamma$ is isomorphic to some subgraph in $K_m$ coloured by
$C_i$). The most well-studied case is the case where $k=2$ (see \cite{HR85}), and there
are many interesting problems in this area.

The study of homogeneous factorisations of complete graphs was first introduced by Li and Praeger
in \cite{LP2003} (see also \cite{LP2002}). In it, the theory of homogeneous factorisations was
established and is fundamental source to the development of this thesis. The interest is motivated
mainly by the investigations into vertex-transitive self-complementary graphs, which arise as the
factors of homogeneous factorisations of complete graphs of index $2$. The study of
vertex-transitive self-complementary graphs began with a construction of a family of
self-complementary circulant graphs by Sachs in \cite{Sachs}. Since then, vertex-transitive
self-complementary graphs have received much attention, see for example \cite{Jaj, NCay, Rao,
Sup}. Also, a lot of effort has been made to determine the positive integers $n$ such that there
exist vertex-transitive, self-complementary graphs with $n$ vertices (see for example,
\cite{AMV,Li-1,Muzychuk}). In \cite{Muzychuk}, Muzychuk completely determined such positive
integers.

\newpage

One of the major results (Theorem 1.1) of the paper by Li and Praeger \cite{LP2003} is a
generalisation of Muzychuk's result. They proved that a cyclic\footnote{A cyclic homogeneous
factorisation of $K_n$ is one which the group $G$ induces a cyclic permutation group on the
partition $\M{E}$} homogeneous factorisation of index $k$ of $K_n$ (where $n =p_1^{d_1} p_2^{d_2}
\cdots p_m^{d_m}$, $p_i$ are distinct primes and $d_i, m \geq 1$) exists if and only if, for all
$i \in \{1, \ldots, m\}$,
\[
\begin{array}{l}
\mbox{$p_i^{d_i}\equiv1$ $(\mbox{mod } 2k)$ if $p_i$ is odd, and}\\
\mbox{$p_i^{d_i}\equiv1$ $(\mbox{mod } k)$ if $p_i=2$.}
\end{array}
\]
We further note that in \cite{GLPS}, the point primitive case for cyclic homogeneous
factorisations of complete graphs is studied (that is, where $G$ acts primitively\footnote{See
Section~\ref{c2:primitive} for definition of a primitive permutation group.} on the vertices $V$).
They are related to the exceptionality of permutation groups studied in \cite{GMS}.

Our investigation of edge-transitive homogeneous factorisations of complete graphs can also be
traced back to the study of arc-transitive self-complementary graphs (which are factors of
arc-transitive homogeneous factorisations of complete graphs of index $2$). In \cite{Zhang}, H.
Zhang gave an algebraic characterisation of all self-complementary, arc-transitive graphs, and
showed that they are all Cayley graphs. One family of examples is the well-known class of Paley
graphs. In \cite{Peisert2001}, W. Peisert, in an attempt to determine if there are any
self-complementary arc-transitive graphs other than the Paley graphs, completely classified, up to
graph isomorphism, all finite self-complementary arc-transitive graphs (see Theorem~\ref{MainPei}).
In his paper \cite{Peisert2001}, Peisert constructed a second infinite family of
self-complementary arc-transitive graphs that are not isomorphic to the Paley graphs (except for
one case), and called them the $\M{P}^*$-graphs. Like the Paley graphs, the $\M{P}^*$-graphs are
defined on finite fields $\F{p^R}$ where $p \equiv 3$ (mod $4$) and $R$ is even. In fact, by
expanding on Peisert's construction, we are able to construct an infinite family of graphs (where
the $\M{P}^*$-graphs occur as a special case), called the twisted generalised Paley graphs, that
arise as factors of arc-transitive homogeneous factorisations of complete graphs.

We shall present Peisert's result below.

\bt \label{MainPei} A graph $\Gamma =(V,E)$ is self-complementary and arc-transitive if and only
if $|V|=p^R$ for some prime $p$, $p^R \equiv 1$ $($mod $4$$)$, and $\Gamma$ is one of the
following$:$ \bnum \item[$1.$] a Paley graph, or
\item[$2.$] a \textit{$\mathcal{P}^*$-graph}, or
\item[$3.$] a graph $G(23^2)$ of order $23^2$ not isomorphic to a Paley graph or
a \textit{$\mathcal{P}^*$-graph}. \enum  \et

\brk \label{rk:pei} \rm{Along with the graph $G(23^2)$, there are two other exceptional graphs,
namely the $G(7^2)$ graph (on $7^2$ vertices) and the $G(9^2)$ graph (on $9^2$ vertices). Even
though both $G(7^2)$ and $G(9^2)$ are isomorphic to $\mathcal{P}^*$-graphs \cite[Lemmas 6.6 and
6.7]{Peisert2001}, the three graphs $G(7^2)$, $G(9^2)$ and $G(23^2)$ are ``exceptional" in the
sense that unlike the Paley graphs and most of the $\M{P}^*$-graphs, none of their full
automorphism groups is permutationally isomorphic to a one-dimensional affine group \cite[p.
217]{Peisert2001}. \footnote{Note that we will not (for now) give detailed descriptions of the
$\M{P}^*$-graphs and the three exceptional graphs; the reason being that the $\M{P}^*$-graphs occur
as a special case of the twisted generalised Paley graphs constructed in Chapter~\ref{c3}, while
the three exceptional graphs arise from our investigation in Section~\ref{c4:Gqk}.} } \erk

This thesis originated in the work of H. Zhang \cite{Zhang} in 1992 on arc-transitive
self-complementary graphs. However, in our investigation, we use and refine the methodology
developed by Li and Praeger in \cite{LP2003} for studying homogeneous factorisations of complete
graphs. Just as Peisert was able to complete Zhang's work, we achieve an almost complete
classification of arc-transitive homogeneous factorisations of complete graphs of all indices $k$.
Our classification relies on the finite simple group classification through its use of the
classification of finite 2-transitive permutation groups. Moreover, our classification specialises
in the index 2 case to the theorem of Peisert (Theorem~\ref{MainPei}). Finally, on a more
ambitious note, it is hoped that by our study of edge or arc-transitive homogeneous factorisations
of $K_n$, we may be able to shed some light on the computation of lower bounds for some
$k$-colouring Ramsey numbers.

\section{Summary of main results} \label{c1:m_results}
In this section, we will give a summary of the important results achieved in this thesis.
The main results are found in Chapters~\ref{c6} to \ref{c5}. Chapter~\ref{c6} is based
on the paper \cite{LP}, Chapters~\ref{c3}, \ref{c4} and \ref{c5} are based on \cite{LLP}
and Chapter~\ref{c7} is based on \cite{Lim}.

\textbf{Chapter~\ref{c6}} deals with the structure of the class of generalised Paley graphs,
denoted as $\GP$, which arise as factors of arc-transitive homogeneous factorisations of complete
graphs on $p^R$ vertices, where $p$ is prime and $R \geq 1$ (see Definition~\ref{defn:gpaley:c6}).
One of the main results in this chapter says that a connected generalised Paley graph is either a
Hamming graph or a normal Cayley graph where its full automorphism group is an affine primitive
group.

\btcA $($Theorems~$\ref{thm:A}$ and $\ref{thm:no1}$$)$. Let $\Gamma$ $=$ $\GP$ be a connected
generalised Paley graph. Then one of the following holds. \bnum

\item[$1.$] $\Gamma \cong H(p^a, \frac{R}{a})$ $($Hamming graph$)$, for some divisor $a$ of $R$
with $1 \leq a < R$, if and only if $k=\frac{a(p^R-1)}{R(p^a-1)}$.

\item[$2.$] If $k \neq \frac{a(p^R-1)}{R(p^a-1)}$ for any divisor $a$ of $R$ with $1 \leq a
< R$, then $Aut(\Gamma)$ is a primitive subgroup of $AGL(R,p)$ containing the translation group $T
\cong \Z_p^R$. \enum \etcA

Furthermore, we determine explicitly the full automorphism group of a class of generalised Paley
graphs under certain conditions which are relevant to the study of arc-transitive homogeneous
factorisations (see Theorem B below; also see Section~\ref{c6:main} for details on the notations
used). This result will be used in Chapter~\ref{c5} to show the non-isomorphism between the class
of generalised Paley graphs and another infinite family of graphs (the twisted generalised Paley
graphs, see Definition~\ref{defn:tgpaley}) that also arise as factors of arc-transitive
homogeneous factorisations of complete graphs.

\btcB $($Theorem~$\ref{thm:no2}$$)$. Suppose $\Gamma = \GP$ is a generalised Paley graph such that
$k \mid (p-1)$, and if $p$ is odd, then $2k \mid (p^R-1)$. Then $\Gamma$ is connected and
$Aut(\Gamma) = T \rtimes \brac{\omb^k, \al} \leq A\Gamma L(1,p^R)$. \etcB

In the last part of the chapter, we will briefly mention the connection between the
generalised Paley graphs and the class of graphs associated with a symmetric cyclotomic
association scheme.

In \textbf{Chapter~\ref{c3}}, we will develop the general theory of edge-transitive homogeneous
factorisations $(M,G,V,\M{E})$ of $K_n=(V,E)$. In particular, we will show that in such a
factorisation, the group $G$ is necessarily $2$-homogeneous (Lemma~\ref{G_2homo}). Then, using the
classification of finite 2-transitive permutation groups as well as Kantor's description of finite
2-homogeneous but not 2-transitive permutation groups \cite{Kantor69}, we are able to employ a
case-by-case approach in our investigation to classify factors arising from edge-transitive
homogeneous factorisations of $K_n$. We first deal with the case where $G$ is an almost simple
group, that is $G$ has a unique minimal normal subgroup that is nonabelian simple, and show that
there exists only one example of an edge-transitive homogeneous factorisations of $K_n$ where
$n=28$ (Proposition~\ref{psl28}). Then we reduce the problem of classifying factors that arise
from edge-transitive homogeneous factorisations $(M,G,V,\M{E})$, where $G$ is an affine group, to
the problem of classifying factors arising from arc-transitive homogeneous factorisations, where
$G$ is affine 2-transitive. Also in the same chapter, we construct two infinite families of
examples: the generalised Paley graphs $\GP$ and the twisted generalised Paley graphs
$TGPaley(p^R,\frac{p^R-1}{2h})$. Both graphs arise as factors of arc-transitive homogeneous
factorisations $(M,G,V,\M{E})$ of $K_n=(V,E)$, where $G \leq A\Gamma L(1,p^R)$ (and so $n=p^R$),
that is, $G$ is a one-dimensional affine permutation group acting on the finite field $\F{p^R}$.

In \textbf{Chapter~\ref{c4}}, we will look in detail at the list of affine 2-transitive
permutation groups in Theorem~\ref{2_trxclass}. For each class of groups $G$ in that list, we
analyse and find out all possible transitive normal subgroups $M$, and use this information to
determine the structure of the $M$-arc-transitive factors. We show that many of the examples are
disconnected generalised Paley graphs. In other cases, we use \textsc{Magma} \cite{Magma} to
construct the factor graphs (as $M$-orbital graphs) and compute their full automorphism groups. By
doing that, we are able to determine their structures more explicitly. The main result in this
chapter classifies all factors arising from $(M,G,V,\M{E})$, where $G \leq A\Gamma L(a,q)$ is an
affine 2-transitive permutation group with $a \geq 2$ (that is, $G$ is not a one-dimensional
affine 2-transitive group). Combining the result (Proposition~\ref{psl28}) from the previous
chapter for the case where $G$ is almost simple, we have the following.

\btcC $($Proposition~$\ref{psl28}$ and Theorem~$\ref{thm:mainc4}$$)$. Let $(M,G,V,\M{E})$ be an
arc-transitive homogeneous factorisation of $K_{n}=(V,E)$ of index $k$ with factors $\Gamma_i$.
Suppose $G$ is a $2$-transitive permutation group on $V$ such that $G \nleqslant A\Gamma
L(1,p^R)$. Then one of the following holds. \bnum

\item[$1.$] $G= P\Gamma L(2,8)$, $|V| =n =28$, $M=PSL(2,8)$, $k=3$ and each $\Gamma_i$ is a non-Cayley,
$M$-arc-transitive graph of valency $9$.

\item[$2.$] $G=T \rtimes G_0$ is an affine $2$-transitive permutation group on $V=V(a,q)$
$($and hence $n =q^a$$)$ such that $G_0$ is in one of the cases $2(b) - 2(h)$ of
Theorem~$\ref{2_trxclass}$. Then $\Gamma_i =Cay(V,S_i)$, where $|S_i| =\frac{q^a-1}{k}$
and $|S_i|$ is even if $q$ is odd, and $M = T \rtimes M_0$ is such that precisely one of
the following holds. \bnum

\item[$(a)$] $M_0 \leq Z=Z(GL(a,q))$, $(q+1) \mid k$, $\M{E} =\M{E}_{GP}(q^a,k)$, and each
$\Gamma_i$ is a $($disconnected$)$ generalised Paley graph $GPaley(q^a,
\frac{q^a-1}{k})$.

\item[$(b)$] $M_0$ is as in Tables~$\ref{t25}$ and $\ref{t23}$, the $\Gamma_i$ are connected, and
one of the following holds. \bnum \item $\Gamma_i \cong G(q^2,k)$ and $(q,k)=(11,5)$, $(23,22)$,
$(23,11)$, $(23,2)$, $(19,3)$, $(29,7)$ or $(59,29)$ $($the graph $G(q^2,k)$ is as in
Definition~$\ref{def:Gqk}$$)$.
\item $\Gamma_i \cong G(q^2,k)$ where $Aut(\Gamma_i) \lnsim A\Gamma L(1,23^2)$ and $(q,k)=(23,66)$ or $(23,6)$.
\item $\Gamma_i \cong GPaley(q^2, \frac{q^2-1}{k})$ and $(q,k) = (5,3)$, $(11,15)$ or
$(11,3)$ $($note that for $(q,k) = (5,3)$, $\Gamma_i \cong H(5,2)$$)$.
\item $\Gamma_i \cong TGPaley(q^2, \frac{q^2-1}{k})$ and $(q,k) = (9,2)$, $(7,6)$ or $(7,2)$.
\enum

\item[$(c)$] $M_0 =\mathbf{E}$ and $k =5$ and the homogeneous factors $\Gamma_i$ are all
isomorphic to the Hamming graph $H(9,2)$. \enum \enum \etcC (See Chapter~\ref{c4} for further
details on the notations used in Theorem C. Also, part $2(c)$ of Theorem C relies on a result from
Chapter~\ref{c7} which shows that the factors are Hamming graphs.)

In \textbf{Chapter~\ref{c7}} we show, using representation theory, that one example of an
arc-transitive homogeneous factorisation of $K_n$ determined in Chapter~\ref{c4} has factors that
are all isomorphic to the Hamming graph $H(9,2)$. We also show that such a factorisation of a
complete graph (into Hamming graphs) gives rise to a related edge partition that corresponds to an
interesting 2-design (see Theorem~\ref{2design}).

\textbf{Chapter~\ref{c5}} deals with the (one) remaining case not covered in Chapter~\ref{c4}:
where $G$ is a one-dimensional affine 2-transitive permutation group. By expressing the groups $G$
and $M$ using a standard set of parameters introduced by Foulser in \cite{Foulser64} (see also
\cite{FK78}), we give a generic construction (Construction~\ref{generic}) that will yield all
possible arc-transitive homogeneous factorisations $(M,G,V,\M{E})$ where $M,G, \leq
A\Gamma(1,p^R)$ (see Theorem~\ref{thm2}).

We also show that the two infinite families of examples, the generalised Paley graphs
$\GP$ and twisted generalised Paley graphs $TGPaley(p^R,\frac{p^R-1}{2h})$, constructed
in Section~\ref{c3:examples} can be derived from this generic construction. Finally, in
the last part of the chapter, we will prove that, except for one small example, the
graphs from the two infinite families are pair-wise non-isomorphic.

\btcD $($Theorem~$\ref{luvsandra}$$)$. For $R$ even, $p \equiv 3$ $($mod $4)$, and $h$ odd such
that $2h \mid (p-1)$, $GPaley(p^R, \frac{p^R-1}{2h}) \ncong TGPaley(p^R, \frac{p^R-1}{2h})$ except
for the case $(R,p,h)=(2,3,1)$, when we have $GPaley(3^2, 4) \cong TGPaley(3^2, 4)$. \etcD


\chapter{Preliminaries} \label{c2}

This chapter is a collection of basic definitions and well-known results pertaining to
graphs and permutation groups that will be used in subsequent chapters. The results and
notations contained here are often standard and can be found in texts such as
\cite{Biggs_b93, Cameron_b99, DM_b96, GR}.

\section{Permutation groups and group actions} \label{c2:perm_gp}
Let $V$ be a finite set of size $n$. A bijection of $V$ onto itself is called a
\textit{permutation} on $V$. The set of all permutations on a finite set $V$ form a group
under the operation of composition of mappings. We call this group the \textit{symmetric
group} on $V$ and denote it by Sym$(V)$ or S$_n$. Any subgroup $G$ of Sym$(V)$ is said
to be a \textit{permutation group} on $V$, and we say that $G$ is a permutation group of
\textit{degree} $n$ if $|V| =n$. We denote the image of the point $v \in V$ under the
permutation $\al$ by $v^{\al}$. Note that by our definition, a permutation group is
always finite. It is of course possible to extend the definition of a permutation group
to the infinite case, but in this thesis, we will only consider finite permutation
groups and finite groups in general.

Let $G$ be a group and $V$ a finite set. An \textit{action} of $G$ on $V$ is a map $V
\times G \mapp V$ written $(v,g) \mapsto v^g$, such that \bi
\item for every $v \in V$, $v^{1_G} =v$, where $1_G$ denotes the identity element
of the group $G$, and
\item for every $g,h \in G$ and $v \in V$, we have $v^{gh} = (v^g)^h$.
\ei The \textit{kernel} of the action of $G$ on $V$ is defined to be the subgroup of all
elements of $G$ which fix each point of $V$. If this kernel consists of just the identity
element of $G$, then we say that $G$ acts \textit{faithfully} on $V$ (or $G$ is
\textit{faithful} on $V$). Otherwise $G$ is said to act on $V$ \textit{unfaithfully} (or
$G$ is \textit{unfaithful} on $V$).

Note that every subgroup $G$ of Sym$(V)$ acts naturally on $V$. It is easy to see that such an
action is faithful. We will always assume that this is the action we are dealing with whenever we
have a permutation group.

When a group (not necessarily a permutation group) $G$ acts on a set $V$, a typical
element $v \in V$ is moved by elements of $G$ to various other points. The set of these
images is called the \textit{orbit} of $v$ under $G$ or the \textit{$G$-orbit} containing
$v$. We usually denote it by $v^G:= \{v^g \mid g \in G\}$ but sometimes, we will use
$\M{O}$ to denote a typical orbit of $G$ in $V$. Finally, the number of elements in a
$G$-orbit $\M{O}$ is called the \textit{length} of $\M{O}$.

We also call the set of elements in $G$ which fix a specified point $v \in V$, the
\textit{stabiliser} of $v$ in $G$ and denoted it by $G_v :=\{ g \in G \mid v^g =v \}$. Note that
$G_v$ forms a subgroup of $G$ and is often called the \textit{point-stabiliser} of $v$.

A group $G$ acting on $V$ is said to be \textit{transitive} if for any $u, v \in V$,
there is some $g \in G$ such that $u^g=v$. In other words, $G$ is transitive on $V$ if
it has only one orbit, and so $v^G=V$ for all $v \in V$. A group $G$ acting on a set $V$
is said to be \textit{semiregular} if the only element fixing a point in $V$ is the
identity. We say that a group $G$ acting on $V$ is \textit{regular} if $G$ is both
transitive and semiregular on $V$.

Let $G$ be a group and let $V = G$. Then each element $g \in G$ induces a permutation
$\hat{g}$ of $V$ by \textit{right multiplication}:\[ \hat{g}: v \mapp vg \ \mbox{ for
all } v \in V. \] (This is the \textit{right regular representation} of $G$; see for
example \cite[Example 1.3.3]{DM_b96}.) All permutations $\hat{g}$, with $g \in G$, form a
permutation group $\widehat{G} <$ Sym$(V)$ (where $V =G$). Furthermore, $\widehat{G}$ is
isomorphic to $G$ and is a regular subgroup of Sym$(V)$.

Let $G$ be a group acting on $V$ and $H$ be a group acting on $U$. Then we say that $G$
on $V$ is \textit{permutationally isomorphic} to $H$ on $U$ if there exist a group
isomorphism $\theta : G \mapp H$ and a bijection $\xi : V \mapp U$ such that
$(\xi(v))^{\theta(g)} = \xi (v^g)$ for all $g \in G$ and $v \in V$.

\section{$2$-homogeneous and $2$-transitive permutation groups} \label{act}

For a positive integer $k$, we use $V^{(k)}$ to denote the set of $k$-tuples of
\textit{distinct} members of a finite set $V$. Suppose a group $G$ acts on $V$. Then $G$
induces a natural action of $G$ on $V^{(k)}$ defined by \[ (v_1, v_2, \ldots, v_k)^g :=
(v_1^g, v_2^g, \ldots, v_k^g) \] for all $(v_1, v_2, \ldots, v_k) \in V^{(k)}$ and $g
\in G$. If $G$ is transitive on $V^{(k)}$ under this action, then we say that $G$ is
\textit{$k$-transitive} on $V$. Now let $V^{\{k\}}$ be the set of all $k$-subsets (that
is, the subsets of size $k$) of $V$ where $k = 1,2 \ldots, |V|$. Then $G$ also induces a
natural action of $G$ on $V^{\{k\}}$ defined by \[ \{v_1, v_2, \ldots, v_k \}^g := \{
v_1^g, v_2^g, \ldots, v_k^g \} \] for all $\{v_1, v_2, \ldots, v_k \} \in V^{\{k\}}$ and
$g \in G$. The group $G$ is \textit{$k$-homogeneous} if its action on $V^{\{k\}}$ is
transitive. It is not too difficult to see that a $k$-transitive group is necessarily
$k$-homogeneous but the converse is not always true.

As a result of the classification of finite simple groups, all finite $2$-transitive permutation
groups have been classified up to permutation isomorphism (see \cite{Cameron81}). Now if $G$ is
2-transitive, then by a result of Burnside \cite[Section 154]{Burnside_b55}, $G$ has a unique
minimal normal subgroup $N$ where either $N$ is a regular normal subgroup that is elementary
abelian ($G$ is an affine group), or $N$ is a nonabelian simple group ($G$ is almost simple).

A group $G$ is \textit{almost simple} if it has a unique minimal normal subgroup $N$
which is nonabelian and simple. Equivalently, $G$ is almost simple if $N \leq G \leq
Aut(N)$ for some nonabelian simple group $N$. This is true since $C_G(N)=1$, so $G$ can
be embedded in $Aut(N)$ and the image of $G$ under the embedding contains $Inn(N) \cong
N$.

We call a group $G$ acting on a set $V$ \textit{affine} if it has a regular normal
subgroup $N$ which is an elementary abelian $p$-group for some prime $p$, say $\Z_p^R$.
Such a regular normal subgroup is often written additively, and we usually call it a
\textit{translation} subgroup and denote it by $T$ (and so $T \cong N = \Z_p^R$). The set
$V$ is identified with a $R$-dimensional vector space over the prime field $\F{p}$ and
the group $T$ is such that for each $x \in V$, there exists a unique $t_x \in T$ which
acts on $V$ by translation: $t_x: v \mapp v+x$ for all $v \in V$. Each affine group $G$
is a semidirect product $T \rtimes G_0$ where $G_0$, the stabiliser of $\0 \in V=V(R,p)$,
is contained in $\Gamma L(R,p)$. However, if $G_0 \leq \Gamma L(R,p)$ preserves on $V$ a
vector space of dimension $a \geq 1$ over a finite field $\F{q}$ (so $q^a =p^R$), then we
may identify $V=V(a,q) =\F{q}^a$ and $G_0$ as a subgroup of $\Gamma L(a,q)$. If such
case arise, then the affine group $G = T \rtimes G_0$ is identified as a subgroup of $A
\Gamma L(a,q)$ acting on $V =V(a,q)$ where $a \geq 1$ and $q^a=p^R$.

The following contains a list of all finite $2$-transitive permutation groups (up to permutation
isomorphism) compiled from \cite[Theorem 2.9]{GLP98} and \cite{Kantor85} (see also
\cite{Liebeck87}). We further refer readers to \cite[Section 7.7]{DM_b96} for definitions and
descriptions of the nonabelian groups $N$ occurring in Theorem~\ref{2_trxclass} (1).

\bt \label{2_trxclass} Let $G$ be a finite $2$-transitive permutation group on a set $V$ with $|V|
=n$. Then $G$ has a unique minimal normal subgroup $N$ where one of the following holds$:$ \bnum
\item[$1.$] $N$ is one of the following nonabelian simple groups and $N \leq G \leq Aut(N)$$:$ {\rm \bnum
\item $N = A_n$, \textit{of degree $n \geq 5$.}
\item $N = PSL(a,q)$, \textit{of degree $(q^a-1)/(q-1)$ with $(a,q) \neq (2,2)$ or
$(2,3)$.}
\item \textit{$N = PSU(3,q^2)$, of degree $q^3+1$ with $q \geq 3$.}
\item \textit{$N = \hspace{1mm}^2B_2(q)$ $($Suzuki$)$, of degree $q^2+1$, with $q=2^{2c+1} \geq
8$.}
\item \textit{$N = \hspace{1mm}^2G_2(q)$ $($Ree$)$, of degree $q^3+1$, with $q=3^{2c+1} > 3$.}
\item \textit{$N = PSp(a,2)$, of degree $2^{2l-1} \pm 2^{l-1}$, with $a=2l$ and $l \geq 3$.}
\item \textit{$N = PSL(2,11)$, of degree $11$; $N = PSL(2,8)$, of degree $28$; $N = A_7$, of
degree $15$; $N = HS$, of degree $176$; $N = Co_3$, of degree $276$.}
\item \textit{$N = M_e$, of degree $e=11,12,22,23,24$ or $N = M_{11}$ of degree $12$.}
\enum }
\item[$2.$] $N = \Z_p^R$ for some prime $p$ and integer $d \geq 1$ and $G \cong N \rtimes
G_0$ is an affine permutation group of degree $p^R$; $G_0 \leq \Gamma L(a,q)$ where $q^a=p^R$ and
$G_0$ is one of the following $($note that the symbol {\rm``$\circ$''} denotes a central
product$)$$:$ {\rm \bnum
\item \textit{$a=1$ and $G_0 \leq \Gamma L(1,q)$,}
\item \textit{$a \geq 2$ and $\Gamma L(a,q) \geq G_0 \normgeq SL(a,q)$,}
\item \textit{$a=2l \geq 4$ and $\Z_{q-1} \circ \Gamma Sp(a,q) \geq G_0 \normgeq Sp(a,q)$,}
\item \textit{$a=6$, $q$ even and $\Z_{q-1} \times Aut(G_2(q)) \geq G_0 \normgeq G_2(q)'$,}
\item \textit{$a=4$, $q=2$ and $G_0 \cong A_6$ or $A_7$,}
\item \textit{$a=6$, $q=3$ and $G_0 =SL(2,13)$,}
\item \textit{$a=2$, $q=p=5,7,11$ or $23$ and $G_0 \normgeq SL(2,3)$, or
$a=2$, $q=9,11,19,29$ or $59$ and $G_0 \normgeq SL(2,5)$,}
\item \textit{$a=4$, $q=3$ and $G_0$ has a normal extraspecial subgroup $\mathbf{E}$ of order $2^5$,
and $G_0/\mathbf{E}$ is isomorphic to a subgroup of $S_5$.} \enum } \enum \et

\brk \label{2_trxclassrk} \rm{For each of the cases in Theorem~\ref{2_trxclass} (2), $q$ is chosen
such that $G$ preserves on $V$ the structure of an $a$-dimensional vector space over the finite
field $\F{q}$. Thus $G_0 \leq \Gamma L(a,q)$ and $V = \F{q}^a = V(a,q)$ where $q^a =p^R$ with $a
\geq 1$. Also, $G_0$ (stabiliser of the zero vector in $V$) is transitive on the set of nonzero
vectors in $V$, denoted as $V^*$.} \erk

The next result characterises finite permutation groups that are $2$-homogeneous but not
$2$-transitive.

\bt \label{2hom} {\rm \cite[pp. 368-369]{Huppert3} \mbox{ {\rm or} } \cite{Kantor69}}. \ Suppose
$G$ is a finite $2$-homogeneous but not $2$-transitive permutation group. Then $G$ is isomorphic
to an affine primitive\footnote{See Section~\ref{c2:primitive} for definition of a primitive
permutation group.} group of semilinear mappings $x \mapp x^{\sigma}b+a$ on $\F{q}$ where $q$ is a
prime power, $b \neq 0$, $a, b\in \F{q}$ and $\sigma \in Aut(\F{q})$, that is, $G \leq A\Gamma
L(1,q)$. In particular, $|G|$ is odd and $q \equiv 3$ $($mod $4$$)$. In fact, a permutation group
$G$ is $2$-homogeneous but not $2$-transitive if and only if $G$ is transitive on $\F{q}$ and the
point stabiliser $G_0$ of the zero element has just $2$-orbits, $X$ and $-X$, in $\F{q} \setminus
\{\0\}$, of equal odd length. \et

\section{Blocks and primitivity} \label{c2:primitive}
A \textit{partition} of a finite set $V$ is a set $\M{B}$ of subsets of $V$ such that $\cup_{B \in
\M{B}} B = V$ and $B_1 \cap B_2 =\emptyset$ for distinct $B_1, B_2 \in \M{B}$. Let $G$ be a group
acting on a finite set $V$. A nonempty subset $B \subseteq V$ is called a \textit{block} for $G$
if for every $g \in G$, either $B \cap B^g = \emptyset$ or $B = B^g$. A block $B$ is said to be
\textit{trivial} if $|B| = 1$ or $B = V$. Otherwise, $B$ is called \textit{nontrivial}. Let
$\M{B}$ be a partition of $V$ and $G$ be a group acting on $V$. If $B^g \in \M{B}$ for any $B \in
\M{B}$ and $g \in G$, then we say that $\M{B}$ is a \textit{$G$-invariant partition} of $V$. It is
easy to see that elements of a $G$-invariant partition $\M{B}$ are blocks for $G$. Moreover, $G$
permutes the elements of $\M{B}$ blockwise and induce a (possibly unfaithful) natural action on
$\M{B}$. We say that a group $G$ acting on $V$ is \textit{primitive} if $G$ is transitive on $V$
and the only blocks for $G$ are the trivial ones. If $G$ is transitive but not primitive on $V$,
then $G$ is said to be \textit{imprimitive}.

The following lemma shows that partitions invariant under a transitive permutation group often
arise as orbits set of a normal subgroup.

\bl \label{lem:DD} Let $G$ be a transitive permutation group on $V$ and $M$ a normal subgroup of
$G$. Then $G$ acts transitively on the set $\M{B}$ of $M$-orbits in $V$. In particular, $\M{B}$ is
a $G$-invariant partition of $V$ and all $M$-orbits in $V$ have equal length. \el \bpf We first
prove that $G$ leaves the set $\M{B}$ of $M$-orbits in $V$ invariant, that is, if $\Or \in \M{B}$,
then $\Or^g  \in \M{B}$ for all $g \in G$. Observe that for all $m \in M$ and $g \in G$, $gmg^{-1}
\in M$ (since $M$ is normal in $G$). Hence $\Or^{gmg^{-1}} =\Or$, and so $(\Or^g)^m = \Or^g$. Thus
$\Or^g$ is $M$-invariant. Suppose $u^g, v^g \in \Or^g$. Then we have $u, v \in \Or$. Since $\Or$
is an $M$-orbit, there is an $m \in M$ such that $u^m =v$. Now $(u^g)^{g^{-1} m g} = u^{mg} =v^g$,
and so it follows (since $g^{-1} m g \in M$) that $M$ is transitive on $\Or^g$. Hence $\Or^g  \in
\M{B}$ and $G$ leaves $\M{B}$ invariant. Further, as elements of $\M{B}$ are $M$-orbits in $V$,
$\M{B}$ is a $G$-invariant of $V$.

Finally, let $\Or_i, \Or_j \in \M{B}$ such that $u \in \Or_i$ and $v \in \Or_i$. Then as $G$ is
transitive on $V$, there exists $g \in G$ such that $u^g=v$. Since $G$ leaves the set $\M{B}$
invariant, it follows that $\Or_i^g =\Or_j$. Thus $G$ is transitive on the set $\M{B}$ of
$M$-orbits. The assertion on all $M$-orbits having equal length then follows immediately. \epf

\newpage

\section{Various types of primitive permutation groups}

The structure of finite primitive permutation groups up to permutational isomorphism are described
by the O'Nan-Scott Theorem (for example see \cite{Praeger97b}; see also \cite{LPS88} or
\cite{Scott80}). The O'Nan-Scott Theorem essentially provides an identification of several types
of finite primitive permutation groups such that, for each type, information is given about the
abstract group structure, or the group action, or both. Here, we will briefly describe only three
types of finite primitive permutation groups (we refer readers to \cite{Praeger97b} for further
details about the remaining types).

A finite primitive permutation group $G$ on $V$ is of type \textbf{HA} (holomorph of an abelian
group) if $G = T \rtimes G_0$ is a subgroup of an affine group $AGL(R,p)$ on $V$, where $T \cong
\Z_p^R$ is the (regular) group of translations (and we may identify $V$ with a $R$-dimensional
vector space over $\F{p}$) and $G_0$ is an irreducible subgroup of $GL(R,p)$. We often call a
primitive group of this type an affine primitive permutation group. (Note that the
\textit{holomorph} Hol$(K)$ of a group $K$ is the semidirect product $K \rtimes Aut(K)$, where
$Aut(K)$ acts naturally on the normal subgroup $K$. The group $AGL(R,p)$ is Hol$(T)$ where $T
\cong \Z_p^R$, so the affine primitive groups are primitive subgroups of holomorphs of elementary
abelian groups.)

A primitive permutation group $G$ is of type \textbf{AS} (almost simple) is such that $G$ is an
almost simple group.

Finally, to describe the next type of primitive group, we define a special type of group action
which we will also often encounter throughout this thesis. Let $H$ be a group, $b >1$ a positive
integer and $K$ be a subgroup of the symmetric group $S_b$. Then the wreath product $H \ wr \ K$
is the semidirect product $H^b \rtimes K$ where elements of $K$ act on $H^b$ by permuting the
``entries" of elements of $H^b$, that is, $(h_1, h_2, \ldots, h_b)^{k^{-1}} = (h_{1^k}, h_{2^k},
\ldots, h_{b^k})$ for all $(h_1, h_2, \ldots, h_b) \in H^b$ and $k \in K$. Now suppose $H \leq$
Sym$(\Delta)$. Then the \textit{product action} of $H \ wr \ K$ on $\Delta^b$ is defined as
follows. Elements of $H^b$ act coordinate-wise on $\Delta^b$ and elements of $K$ permute the
coordinates: for $(h_1, \ldots, h_b) \in H^b$, $k \in K$, and $(\delta_1, \ldots, \delta_b) \in
\Delta^b$, \be (\delta_1, \ldots, \delta_b)^{(h_1, \ldots, h_b)} & := &
(\delta_1^{h_1}, \ldots, \delta_b^{h_b}) \\
(\delta_1, \ldots, \delta_b)^{k^{-1}} & := & (\delta_{1^k}, \ldots, \delta_{b^k}) \ee

A primitive permutation group $G$ on $V$ is of type \textbf{PA} (product action) if $V=\Delta^b$
($b > 1$) and $N^b \leq G \leq H \ wr \ S_b \leq$ Sym$(\Delta) \ wr \ S_b$ in its product action,
where $H$ is a primitive permutation group on $\Delta$ of type \textbf{AS} with socle $N$. (The
\textit{socle} of a group $H$, usually denoted as soc$(H)$, is the product of its minimal normal
subgroups.)

We further note that all finite $2$-homogeneous permutation groups, which also include all finite
2-transitive permutation groups, are primitive. (This fact is standard and can be found in texts
like \cite{Cameron_b99, DM_b96}; see also Theorems~\ref{2_trxclass} and \ref{2hom}.) It follows
from an old result of Burnside \cite[Section 154]{Burnside_b55} that the $2$-transitive
permutation groups found in Theorem~\ref{2_trxclass} are of either (1) the almost simple type
\textbf{AS}, or (2) the affine type \textbf{HA} respectively.

For other terminology and notations on permutation groups not defined here, we refer readers to
standard texts such as \cite{Cameron_b99, DM_b96}.

\section{Orbitals for transitive permutation groups}
Let $G$ be a transitive permutation group on $V$. Then $G$ acts faithfully on the set $V \times V$
via the action as in Section~\ref{act}, namely, $(v_1, v_2)^g =(v_1^g, v_2^g)$ where $g \in G$.
The orbits of $G$ on $V \times V$ are known as the \textit{orbitals} of $G$ in $V$ (or simply the
\textit{$G$-orbitals}). The subset $\{(v,v) \mid v \in V \}$ is a $G$-orbital and is called a
\textit{diagonal orbital} or \textit{trivial orbital} of $G$ in $V$. Furthermore if $\M{O} = (u,
v)^G$ is an orbital of $G$ in $V$, then the set $\M{O}^* = (v, u)^G$ is also a $G$-orbital, and is
called the \textit{paired orbital} of $\M{O}$. If $\M{O} =\M{O}^*$, then $\M{O}$ is called a
\textit{self-paired} orbital of $G$ in $V$.

The following is an observation (similar to Lemma~\ref{lem:DD}) about the set of
nontrivial $M$-orbitals, where $M$ is a transitive normal subgroup of a transitive
permutation group.

\bl \label{gin} Let $G$ be a transitive permutation group on $V$. Let $M$ be a
transitive normal subgroup of $G$. Then $G$ leaves invariant the set of nontrivial
$M$-orbitals and the set of self-paired $M$-orbitals in $V$. Furthermore, if $G$ is
$2$-transitive on $V$, then $G$ acts transitively on the set of nontrivial $M$-orbitals
in $V$. \el \bpf Let $\M{O}$ be a nontrivial $M$-orbital in $V$ and let $g \in G$. The
proof that $\M{O}^g = \{(u^g,v^g) \mid (u,v) \in \M{O} \}$ is an $M$-orbital is similar
to the proof of Lemma~\ref{lem:DD} (and so we do not give details here). Suppose now
that $\M{O} =\M{O}^*$ and let $(u,v) \in \M{O}$. Then $(v,u) \in \M{O}$ and hence $(u^g,
v^g), (v^g,u^g) \in \M{O}^g$. Thus $\M{O}^g$ is a self-paired $M$-orbital and so the
first part of the lemma is proved.

Suppose now that $G$ is $2$-transitive on $V$. Let $\Or_i$ and $\Or_j$ be nontrivial
$M$-orbitals in $V$ such that $(u,v) \in \Or_i$ and $(x,y) \in \Or_j$. Since $G$ is
2-transitive, there exists $g \in G$ such that $(u,v)^g=(x,y)$. As $G$ leaves the set of
nontrivial $M$-orbitals invariant, it follows that $\Or_i^g = \Or_j$. Thus $G$ is
transitive on the set of nontrivial $M$-orbitals in $V$. \epf

\section{Graphs: Basic concepts} \label{c2:graphs}
All \textit{graphs} considered in this thesis are \textit{finite, undirected} and \textit{without
loops} or \textit{multiple edges}.

By an (undirected) graph, we mean a pair $\Gamma =(V, E)$ such that $V$ is a set whose element are
called \textit{vertices} of $\Gamma$, and $E$ is a subset of the set $V \times V \setminus
\{\{v,v\} \mid v \in V \}$ whose elements are called \textit{edges} of $\Gamma$. Thus a typical
element or edge in $E$ is a distinct \textit{unordered pair} $\{u,v\}$ where $u,v \in V$. We say
that $u$ is \textit{adjacent} to $v$ in $\Gamma$ if $\{u,v\} \in E$.

By an \textit{arc} in an undirected graph $\Gamma$, we mean an ordered pair $(u,v)$ (or $(v,u)$)
such that $\{u,v\} \in E\Gamma$. Also in an (undirected) graph $\Gamma=(V,E)$, each $\{u,v\} \in E$
always gives rise to two corresponding ordered arc-pairs, $(u,v)$ and $(v,u)$. Thus if the arcs
set of a graph $\Gamma$ always comes in pairs of the form $\{(u,v), (v,u)\}$, then $\Gamma$ is an
undirected graph. We shall use $A\Gamma$ to denote the set of all arcs in $\Gamma$. We are now
ready to present other graph-related definitions and notations.

For $v \in V\Gamma$, we use $\Gamma (v)$ to denote the set of vertices of $\Gamma$
adjacent to $v$. The \textit{valency} of $v$ in $\Gamma$ is defined to be the size of
$\Gamma (v)$. If all vertices of $\Gamma$ has the same valency, then we say that the
graph $\Gamma$ is \textit{regular}.

A \textit{subgraph} of $\Gamma$ is a graph $\Sigma=(V\Sigma, E\Sigma)$ such that $V\Sigma \subseteq
V\Gamma$ and $E\Sigma \subseteq E\Gamma$. We say that a subgraph $\Sigma$ of $\Gamma$ is a
\textit{spanning subgraph} if $V\Sigma = V\Gamma$. Often, such a graph $\Sigma$ is called a
\textit{factor} of $\Gamma$. Two graphs $\Gamma_1$ and $\Gamma_2$ are said to be
\textit{vertex-disjoint} if they have no vertex in common. Similarly two graphs $\Gamma_1$ and
$\Gamma_2$ are said to be \textit{edge-disjoint} if they have no edges in common. However, any two
edge-disjoint graphs may share common vertices. Suppose $\Gamma_1$ and $\Gamma_2$ are two (not
necessarily vertex or edge-disjoint) graphs. Then the \textit{union} $\Gamma_1 \cup \Gamma_2$ of
$\Gamma_1$ and $\Gamma_2$ is the graph having vertex set $V\Gamma_1 \cup V\Gamma_2$ and edge set
$E\Gamma_1 \cup E\Gamma_2$. The union of finitely many graphs is defined similarly.

A \textit{path} of a graph $\Gamma$ of length $h$ is a sequence of $h+1$ vertices $v_1,
v_2, \ldots, v_{h+1}$ such that $v_{i}$ and $v_{i+1}$ are adjacent for $i =1, 2, \ldots,
h$. Such a path is said to \textit{connect} $v_1$ and $v_{h+1}$. A graph $\Gamma$ is
\textit{connected} if for any two distinct vertices $v_1$ and $v_2$ in $\Gamma$, there is
a path which connects $v_1$ and $v_2$. Otherwise, $\Gamma$ is said to be
\textit{disconnected}, and in this case $\Gamma$ is a union of some connected
vertex-disjoint subgraphs which are called the \textit{connected components} of $\Gamma$.

\section{Automorphisms of graphs}
Let $\Gamma_1$ and $\Gamma_2$ be two graphs. A mapping $\phi: V\Gamma_1 \mapp V\Gamma_2$ is called
a (graph) \textit{homomorphism} from $\Gamma_1$ to $\Gamma_2$ if $\phi$ preserves edges, that is,
for any edge $\{u,v\}$ of $\Gamma_1$, $\{\phi(u), \phi(v)\}$ is an edge of $\Gamma_2$. If $\phi$
is a bijection with $\phi^{-1}$ a (graph) homomorphism from $\Gamma_2$ to $\Gamma_1$, then $\phi$
is said to be a (graph) \textit{isomorphism} from $\Gamma_1$ to $\Gamma_2$. In this case, we write
$\Gamma_1 \cong \Gamma_2$. An isomorphism from $\Gamma$ to itself is called an
\textit{automorphism} of $\Gamma$.

An automorphism of $\Gamma =(V\Gamma,E\Gamma)$ can be seen as a permutation of the vertex set
$V\Gamma$ that leaves the edge set $E\Gamma$ invariant. The set of all automorphisms form a
subgroup $Aut(\Gamma)$ of Sym$(V)$ called the \textit{automorphism group} of $\Gamma$. Let $G$ be
a subgroup of $Aut(\Gamma)$. Then we say that a graph $\Gamma$ is $G$-\textit{vertex-transitive},
$G$-\textit{edge-transitive} or $G$-\textit{arc-transitive} if $G$ acts transitively on the
vertices, edges or arcs of $\Gamma$ respectively. Note that a $G$-arc-transitive graph is also
$G$-vertex and $G$-edge-transitive. However the converse is not true in general. There are graphs
that are $G$-vertex and $G$-edge-transitive but not $G$-arc-transitive (see example below).

\bexp \rm{Let $\Gamma =C_3$ (a cycle on $3$ vertices) with vertex set $V\Gamma = \{1,2,3
\}$ and edge set $E\Gamma=\{\{1,2\}, \{2,3\}, \{1,3\} \}$. Let $G =\brac{(1 2 3)} \cong
\Z_3$. Then it is easily seen that $G$ acts transitively on both $V\Gamma$ and
$E\Gamma$. However, $G$ is not arc-transitive on $\Gamma$ as there is no element in $G$
that maps the arc $(1,2)$ to the arc $(1,3)$.}\eexp

The following is a useful characterisation of an arc-transitive graph which can be taken
from standard graph theory texts such as \cite{Biggs_b93} or \cite{GR}.

\bl \label{arcgraph} A connected graph $\Gamma$ is $G$-arc-transitive if and only if $G$ is
transitive on $V\Gamma$ and for any $v \in V\Gamma$, $G_{v}$ is transitive on $\Gamma(v)$. \el


Let $G$ be a transitive permutation group on a finite set $V$. Let $\M{O}$ be a nontrivial
self-paired $G$-orbital. Then an \textit{orbital graph} of $G$ is a graph $\Gamma$ with vertex set
$V$ and edge set $E =\{ \{x,y\} \mid (x,y) \in \M{O}\}$. In the next result, we shall see that an
orbital graph is always arc-transitive, and that the edge set of a $G$-edge and
$G$-vertex-transitive graph can be characterised in terms of its corresponding nontrivial
$G$-orbital.

\bl \label{lem:orbgraph} Let $G$ be a transitive permutation group on $V$. \bnum
\item $G$ is an arc-transitive group of automorphisms of a connected graph $\Gamma=(V,E)$
if and only if $\Gamma$ is an orbital graph for $G$ where $E =\{ \{x,y\} \mid (x,y) \in
\M{O} \}$ and $\M{O}$ is a nontrivial self-paired $G$-orbital in $V$. Furthermore, the
arc set of $\Gamma$ is $\M{O}$.
\item $G$ is an edge-transitive but not an arc-transitive group of automorphisms of a connected
graph $\Gamma=(V, E)$ if and only if the edge set $E =\{ \{x,y\} \mid (x,y) \in \Or \cup
\Or^* \}$, where $\Or$ is a nontrivial $G$-orbital in $V$ with $\Or^*$ as its paired
orbital and $\Or \neq \Or^*$. Furthermore, the arc set of $\Gamma$ is $\M{O} \cup
\M{O}^*$. \enum \el \bpf (1) \ This is well-known, see for example \cite[Theorem 2.1
(b)]{Praeger97b}.

(2) \ See for example \cite[7.53, p. 59]{Tutte}. \epf

A \textit{complete graph} is a graph in which every two distinct vertices are adjacent. We use
$K_n$ to denote a complete graph with $n$ vertices. Given a graph $\Gamma=(V\Gamma, E\Gamma)$, its
\textit{complement}, usually denoted as $\oline{\Gamma}$, is the graph with the same vertex set as
$\Gamma$ in which $u,v \in V$ are adjacent if and only if they are not adjacent in $\Gamma$. Thus
$\Gamma$ and $\oline{\Gamma}$ are edge-disjoint and it is easy to see that $\Gamma \cup
\oline{\Gamma} \cong K_{|V\Gamma|}$. If a graph $\Gamma$ is isomorphic to its complement
$\oline{\Gamma}$, that is $\Gamma \cong \oline{\Gamma}$, then we say that $\Gamma$ is a
\textit{self-complementary} graph.

\section{Cayley graphs}

We now introduce some classes of vertex-transitive graphs that we shall often encounter
throughout this thesis.

\bdeff \label{def:cayley} \rm{\textbf{(Cayley graph)} \ Let $G$ be a finite group and $S=S^{-1}$ be
a nonempty subset of $G$ such that $1_G \notin S$. The \textit{Cayley graph} $\Gamma =Cay(G,S)$ of
$G$ relative to $S$ is the graph with vertex set $G$ such that $\{u,v\}$ is an edge if and only if
$vu^{-1} \in S$. Finally, the set $S$ is called the \textit{connection set} of $Cay(G,S)$. (Note
that for $Cay(G,S)$ to be undirected, we must have $S= S^{-1}$.) } \edeff

The following are some standard results regarding Cayley graphs.

\bl \label{some} Let $\Gamma =Cay(G,S)$ be the Cayley graph of $G$ relative to $S$.
Then$:$ \bnum
\item $\Gamma$ is $\widehat{G}$-vertex-transitive\footnote{Note that $\widehat{G}$ acts
on $V\Gamma=G$ by right multiplication; see also definition of $\widehat{G}$ in
Section~\ref{c2:perm_gp}}, and so $\widehat{G} \leq Aut(\Gamma)$.
\item $Aut(G;S) :=\{\sigma \in Aut(G) \mid S^{\sigma}=S\}$ is a subgroup of
$Aut(\Gamma)$.
\item $\Gamma$ is connected if and only if $\brac{S} = G$,
where $\brac{S}$ is the group generated by elements of $S$. \enum \el \bpf (1) This is
well-known, for example see proof of \cite[Proposition 16.2 (1)] {Biggs_b93}. We shall
present the proof here. Take any $g \in G$, then there is a corresponding $\hat{g} \in
\widehat{G}$ such that $v^{\hat{g}} = vg$ for all $v \in G$. Now notice that
\be \{u,v\} \in E\Gamma & \ifff & vu^{-1} \in S \\
& \ifff & vg (ug)^{-1} \in S \\
& \ifff & \{ug,vg\} \in  E\Gamma \\
& \ifff & \{u^{\hat{g}},v^{\hat{g}}\} \in E\Gamma. \ee Thus $\widehat{G} \leq
Aut(\Gamma)$, and since $\widehat{G}$ acts regularly on $G$, it follows that $\Gamma$ is
$\widehat{G}$-vertex-transitive.

(2) Suppose $\sigma \in Aut(G;S)$. Then by definition, $\sigma \in Aut(G)$. Now $\{u,v\}
\in E\Gamma \ifff vu^{-1} \in S \ifff (vu^{-1})^{\sigma} \in S \ifff v^{\sigma}
(u^{-1})^{\sigma} \in S \ifff v^{\sigma} (u^{\sigma})^{-1} \in S \ifff
\{u^{\sigma},v^{\sigma}\} \in E\Gamma$. Thus $\sigma \in Aut(\Gamma)$.

(3) Suppose $S$ generates $G$. Then each element $g \in G$ can be written as $g =s_m
\cdots s_1$ where $s_i \in S$ for $i \in \{1,\ldots,m\}$. It is easy to see that $1_G, \
s_1, \ s_2s_1, \ s_3 s_2 s_1, \ldots, s_m \cdots s_1 =g$ is a path from $1_G$ to $g$.
Thus every element of $G$ is connected to $1_G$ with a path, and hence $\Gamma$ is
connected.

Conversely, suppose $\Gamma =Cay(G,S)$ is connected. Then for each $g \in G$, there is a
path of length $m$, $1_G=g_0, \ g_1, \ g_2, \ldots, g_m=g$, connecting $1_G$ to $g$. We
claim that each element in this path belongs to $\brac{S}$, that is, $g_i \in \brac{S}$
for $i \in \{0, 1, \ldots, m\}$. We prove this by induction on $i$. For $i=0$, the
element $1_G$ is clearly in $\brac{S}$. Suppose the claim holds for some $i \in \{1,
\ldots, m\}$. Then we show that it also holds for $i+1$. Now $\{g_i, g_{i+1}\}$ is an
edge in $\Gamma$ and so $g_{i+1} g_{i}^{-1} = s$ for some $s \in S$. Since $s$ and $g_i$
is contained in $\brac{S}$, it follows that $g_{i+1} \in  \brac{S}$. Thus by the
induction process, we have $g_i \in \brac{S}$ for all $i \in \{0, 1, \ldots, m\}$; and
as the claim holds for each $g \in G$, we have $\brac{S}=G$. \epf

Furthermore, since $\Gamma = Cay(G,S)$ is vertex-transitive (Lemma~\ref{some} (1)), it follows
that if $\Gamma$ is disconnected, then each connected component is isomorphic to $Cay(\brac{S},S)$
(see for example \cite[Section 1]{survey}). By observing that the vertex set of each connected
component is a (right) coset of $\brac{S}$, it is easily seen that the number of connected
components equals $\frac{|G|}{|\brac{S}|}$.

Notice that by Lemma~\ref{some} (1), a Cayley graph $\Gamma= Cay(G,S)$ contains
$\widehat{G} \leq Aut(\Gamma)$ acting regularly on its vertices. It follows that we can
characterise a Cayley graph through the existence of a subgroup of its full automorphism
group which acts regularly on its vertices.

\bl {\rm \cite[Lemma 16.3]{Biggs_b93}}. \label{lem:cayreg} A graph $\Gamma$ is isomorphic to a
Cayley graph for some group if and only if its automorphism group $Aut(\Gamma)$ has a subgroup
which is regular on vertices. \el

In subsequent chapters, we will often come across edge-transitive and arc-transitive Cayley graphs
$\Gamma = Cay(V,S)$ of a vector space $V=V(a,q)$ (where $a \geq 1$ and $q$ a prime power), which
admit an affine transitive permutation group $G =T \rtimes G_0$ as a subgroup of automorphisms.
(Note that the translation group $T$ acts regularly on $V$.) In the next result, we characterise
(in a similar fashion to Lemma~\ref{lem:orbgraph}) such a Cayley graph with respect to the
$G_0$-orbits in its connecting set $S$.

\bl \label{fact1} Let $G =T \rtimes G_0$ be an affine transitive permutation group on the vector
space $V =V(a,q)$ where $a \geq 1$ and $q$ is a prime power. Let $\Gamma = Cay(V,S)$ be a Cayley
graph on $V$ $($with $S =-S$$)$ and suppose $G$ is a subgroup of $Aut(\Gamma)$. Then the following
hold. \bnum
\item[$1.$] $\Gamma$ is $G$-arc-transitive if and only if the point stabiliser $G_{0}$ $($of the
zero-vector$)$ is transitive on $S$.
\item[$2.$] $\Gamma$ is $G$-edge-transitive but not $G$-arc-transitive if and only if
$G_{0}$ has $2$ orbits on $S = D \cup \widehat{D}$, where $D=s^{G_0}$ for some $s \in S$, and
$\widehat{D}=-D$. In particular, in this case $|V|$ is odd. \enum \el \bpf (1) Observe that the
connecting set $S$ of a Cayley graph $\Gamma =Cay(V,S)$ is precisely the set of all vertices in
$\Gamma$ that are adjacent to the zero-vector $\0 \in V\Gamma =V$. The result then follows from
Lemma~\ref{arcgraph}.

(2) Recall from Lemma~\ref{lem:orbgraph} (2) that $\Gamma =Cay(V,S)$ is $G$-edge-transitive but
not $G$-arc-transitive if and only if the edge set $E\Gamma =\{\{x,y\} \mid (x,y) \in \Or \cup
\Or^*\}$ for some nontrivial $G$-orbital $\Or$ in $V$ with corresponding paired orbital $\Or^*$
and $\Or \neq \Or^*$. Furthermore, $E\Gamma =\{\{x,y\} \mid (x,y) \in \Or \cup \Or^*\} \ifff
A\Gamma =\Or \cup \Or^*$, where $A\Gamma$ is the arc set of $\Gamma$. Now there is a $1$--$1$
correspondence between the set of nontrivial $G$-orbitals and the set of $G_0$-orbits in $V^* :=V
\setminus \{\0\}$ (since for a nontrivial orbital $\Or$, the corresponding $G_0$-orbit is
$\Or(\0):= \{v \in V \mid (\0,v) \in \Or\}$). Thus it follows that

\be \Gamma \mbox{ is $G$-edge-transitive
but not $G$-arc-transitive} & \ifff & E\Gamma =\{\{x,y\} \mid (x,y) \in \Or \cup \Or^*\}, \\
& & \mbox{for some nontrivial orbital} \\ & & \Or \neq \Or^* \\
& \ifff & A\Gamma =\Or \cup \Or^*, \  \Or \neq \Or^*  \\
& \ifff & \mbox{$G_0$ has 2 paired orbits $\Or(\0)$} \\ & & \mbox{and $\Or^*(\0)$ in
$S$.} \ee

Now we need to show that $\Or(\0)^*= - \Or(\0)$. Let $D =\Or(\0)$. Since $\Or(\0)$ is a
nontrivial orbit of $G_0$ in $V^*$, we have $D = s^{G_0}$ where $\0 \neq s \in \Or(\0)$.
Similarly, let $\widehat{D} =\Or^*(\0)$. Then $\widehat{D} = s'^{G_0}$ for some nonzero
$s' \in \Or^*(\0)$. Clearly, $S=D \cup \widehat{D}$ and $D \neq \widehat{D}$. Since $s,
s' \in S$, it follows that $\{\0,s\}$ and $\{\0,s'\}$ are edges of $\Gamma$. As $\Gamma$
is $G$-edge-transitive, there exists $g \in G$ such that $\{\0,s\}^g = \{\0,s'\}$. Now
$g =t_v \mu$ where $t_v \in T$ for some $v \in V$ and $\mu \in G_0$. It follows that
$\{\0,s\}^g = \{\0,s\}^{t_v \mu} = \{v,s+v\}^{\mu} =\{v^{\mu}, (s+v)^{\mu}\} =
\{\0,s'\}$. If $v^{\mu} =\0$, then $v= \0$ and $s'=s^{\mu}$, which is a contradiction
since $s' \notin D$. Hence $(s+v)^{\mu} = \0$, and $v^{\mu} =s'$. Thus $\0 =s^{\mu} +
v^{\mu} = s^{\mu} +s'$, and so $s' = -(s^{\mu}) \in -D$ and thus $\widehat{D} \subseteq
-D$. Furthermore since $-D =(-s)^{G_0}$ is a $G_0$-orbit in $S$ and $\widehat{D} \neq
D$, it follows that $D \neq -D$ and so $\widehat{D} = -D$.

Conversely, suppose $S = D \cup \widehat{D}$ with $D=s^{G_0} =\Or(\0)$ for some $s \in S$, and
$\widehat{D}=-D$. Then a similar argument shows that $-D =\Or^*(\0)$, and with that, the first
part of Lemma~\ref{fact1} (2) follows.

Finally, suppose $\Gamma$ is $G$-edge-transitive but not $G$-arc-transitive and $|V|$ is
even. Then $V=V(a,q)=\F{q}^a$ has characteristic 2, that is, every element of $V$ is its
own (additive) inverse. It follows that $D = s^{G_0} =(-s)^{G_0} = -(s^{G_0}) =-D$, a
contradiction since $D \neq -D$. Thus $|V|$ is odd. \epf

\bdeff \label{def:Paley} \rm{\textbf{(Paley graph)} \ Let $V=\F{q}$ be a finite field with $q=p^R$
elements such that $q \equiv 1$ (mod $4$). Let $\om$ be a primitive element in $V$ and $S$ the set
of nonzero squares in $V^*=\F{q} \setminus \{\0\}$, so $S=\{\om^2, \om^4, \ldots, \om^{q-1}=1 \}
=-S$. The \textit{Paley graph}, denoted by $Paley(q)$, is defined to be the Cayley graph
$Cay(V,S)$ of $V$ relative to $S$. (Note that we take $V$ as the additive group of $\F{q}$.) }
\edeff

The Paley graphs are a special class of Cayley graphs based on the additive groups of
finite fields. They were first defined by Paley in \cite{Paley}, and are known to be
self-complementary and arc-transitive. Furthermore, the full automorphism group of a
Paley graph is contained in $A\Gamma L(1,q)$, see for instance \cite{Peisert2001}.

\bdeff \label{def:Hamming} \rm{\textbf{(Hamming graph)} \ The \textit{Hamming graph}
$H(a,b)$ is a connected graph with vertices the $b$-tuples ($b >1$) with entries from a
set $\Delta$ of size $a$. Two vertices are adjacent in $H(a,b)$ if and only if the two
$b$-tuples differ in exactly one component.} \edeff

The full automorphism group of the Hamming graph $H(a,b)$ is the wreath product $S_a \
wr \ S_b$ (of product action type on $\Delta^b$; see Section~\ref{c2:perm_gp} for
details on the product action) \cite[Theorem 9.2.1]{bcn}. If $A$ is a regular subgroup
of $S_a$, then $A^b < S_a \ wr \ S_b$ and $A^b$ acts regularly on the vertex set of
$H(a,b)$. Thus by Lemma~\ref{lem:cayreg}, each $H(a,b)$ is a Cayley graph. If $|\Delta|$
is a prime power $q$, then $\Delta$ can be identified with a finite field $\F{q}$ of
order $q$. In this case, we may choose $A$ to be the additive group of the field $\F{q}$.

Finally, we refer readers to \cite{Biggs_b93} or \cite{GR} for further graph-theoretic
terminology and notations not mentioned in this chapter.
\chapter{Generalised Paley graphs} \label{c6}


The generalised Paley graphs are a special class of Cayley graphs based on the additive
group of a finite field. They are a generalisation of the Paley graphs, first defined by
Paley in 1933 (see \cite{Paley}), which are well-known to be self-complementary and
arc-transitive. They also arise as factors of arc-transitive homogeneous factorisations
$(M,G,V,\M{E})$ of complete graphs where $G$ is a one-dimensional affine group (see
Section~\ref{c3:examples}). In this chapter, we study the structure of generalised Paley
graphs and their automorphism groups. One of the main results here is important for its
use in Chapter~\ref{c5} to distinguish the generalised Paley factors from another class
of factor graphs (the twisted generalised Paley graphs defined in Chapter~\ref{c3}). This
chapter is based on the work done by Praeger and the author in \cite{LP}.

\section{Main results on structure of generalised Paley graphs} \label{c6:main}
We will state the major results of this chapter and will prove them in subsequent
sections. First, we introduce some notations.

Let $V=\F{q}$ be a finite field of order $q=p^R$ where $p$ is a prime and $R \geq 1$. Let $V^*= V
\setminus \{\0\}$ be the set of nonzero elements in $V$. For a fixed primitive element $\om$ in
$V$, we let $\omb$ be the corresponding scalar multiplication $\omb: x \mapp x\om$ for $x \in V$.
Then let $\brac{\omb^i}$ be the multiplicative subgroup of $V^*$ generated by the element $\omb^i$
where $i \geq 1$ and $i \mid (q-1)$. (Note that $\omb^i: x \mapp x\om^i$, $x \in V$.) Since for
each $\omb^i$, there is a corresponding $\om^i \in V$, we shall use $\brac{\om^i}$ to denote the
set $\{1,\om^i,\om^{2i}, \ldots, \om^{((q-1)/i )-i} \} \subseteq V^*$. Note that for $1, -1 \in
V^*$, we will use $\pone$ and $\mone$ to denote the corresponding scalar multiplications. Also, we
let $\al$ denote the Frobenius automorphism of $\F{p^R}$, that is, $\al: x \mapp x^p$. Finally, we
note that the group $\Gamma L(1,p^R)$ is generated by the scalar multiplication map $\omb$ and the
Frobenius automorphism $\al$. Thus $\Gamma L(1,p^R) =\brac{\omb, \al}$ and hence $A\Gamma L(1,p^R)
= T \rtimes \brac{\omb, \al}$, where $T$ is the translation group acting additively (and
regularly) on $V =\F{p^R}$.

\bdeff \label{defn:gpaley:c6} \rm{\textbf{(Generalised Paley graph)} \ Let $V$ and $V^*$
be as defined above. Let $k \geq 2$ be an integer which divides $q-1$ such that if $p$
is odd, then $\frac{q-1}{k}$ is even. Let $\om$ be a fixed primitive element in $V$. The
graph $\GPq$ is the Cayley graph $Cay(V,S)$ with connecting set $S = \brac{\om^k} =\{ 1,
\om^k, \om^{2k}, \ldots, \om^{q -1 -k} \} \varsubsetneq V^*$ and is called a
\textit{generalised Paley graph} with respect to $V$ (that is to say it is the graph
with vertex set $V$ such that $\{u,v\}$ is an edge if and only if $v - u \in S$).} \edeff

\brk \label{rk:gpaley:c6} \rm{Since $\frac{q-1}{k}$ is even when $p$ is odd, we have
$S=-S$, and so $\GPq$ is an undirected Cayley graph (see also
Proposition~\ref{prop:gpaley}). Also, $\GPq$ has valency $\frac{q-1}{k} =|S|$. Observe
that if $k=2$, then $\GPq$ is the familiar Paley graph. In this case, since we require
$\frac{q-1}{2}$ to be even for $q$ odd, it follows that $q \equiv 1$ (mod 4). Finally, a
generalised Paley graph is not necessarily connected. In Theorem~\ref{thm:add} (see also
Lemma~\ref{lem:connect}), we give a necessary and sufficient condition for a $\GPq$ to
be connected.} \erk

The following simple observation tells us that $T \rtimes \brac{\omb^k,\al} < A\Gamma L(1,p^R)$
occurs as a subgroup of automorphisms in a generalised Paley graph. (Note that $k$ is assumed to
satisfy the conditions in Definition~\ref{defn:gpaley:c6}.)

\bl \label{fel} Let $\Gamma =\GP=Cay(V,S)$ be a generalised Paley graph. Then $T \rtimes
\brac{\omb^k, \al} \leq Aut(\Gamma)$. Further, $T \rtimes \brac{\omb^k}$ is arc-transitive on
$\Gamma$. \el \bpf The first part of the lemma follows from Lemma~\ref{some} (1) -- (2) by
identifying the group $\widehat{G}$ acting by right multiplication on $G$ (in Lemma~\ref{some}
(1)) with the translation subgroup $T$ acting (additively) on $V$, and by observing that
$\brac{\omb^k, \al}$ is a subgroup of $Aut(V)$ that leaves $S$ invariant. The second assertion in
the lemma then follows from Lemma~\ref{fact1} (1) by observing that $\brac{\omb^k}$ acts regularly
on $S =\brac{\om^k}$. \epf

\bt \label{thm:add} Let $V=\F{p^R}$ be a finite field of order $p^R$ and $S$ be as
defined above. Suppose $\Gamma =\GP=Cay(V,S)$ is a disconnected generalised Paley graph.
Then $S$ is contained in a proper subfield $\F{p^a}$ of $V$ $($where $a$ is a proper
divisor of $R$$)$ and each connected component of $\Gamma$ is isomorphic to $GPaley(p^a,
\frac{p^a-1}{k'}) = Cay(\F{p^a},S)$, where $k = \frac{p^R-1}{p^a-1} \cdot k'$. \et

The above theorem follows immediately from Lemma~\ref{lem:connect} which we prove in the
next section. Since (by Theorem~\ref{thm:add}) the components of a disconnected
generalised Paley graph are generalised Paley graphs of subfields, we will, in the
remainder of the chapter, only consider connected generalised Paley graphs.

Our next main result below gives a necessary and sufficient condition for a connected generalised
Paley graph to be isomorphic to a Hamming graph.

\bt \label{thm:A} Let $V=\F{p^R}$ be a finite field of order $p^R$. Let $\Gamma$ $=$
$\GP$ $=$ $Cay(V,S)$ be a connected generalised Paley graph. Then $\Gamma \cong H(p^a,
\frac{R}{a})$ $($Hamming graph$)$, for some divisor $a$ of $R$ with $1 \leq a < R$, if
and only if $k=\frac{a(p^R-1)}{R(p^a-1)}$. \et

As we noted in Chapter~\ref{c2}, the automorphism group of a Hamming graph is a wreath product of
symmetric groups. By contrast, we will show in our next result, that whenever a connected
generalised Paley graph is not isomorphic to a Hamming graph, then its full automorphism group is a
primitive affine group (we refer readers to Section~\ref{c2:primitive} for details on primitive
groups). Note that the result below is crucial for the proof of Theorem~\ref{thm:no2}, needed in
Chapter~\ref{c5}.

\bt \label{thm:no1} Let $V=\F{p^R}$ be a finite field of order $p^R$. Let $\Gamma =\GP=Cay(V,S)$
be a connected generalised Paley graph. Suppose $k \neq \frac{a(p^R-1)}{R(p^a-1)}$ for every
divisor $a$ of $R$ with $1 \leq a < R$. Then $Aut(\Gamma)$ is a primitive subgroup of $AGL(R,p)$
containing the translation group $T \cong \Z_p^R$. \et

Theorem~\ref{thm:no1} proves that the graph $\Gamma$ is a \textit{normal Cayley graph}, that is,
the regular subgroup $T$ of automorphisms is normal in the full automorphism group of $\Gamma$.
Its proof uses results from \cite{Praeger90} that rely on the finite simple groups classification.

Although Theorems~\ref{thm:A} and \ref{thm:no1} provide a great deal of information about the
generalised Paley graphs, we do not yet have a complete knowledge of their full automorphism
groups. In particular, for Theorem~\ref{thm:no1}, we do not know the precise conditions under which
$Aut(\Gamma)$ is contained in the one-dimensional affine group $A\Gamma L(1,p^R)$ (much like in
the case of a Paley graph where the full automorphism group is a subgroup of index 2 in $A\Gamma
L(1,p^R)$). We find it worthwhile to pose the following problem.

\bpro \label{prob:p1} \rm{Determine the full automorphism group of the generalised Paley
graph \ $\GP$ for all $k \geq 2$. }\epro

We are able to provide a partial solution to this problem. For the case when $k$ divides
$p-1$ (and also $2k \mid (p^R-1)$ if $p$ is odd), we will prove that $Aut(\Gamma)$ is a
one-dimensional affine group. This result depends heavily on results in \cite{GPPS}, and
hence on the finite simple group classification. This result was applied in
Chapter~\ref{c5} to distinguish between the generalised Paley graphs and their
``twisted" versions (see Chapter~\ref{c3} for definition of a twisted generalised Paley
graph).

\bt \label{thm:no2} Suppose $\Gamma = \GP$ is a generalised Paley graph such that $k \mid (p-1)$,
and if $p$ is odd, then $2k \mid (p^R-1)$. Then $\Gamma$ is connected and $Aut(\Gamma) = T \rtimes
\brac{\omb^k, \al} \leq A\Gamma L(1,p^R)$. \et

\brk \rm{Suppose $R=1$. Then by definition of a generalised Paley graph, $k$ always
divides $p-1$. Thus for $R=1$, Theorem~\ref{thm:no2} always applies.}\erk

We note that apart from the disconnected case and the Hamming graph case already dealt
with in Theorems~\ref{thm:add} and \ref{thm:A}, there are other cases where $Aut(\GP)$
is not a one-dimensional affine group. We give an explicit example below.

\bexp \label{ex:1} \rm{Take $R=4$, $p=3$ and $k =4$. Then $\frac{p^R-1}{k} =20 $ and $k \neq
\frac{a(p^R-1)}{R(p^a-1)}$ for any divisor $a$ of $R$ with $1 \leq a < R$. Let $\Gamma =
GPaley(3^4, 20)$ be a generalised Paley graph with parameters $(R,p,k) = (4,3,4)$ as given above.
By Theorems~\ref{thm:A} and \ref{thm:no1}, $\Gamma$ is not a Hamming graph and $Aut(\Gamma)$ is a
primitive subgroup of $AGL(4,3)$. Using \textsc{Magma} \cite{Magma}, we computed $Aut(\Gamma)$,
which has order $|Aut(\Gamma)| = 233 280$, greater than $|T \rtimes \brac{\omb^4, \al}|=6480$. Thus
$Aut(\Gamma)$ is not a one-dimensional affine group.

A further check using \textsc{Magma} shows that the point-stabiliser $A_0$ of $A:=Aut(\Gamma)$ has
order $2880$. Further, $A_0$ contains a normal subgroup $B$ isomorphic to $A_6$, and by taking the
quotient, we have $A_0/B \cong D_4$, where $D_4$ is the dihedral group (order 8) of degree $4$. So
$A_0$ is an extension of $A_6$ by $D_4$; written as $A_0 = A_6 \cdot D_4$. Thus it follows that
$Aut(\Gamma) = T \rtimes ( A_6 \cdot D_4)$.}\eexp

Thus in spite of the various results above, Problem~\ref{prob:p1} is still largely open.
It would be interesting to determine a necessary and sufficient condition, based on the
parameters $(R,p,k)$, for $Aut(\GP)$ to be a one-dimensional affine group.

\section{Conditions for generalised Paley graphs to be connected} \label{c6:conn}
Let $\Gamma = \GP = Cay(V,S)$. Recall that by Lemma~\ref{fel}, $Aut(\Gamma)$ contains $T
\rtimes W$ as an arc-transitive subgroup, where $W := \brac{\omb^k}$. The following
result (Lemma~\ref{lem:connect}) relates the connectedness of $\Gamma$ with the action of
$W$ on $V$; Theorem~\ref{thm:add} then follows immediately from part (2) of
Lemma~\ref{lem:connect}.

\bl \label{lem:connect} Let $\Gamma =\GP =Cay(V,S)$ be a generalised Paley graph with
vertex set $V =\F{p^R}$ and $S =\brac{\om^k}$. Then the following hold$:$ \bnum
\item[$1.$] $\Gamma$ is connected if and only if $\brac{\omb^k}$ acts irreducibly on
$V$.
\item[$2.$] Suppose $\Gamma$ is not connected. Then the connected components of $\Gamma$ are
isomorphic to $\Gamma_0 := GPaley (p^a, \frac{p^a-1}{k'})$ for some proper divisor $a$ of $R$,
where $k' \mid (p^a-1)$ and $\frac{p^R-1}{k} = \frac{p^a-1}{k'}$. Furthermore, $Aut(\Gamma) =
Aut(\Gamma_0) \ wr \ S_{p^{R-a}}$. \enum \el \bpf (1) \ Let $U$ be the $\Fp$-span of $S$, that is,
$U = \{ \sum_{\om^{ik} \in S} \lambda_i \om^{ik} \mid \lambda_i \in \Fp \}$. Since $W
=\brac{\omb^k}$ leaves $S$ invariant, it also leaves invariant the $\Fp$-span $U$ of $S$. Also, we
note that $\Gamma$ is connected if and only if $U = V$ (which may be seen by regarding $V$ as an
$R$-dimensional vector space over $\F{p}$).

Suppose $W$ acts irreducibly on $V$. Then as $U$ is $W$-invariant and nonzero, we have
$U = V$ and so $\Gamma$ is connected. Conversely suppose $\Gamma$ is connected. Then $S$
is a spanning set for $V$ (that is $U=V$). Now $S$ is the orbit $1^W = \brac{\om^k}$
(note that $W$ acts by field multiplication). Also, for each $\om^i \in V^*$, $\omb^i$
maps $S$ to $S \om^i$, and as $\omb^i \in GL(1,p^R)$, it follows that $S \om^i$ is also
a spanning set for $V$ (as an $R$-dimensional vector space over $\Fp$). However $S
\om^i$ is the $W$-orbit containing $\om^i$. Thus it follows that every $W$-orbit in
$V^*$ is a spanning set for $V$, and hence $W$ is irreducible on $V$.

(2) \ Suppose $\Gamma$ is disconnected. Let $U$ be the vertex set of the connected
component of $\Gamma$ containing $1 \in \F{p^R}$. Then $U$ is the $\Fp$-span of $S$.
Since $\Gamma$ is vertex-transitive, all the connected components of $\Gamma$ are
isomorphic to $Cay(U,S)$. We claim that $U$ is a subfield of $V =\F{p^R}$.

Since $U$ is $W$-invariant, $U^{\omb^{ik}} =U$ for each $\omb^{ik} \in W$, and hence $U
\om^{ik} = U$ for each $\om^{ik} \in S$. Thus $U$ is closed under multiplication by
elements of $S$. Now each element of $U$ is of the form $\sum_{\om^{ik} \in S} \lambda_i
\om^{ik}$ for some $\lambda_i \in \Fp$, and (by regarding $\lambda_i$ as an integer in
the range $0 \leq \lambda_i \leq p-1$) each $\lambda_i \om^{ik}$ is equal to the sum
$\om^{ik} + \cdots + \om^{ik}$ ($\lambda_i$-times). Thus each element of $U$ is a sum of
a finite number of elements of $S$. Since $U$ is closed under addition and under
multiplication by elements of $S$, it follows that $U$ is closed under multiplication.
Thus $U$ is a subring of $V$. Also $U$ contains the identity $1$ of $V$ (since $1 \in
S$). Let $u \in U$. Then since $V$ is finite, we have that $u^j =1$ for some $j$ and,
since $U$ is closed under multiplication, $u^{-1} = u^{j-1} \in U$. Thus $U$ is a
subfield of $V=\F{p^R}$ as claimed.

Since $U$ is a proper subfield of $V$, we have $|U| =p^a$ where $a \mid R$ and $a \neq
R$. Also, since $S \leq U^*$, it follows that $|S| = \frac{p^R-1}{k}$ divides $p^a-1$.
Let $k' = \frac{(p^a-1)k}{p^R-1}$. Then $|S| = \frac{p^R-1}{k} =\frac{p^a-1}{k'}$, and
by definition of a generalised Paley graph, we have $Cay(U,S) = GPaley(p^a,
\frac{p^a-1}{k'})$. Thus the connected components of $\Gamma$ are isomorphic to
$GPaley(p^a, \frac{p^a-1}{k'})$ and since there are $p^{R-a}$ components in $\Gamma$, it
follows that $Aut(\Gamma) = Aut(GPaley(p^a, \frac{p^a-1}{k'})) \ wr \ S_{p^{R-a}}$. \epf

From Lemma~\ref{lem:connect} (2) (see also Theorem~\ref{thm:add}), we see that whenever
a generalised Paley graph $\Gamma = \GP$ is disconnected, the connected components are
generalised Paley graphs for subfields. Thus without loss of generality we will, from
now on, assume that $\Gamma =\GP$ is connected.

\section{Proof of Theorem~\ref{thm:A}} \label{sec:A} In this section, we shall prove
Theorem~\ref{thm:A}.

\vs \noindent \textit{Proof of Theorem~\ref{thm:A}.} \ Suppose $\Gamma = \GP \cong
H(q,b)$. Then the number of vertices is $p^R= q^b$ and so $q=p^a$ where $R =ab$. Now the
valency of $\Gamma$ is $\frac{p^R -1}{k} = b(q-1)$. So $k = \frac{p^R-1}{b(q-1)} =
\frac{p^R-1}{b(p^a-1)} =  \frac{a(p^R-1)}{R(p^a-1)}$ as required.

Conversely suppose that $k = \frac{a(p^R-1)}{R(p^a-1)}$ where $a \mid R$ and $1 \leq a <
R$, and set $R=ab$. Then since $\Gamma = \GP$ is connected, the $\Fp$-span of $S
=\brac{\om^k}$ equals $V$, that is, $\{ \sum_{\om^{ik} \in S} \lambda_i \om^{ik} \mid
\lambda_i \in \Fp \} =V$. Let $U$ be the $\F{p^a}$-span of the set $X := \{1, \om^k,
\om^{2k}, \ldots, \om^{(b-1)k} \}$. We claim that $U=V$.

Observe that $\F{p^a}, X \subseteq V =\F{p^R}$. So it follows that $U \subseteq V$.
Moreover, since $k = \frac{a(p^R-1)}{R(p^a-1)}$, $S =\brac{\om^k}$ has order
$\frac{p^R-1}{k} = \frac{R}{a} \cdot (p^a-1) = b(p^a-1)$, it follows that $\om^{bk}$ has
order $p^a-1$.  Thus $\brac{\om^{bk}} = \F{p^a}^*$. Now suppose $v \neq 0$ and $v \in
V$. Then $v = \sum \lambda_i \om^{ik}$ where $\lambda_i \in \Fp$ (not all zero) and the
sum is over all $\om^{ik} \in S$. Let $i =bx_i +r_i$ where $0 \leq r_i < b$. Then
$\om^{ik} = \om^{(bx_i +r_i)k} = \om^{bx_ik} \cdot \om^{r_ik}$ and we have \be v \ = \
\sum \lambda_i \om^{ik} & = & \sum \lambda_i \om^{bx_ik} \cdot \om^{r_ik}. \ee Note that
$\lambda_i \in \Fp \subseteq \F{p^a}$. Also, since $\F{p^a}^* =\brac{\om^{bk}}$, it
follows that $\om^{bx_ik} \in \F{p^a}$. Thus $\lambda_i \om^{bx_ik} \in \F{p^a}$ and so
$v \in U$. Thus $U=V$ as claimed.

So from now on, we shall regard $V$ as a vector space over $\F{p^a}$. We have shown that
$V$ is spanned by the set $X = \{1, \om^k, \om^{2k}, \ldots, \om^{(b-1)k} \}$, and as
dim$_{\F{p^a}}(\F{p^R}) =b$, it follows that $X$ is an $\F{p^a}$-basis for $V$. Define
$\Theta: V \mapp \F{p^a}^b$ as follows. For $u  = \sum_{j=0}^{b-1} \mu_j \om^{jk}$ with
$\mu_j \in \F{p^a}$, let \[ \Theta(u) = (\mu_0, \mu_1, \ldots, \mu_{(b-1)}) \in
\F{p^a}^b. \] Note that since $X$ is an $\F{p^a}$-basis for $V$, $\Theta$ is a bijection.

Next, we determine the image of the connecting set $S \subset V$ under $\Theta$. As we
observed above, $|S| = |\brac{\om^k}| =b(p^a-1)$. Thus each element of $S$ can be
expressed uniquely as $\om^{ik}$ for some $i$ such that $0 \leq i \leq b(p^a-1)-1$. As
before, we write $i =bx_i +r_i$ where $0 \leq r_i \leq b-1$. So $\om^{ik} =\om^{bx_ik}
\cdot \om^{r_ik}$. Now $\om^{bx_ik} \in \F{p^a}^* =\brac{\om^{bk}}$, and so
$\Theta(\om^{ik}) = \Theta(\om^{bx_ik} \cdot \om^{r_ik}) = (0, \ldots, 0,
\underbrace{\om^{bx_ik}}_{r_i\mbox{th}}, 0, \ldots, 0)$. Observe that $x_i$ can be any
integer satisfying $0 \leq x_i \leq p^a-2$, and that $\om^{bk x_i}$ takes on each of the
values in $\F{p^a}^*$. Moreover each of these values occurs exactly once in each of the
positions $r$, for $0 \leq r \leq b-1$. Thus $\Theta (S)$ is the set of all elements of
$\F{p^a}^b$ with exactly one component non-zero, that is, the set of ``weight-one"
vectors.

Now $\Theta$ determines an isomorphism from $\Gamma$ to the Cayley graph for $\F{p^a}^b$
with connecting set $\Theta(S)$. In this Cayley graph, two $b$-tuples $u, v \in
\F{p^a}^b$ are adjacent if and only if $u-v \in \Theta(S)$, that is, if and only if
$u-v$ has exactly one non-zero component. Thus $\Gamma =\GP$ is mapped under the
isomorphism $\Theta$ to the Hamming graph $H(p^a, b)$ where $b =\frac{R}{a}$. \qedd

\section{Proof of Theorem~\ref{thm:no1}} \label{sec:pfthm1}
Recall that by Lemma~\ref{lem:connect} (1), if $\Gamma=\GP$ is connected, then
$W=\brac{\omb^k}$ acts irreducibly on $V$. It follows that the group $G = T \rtimes W$
is a vertex-primitive subgroup of $Aut(\Gamma)$ of affine type (or of type HA). Thus
$Aut(\Gamma)$ is a primitive permutation group on $V$ containing $G$. In
\cite{Praeger90} the possible types for $Aut(\Gamma)$ were determined. The proofs in
\cite{Praeger90} depend on the finite simple group classification, and therefore
Theorem~\ref{thm:no1} also depends on this classification since we use results from
\cite{Praeger90} and \cite{PS92} in its proof.

\vs \noindent \textit{Proof of Theorem~\ref{thm:no1}.} \ Suppose $\Gamma=\GP =Cay (V,S)$
is connected and let $G = T \rtimes W $ where $T \cong \Z_p^R$ and $W =\brac{\omb^k}$.
Suppose also that $k \neq \frac{a(p^R-1)}{R(p^a-1)}$ for any $a \mid R$ with $1 \leq a <
R$. Note that $G \leq A:=Aut(\Gamma)$, and since $W$ is irreducible on $V$ (by
Lemma~\ref{lem:connect}), it follows that $G$ is a primitive subgroup of $AGL(R,p)$.

Now suppose the group $A$ is not contained in $AGL(R,p)$ (that is $A$ is not of type HA).
Since $k \geq 2$, $\Gamma$ is not a complete graph and so $A \neq S_{p^R}$ or $A_{p^R}$.
By \cite[Proposition 5.1]{Praeger90}, it follows that $A$ is primitive of type PA. Thus
(see Section~\ref{c2:primitive}) $R=ab$ with $b \geq 2$ (and so $1 \leq a <R$), $V =
\Delta^b$ where $|\Delta|=p^a$, and $N^b \leq A \leq H \ wr \ S_b$ where $H$ is a
primitive permutation group on $\Delta$ that is almost simple with socle $N$. Moreover
from \cite[Proposition 5.1]{Praeger90} and \cite[Proposition 2.1]{PS92} we have that $N
= A_{p^a}$ or $N$ is one of the groups listed in \cite[Table 2]{Praeger90} (denoted as
$L_1$ in \cite{Praeger90}). Note that in all cases, $N$ is nonabelian simple with $p^a
\geq 5$, $N^b$ is transitive on $V = \Delta^b$ and $N$ acts 2-transitively on $\Delta$.

We will prove that $\frac{p^R-1}{k} = b(p^{a}-1)$, thus contradicting the assumption on
$k$ and therefore proving the theorem. To do this, we need to know how a point-stabiliser
subgroup of $A$ acts on $V$.

Let $u:=(\underbrace{\gamma, \ldots, \gamma}_{b-tuple}) \in V$ and let $A_{u}$ be its
corresponding point-stabiliser in $A$. It follows that since $N^b$ is transitive on $V$,
we have $A=N^b A_{u}$. Now $A_{u}$ contains $(N^b)_{u} = (N_{\gamma})^b$, and each
$N_{\gamma}$ acts transitively on $\Delta - \{\gamma\}$. Thus if $v \neq u$ and $v :=
(\delta_1, \ldots, \delta_b)$ such that $v$ has $\ell$ entries different from $\gamma$
(where $1 \leq \ell \leq b$), then $|v^{(N^b)_{u}}| = (p^{a}-1)^{\ell}$.

Next, we need to find out how $A_{u}$ permutes the entries of a point in $V$. Let $\pi :$
Sym$(\Delta) \ wr \ S_b \mapp S_b$ be the projection map with $ker(\pi) =$
(Sym$(\Delta))^b$. Since $A$ is a primitive subgroup of Sym$(\Delta) \ wr \ S_b$, $\pi
(A)$ is a transitive subgroup of $S_b$ (for instance, see \cite[Theorem
4.5]{Cameron_b99}). As $A=N^b A_{u}$, it follows that $\pi(A) = \pi(A_{u})$. Thus, the
subgroup of $S_b$ induced by $A_{u}$ permutes the $b$ entries of points of $\Delta^b$
transitively. Consider the action of $A_{u}$ on $v$ where $v = (\delta_1, \ldots,
\delta_b)$ is as defined above. Since $\pi(A_{u})$ permutes the $b$ entries of a point
transitively, the $\ell$-subset of entries $\delta_i$ of $v$ such that $\delta_i \neq
\gamma$ has at least $b/ \ell$ distinct images under $\pi (A_u)$. Let $n_{\ell}$ be the
number of distinct images under $\pi(A_u)$ of this $\ell$-subset. It follows that the
length of the $A_{u}$-orbit containing $v$ is at least $n_{\ell} \cdot (p^{a}-1)^{\ell}
\geq \frac{b}{\ell} \cdot (p^{a}-1)^{\ell}$. Recall that the graph $\Gamma$ has valency
$\frac{p^R-1}{k}$, and choose $v \in \Gamma(u)$. As $A_{u}$ fixes setwise the set of
points adjacent to $u$, we must have \ben \label{Eqa} \frac{p^R-1}{k} \geq
\frac{b}{\ell} \cdot (p^{a}-1)^{\ell}. \een

Observe that $A_{u}$ contains $G_u=W=\brac{\omb^k}$. Now all $W$-orbits in $V \setminus
\{u \}$ have length $\frac{p^R-1}{k}$. It follows that all orbits of $A_{u}$ in $V
\setminus \{ u \}$ have length greater than or equal to $\frac{p^R-1}{k}$. But there
exists an orbit of $A_{u}$ in $V \setminus \{ u \}$ of length $b(p^{a}-1)$ (the set of
$b$-tuples with exactly one entry different from $\gamma$), and so \ben \label{Eqb}
b(p^{a}-1) \geq \frac{p^R-1}{k}. \een Combining inequalities (\ref{Eqa}) and
(\ref{Eqb}), we obtain \ben \label{Eqc} \ell \geq (p^{a}-1)^{\ell - 1}. \een Since $p^a
\geq 5$, it is easily checked that the inequality (\ref{Eqc}) holds if and only if $\ell
=1$, and it follows that $\frac{p^R-1}{k} = b(p^{a}-1)$ as claimed.

However, this implies that $k = \frac{p^R-1}{b(p^a-1)} = \frac{a(p^R-1)}{R(p^a-1)}$,
contradicting the fact that $k \neq \frac{a(p^R-1)}{R(p^a-1)}$. Thus $A$ is a primitive
subgroup of $AGL(R,p)$. \qedd

\section{The case when $k \mid (p-1)$ and proof of Theorem~\ref{thm:no2}} \label{sec:subcase1}
Let $\Gamma = \GP =Cay(V,S)$ be a generalised Paley graph (where $V$ and $S$ are as
defined in Definition~\ref{defn:gpaley}). We first show that if $k$ divides $p-1$, then
$\Gamma$ is connected.

\bl \label{lem:conns} Suppose $k$ divides $p-1$. Then $\Gamma =\GP$ is connected. \el
\bpf Suppose not. Then by Theorem~\ref{thm:add}, $k$ is a multiple of
$\frac{p^R-1}{p^r-1}$ for some proper divisor $r$ of $R$. Let $k =\frac{p^R-1}{p^r-1}
\cdot k'$ for some integer $k'$. Then since $k \mid (p-1)$, it follows that
$\frac{(p-1)(p^r-1)}{k'(p^R-1)} = \frac{p-1}{k'(p^{R-r} + p^{R-2r} + \cdots + p^r +1)}$
is an integer, which is impossible for any proper divisor $r$ of $R$. \epf

Recall that by Lemma~\ref{fel}, $\Gamma$ admits $T \rtimes \brac{\omb^k,\al}$ as an
arc-transitive subgroup of automorphisms. By Lemma~\ref{lem:connect}, $\Gamma$ is
connected if and only if $W =\brac{\omb^k}$ acts irreducibly on $V$. Thus by
Lemma~\ref{lem:conns}, we always have $W$ acting irreducibly on $V$.

In this section we will prove Theorem~\ref{thm:no2}, that is, we will show that when $k$
divides $p-1$, then the full automorphism group of $\Gamma$ is $Aut(\Gamma) = T \rtimes
\brac{W,\al} \leq A\Gamma L(1,p^R)$ where $W =\brac{\omb^k}$.

Observe that when $k=2$, $\Gamma = \GPP{p^R}{2}$ is a Paley graph whose automorphism
group is well known to be $Aut(\Gamma) = T \rtimes \brac{\omb^2, \al} \leq A\Gamma
L(1,p^R)$ (see for example \cite{Peisert2001}). Thus we may assume that $k \geq 3$. Note
also that if $k \geq 3$ and $k \mid (p-1)$, then the prime $p$ must be at least $5$.
Consequently, we shall assume that $k \geq 3$ and $p \geq 5$.

First we apply Theorem~\ref{thm:no1}. Let $A:= Aut(\Gamma)$.

\bl \label{lem:newstuff} For $k \mid (p-1)$,  $A \leq AGL(R,p)$. \el \bpf Suppose this is
not the case. By Theorem~\ref{thm:no1}, $R=ab$ with $b >1$ and $k =
\frac{p^R-1}{b(p^a-1)}$. Now $\frac{p^R-1}{p^a-1} = 1 + p^a + \cdots + (p^{a})^{b-1} >
(p^a)^{b-1}$ and as $ k < p$ we have $b = \frac{p^R-1}{k(p^a-1)} > p^{ab-a-1} \geq
p^{a(b-2)}$. It follows, since $p \geq 5$, that $b=2$, $a=1$ and $k =\frac{p+1}{2}$.
However, since $ k \mid (p-1)$, this is not possible. \epf

Thus by Lemma~\ref{lem:newstuff}, we have $A = T \rtimes A_0$ where $\brac{W, \al} \leq A_0 \leq
GL(R,p)$. Note that $A_0$ preserves $S \subset V$, and hence $A_0$ does not contain $SL(R,p)$. We
identify $V=\F{p^R}$ with an $R$-dimensional vector space over the prime field $\Fp$. The group
$\Gamma L(1,p^R) =\brac{\omb, \al}$ is then identified with a subgroup of $GL(R,p)$ acting on the
vector space $V=V(R,p)$. Moreover the cyclic group $\brac{\omb}$ is called a \textit{Singer cycle}
(acting irreducibly on $V$) and contains $Z:=Z(GL(R,p)) \cong \Z_{p-1}$ (see \cite[Satz
\textbf{II}.7.3 (Theorem 7.3)]{Huppert1}). We next observe that if $A_0$ is a one-dimensional
group then $A_0 =\brac{W, \al}$.

\bl \label{lem:pre} Suppose $A = Aut(\Gamma) = T \rtimes A_0 \leq AGL(R,p)$, and that
$A_0 \leq \Gamma L(1,p^R) = \brac{\omb, \al}$. Then $A = T \rtimes \brac{W, \al}$ where
$W =\brac{\omb^k}$. \el \bpf Since $\brac{W, \al} \leq A_0 \leq \brac{\omb, \al}$, it
follows that $|A_0|= R \cdot |A_0 \cap \brac{\omb}| \geq R \cdot |W|$. Thus in order to
prove the lemma, it suffices to prove that $A_0 \cap \brac{\omb} =W$. Since
$\brac{\omb}$ is semiregular on $V^*$, and since $A_0$ leaves $S=\brac{\om^k}$
invariant, it follows that $A_0 \cap \brac{\omb} \leq \brac{\omb^k}$. As $|A_0| \geq R
\cdot |W|$, we have $A_0 \cap \brac{\omb} =W$ as required. \epf

Thus for the case where $k \mid (p-1)$, if we are able to show that $A_0 \leq \Gamma
L(1,p^R)$, then by Lemma~\ref{lem:pre}, $Aut(\Gamma) = T \rtimes \brac{W, \al}$. We
shall adopt this approach in subsequent parts of this section to prove
Theorem~\ref{thm:no2}. First we consider the case when $R=2$.

\bl \label{lem:R2} If $R=2$, then $Aut(\Gamma) = T \rtimes \brac{W,\al}$ where $W
=\brac{\omb^k}$.\el \bpf Consider the canonical homomorphism $\varphi: GL(2,p) \mapp
PGL(2,p)$. Let $\oline{\brac{W,\al}} = \varphi(\brac{W,\al})$, $\oline{W} = \varphi(W)$
and $\oline{A_0} = \varphi(A_0)$. Since $A_0 \nsupseteq SL(2,p)$ we have that
$\oline{A_0} \nsupseteq PSL(2,p)$. We want to find the maximal subgroups of $PGL(2,p)$
that contain $\oline{A_0}$

Now $\oline{W} =  WZ/Z$, and $WZ =\brac{\omb^k, \omb^{p+1}}$ has order $p^2-1$ if $k$ is
odd and $\frac{p^2-1}{2}$ if $k$ is even. It follows that \[ \oline{\brac{W, \al}} \cong
\left\{
\begin{array}{ll}
                D_{2(p+1)} &\mbox{when $k$ is odd}\\
                D_{p+1}     &\mbox{when $k$ is even.}
                 \end{array}
                            \right.\]

Also, it follows from the classification of the subgroups of $PGL(2,p)$ and $PSL(2,p)$
(see \cite[p. 417]{Suzuki1}) that either $\oline{A_0} \leq \oline{\brac{\omb, \al}} \cong
D_{2(p+1)}$ (where $\oline{\brac{\omb, \al}} = \varphi(\brac{\omb,\al})$), or
$\oline{A_0} \in \{A_4, S_4, A_5 \}$. In the former case, $A_0 \leq \brac{\omb, \al}$
and so by Lemma~\ref{lem:pre}, $A_0 =\brac{W, \al}$ as required. So suppose that
$\oline{A_0} \in \{A_4, S_4, A_5 \}$. Then, since $\oline{A_0} \geq \Z_{(p+1)/2}$ and $p
\geq 5$, it follows that $p=5$ or 7. Moreover since $\oline{A_0} \geq D_{p+1}$ and
$\oline{A_0} \ngeqslant PSL(2,p)$ if follows that $\oline{A_0} =S_4$ in both cases, and
hence $\oline{A_0}$ is transitive on 1-spaces. Thus $S$ consists of, say, $s$ points
from each 1-space and hence $|S| = (p+1)s$. Since $S=-S$, it follows that $s$ is even.
Further, as $|S| = |V^*|/k$, we have $k=|V^*|/((p+1)s) = (p-1)s$ which divides
$(p-1)/2=2$ or $3$. However since $k \geq 3$, this implies that $p =7$, $k =3$ and
$s=2$. This is impossible since in this case $W=\brac{\omb^3} \cong \Z_{16}$ projects to
$\oline{W} \cong \Z_8$ and $\oline{A_0} =S_4$ has no such subgroup. \epf

Now suppose $R \geq 3$ and $k \mid (p-1)$ with $k \geq 3$, $p \geq 5$. Let $G = T
\rtimes \brac{W, \al}$ and $G_0 = \brac{W,\al} \leq A_0 \leq GL(R,p)$. We shall start
off by defining a primitive prime divisor.

\bdeff \label{def:ppd} \rm{Let $b$ and $x$ be positive integers with $b > 1$. A prime
number $r$ is called a \textit{primitive prime divisor} of $b^x-1$ (or $\ppd$ for short)
if $r$ divides $b^x-1$ but does not divide $b^i-1$ for any integer $i$ with $1 \leq i
\leq x-1$.} \edeff

It was shown by Zsigmondy in \cite{Zsig} that for $x \geq 3$, $b^x-1$ has a primitive
prime divisor unless $(x,b)=(6,2)$. Since $p \geq 5$ and $R \geq 3$, it follows that
$p^R-1$ has a primitive prime divisor, say $r$. By Definition~\ref{def:ppd}, $r$ does
not divide $p-1$, and so, since $k \mid (p-1)$, the prime $r$ divides $|W| =
\frac{p^R-1}{k}$, and hence $r$ divides $|A_0|$.

In \cite{GPPS}, a classification is given of all subgroups of $GL(d,q)$ with order divisible by a
primitive prime divisor of $q^x-1$ for some $x$ satisfying $\frac{d}{2} < x \leq d$. In
particular, the results of \cite{GPPS} specify all subgroups of $GL(R,p)$ with order divisible by
$r$.

Let $a$ be minimal such that $a \geq 1$, $a \mid R$ and $A_0$ preserves on $V$ the
structure of an $a$-dimensional vector space over a field of order $q =p^{R/a}$. Then
$A_0 \leq \Gamma L(a,q)$ acts on $V =V(a,q)$.

\bl \label{lem:mc} If $a=1$ then $Aut(\Gamma) = T \rtimes \brac{W, \al}$ where $W =
\brac{\omb^k}$. \el \bpf If $a=1$ then $A_0 \leq \Gamma L(1,p^R)$, and by
Lemma~\ref{lem:pre}, $A=Aut(\Gamma)=T \rtimes A_0$ with $A_0 = \brac{W, \al}$. \epf

So in what follows, we assume $a \geq 2$. Now $r$ is also a $\ppd$ of $q^a-1=(p^{R/a})^a
-1 =p^R-1$. Let $A_1 := A_0 \ \cap \ GL(a,q)$. Since $W \subseteq GL(a,q)$, it follows
that $W \subseteq A_1$ and $r$ divides $|A_1|$. Also, since $k \mid (p-1)$, the order of
$W$ is divisible by $\frac{p^R-1}{p-1}$, and so $A_1$ contains an element of order
$\frac{p^R-1}{p-1}$. We apply the results in \cite{GPPS} to find the possibilities for
$A_1 \leq GL(a,q)$ ($a \geq 2$), where $A_1$ contains an element of order $r$ and an
element of order $\frac{p^R-1}{p-1}=\frac{q^a-1}{p-1}$. By \cite[Main Theorem]{GPPS},
either (here, $Z =Z(GL(a,q))$)
\newcounter{tt}
\begin{list} {{\rm \textbf{(\Alph{tt})}}}{\usecounter{tt}}
\item $A_1$ belongs to one of the families of Examples 2.1 - 2.5 in
\cite{GPPS}, or
\item $A_1$ is \textit{nearly simple}, that is, $L \leq A_1/(A_1 \cap Z)
\leq Aut(L)$ for some nonabelian simple group $L$, and the examples are listed in
Examples 2.6 - 2.9 in \cite{GPPS}.
\end{list}

We first show that $A_1 \leq GL(a,q)$ is not one of the examples in case (\textbf{A}).

\bl \label{lem:CaseA} Let $a\geq 2$ and suppose $A_1 = A_0 \cap GL(a,q)$. Then case
$\mathbf{(A)}$ does not hold. \el \bpf Suppose $a \geq 2$ and $A_1 \leq GL(a,q)$ belongs
to one of the families of Examples 2.1 - 2.5 in \cite{GPPS}. We shall refer to these
families of examples as E1 - E5, and consider each of them in turn.

(E1: \emph{Classical Examples}) \ Here $A_1$ contains a classical group. Since $A_1$ is
not transitive on the set of non-zero vectors $V^*$, $A_1$ cannot contain $SL(a,q)$ or
$Sp(a,q)$ (as in \cite[Examples 2.1 (a) - (b)]{GPPS}), and since $A_1$ contains the
irreducible element $\omb^k$, either $A_1 \leq Z \circ GO^-(a,q)$ (with $a$ even) or
$A_1 \leq Z \circ GU(a, \sqrt{q})$ (with $a$ odd and $q$ a square) (see \cite{Aron} or
\cite{Huppert70}). Now (see \cite{Aron}) $|W \ \cap \ (Z \circ GO^-(a,q))| \leq
(q^{a/2}+1)(q-1)$ and $|W \ \cap \ (Z \circ GU(a,\sqrt{q}))| \leq (q^{a/2}+1)(q-1)$,
whereas $|W \cap A_1| = |W| \geq \frac{p^R-1}{p-1}$. Thus $p^{R/a} -1 \geq \frac{p^{R/2}
-1}{p-1}$. It follows that $a=2$ and $A_1 \leq Z \circ GO^-(2,q)$. However in this
``degenerate case" of E1, $A_0 \leq Z \circ \Gamma O^-(2,q) = \Gamma L(1,q^2)$,
contradicting the minimality of $a$.

(E2: \emph{Reducible Examples} and E3: \emph{Imprimitive Examples}) \ These families of
groups contain no groups with  order divisible by a $\ppd$ of $q^a-1$ (where $a \geq 2$).

(E4: \emph{Extension Field Examples}) \ In this class of examples, since $r$ is a $\ppd$
of $q^a-1$, there exists a proper divisor $b$ of $a$ such that $A_1 \leq GL(a/b, q^b)
\cdot b$ and so $A_0 \leq \Gamma L(a/b, q^b)$. However this contradicts the minimality
of $a$.

(E5: \emph{Symplectic Type Examples}) \ For this class of examples, we have $a =2^m$ and
$r=a+1$. Hence, $m$ is a power of $2$. Recall that $q =p^{\frac{R}{a}}$, and since $R
\geq 3$ and $a=2^m \mid R$, it follows that $R \geq 4$ and $R$ is even. The group $A_1$
is contained in $Z \circ (S \cdot M_0)$ where $S$ and $M_0$ are as listed in
Table~\ref{tab:no2}.

\begin{table}[ht]
\begin{center}
\begin{tabular}{c c c}
\hline $S$ & $M_0$ & $a,p$ \\ \hline  $4 \circ 2^{1+2m} = \Z_4 \circ D_8 \circ \cdots
\circ D_8$ & $Sp(2m,2)$ & $p \equiv 1$
(mod 4), $a \geq 4$ \\
$2_{-}^{1+2m} =D_8 \circ \cdots \circ D_8 \circ Q_8$ & $O^{-}(2m,2)$ & $a \geq 2$ \\
$2_{+}^{1+2m} =D_8 \circ \cdots \circ D_8$ & $O^{+}(2m,2)$ & $a \geq 4$ \\
\hline
\end{tabular}
\vs \caption{} \label{tab:no2}
\end{center}
\end{table}

Note that the largest element orders in $S$ and $M_0$ are less than or equal to 4 and
$2^{2m}-1 =a^2-1$ respectively (see \cite[Proof of Lemma 2]{Aron}), and $S \cap Z \cong
\Z_2$. Since $A_1$ contains an element of order $\frac{p^R-1}{p-1}= \frac{q^a-1}{p-1}$,
we have \ben \label{eq:N1} \frac{q^a-1}{p-1} \leq 2(a^{2}-1)(q-1). \een Now $q^{a-2} -1
\ < \ \frac{q^a-1}{q^2} \ < \ \frac{q^a-1}{(q-1)^2}  \ \leq \ \frac{q^a-1}{(p-1)(q-1)}$,
and since $q \geq p \geq 5$, we have \be 5^{a-2} -1 \ \leq \ q^{a-2} -1 & < & 2(a^2-1).
\ee Since $a=2^m$ and $m$ is a power of $2$, this implies that either $a=2$, or $(a,q) =
(R,p) = (4,5)$. In the latter case, $r=a+1 =5$ is not a $\ppd$ of $5^4-1$, and we have a
contradiction. Hence $a=2$ and so $r=a+1=3$.

Now for $a =2$, inequality~(\ref{eq:N1}) yields \be p^{R/2}+1 = q+1 & \leq & 6(p-1). \ee
However since $p \geq 5$ and $R= 2a \geq 3$, there are no values of $p$, $R$ for which
this holds.

This completes the proof of Lemma~\ref{lem:CaseA}. \epf

Thus case (\textbf{B}) holds. A subgroup of $GL(a,q)$ is said to be \textit{realisable
modulo scalars} over a proper subfield $\F{q_0}$ of $\F{q}$ if it is conjugate by an
element of $GL(a,q)$ to a subgroup of $GL(a,q_0) \circ Z$ where $Z=Z(GL(a,q))$. When
this occurs, we may replace the subgroup by a conjugate if necessary, and assume that it
lies in $GL(a,q_0) \circ Z$. Choose $q_0$ minimal such that $A_1$ is realisable modulo
scalars over $\F{q_0}$ and as above, assume $A_1 \leq GL(a,q_0) \circ Z$ (see \cite[p.
172]{GPPS}). Note that $q_0 = p^{R/aa_0}$ for some divisor $a_0$ of $\frac{R}{a}$, where
$a_0 \geq 1$.

\bl \label{caseB1} Suppose $A_1 \leq GL(a,q_0) \circ Z$. Then $q =q_0$, that is, $a_0
=1$. \el \bpf The maximum order of the elements in the group $GL(a,q_0)$ is $q_0^a-1$ and
hence the maximum order of the elements of $GL(a,q_0) \circ Z$ is less than or equal to
$\frac{(q_0^a-1)(q-1)}{q_0-1}$. Since $A_1$ contains an element of order
$\frac{q^a-1}{p-1}$, we have \be \frac{q^a-1}{p-1} \leq \frac{(q_0^a-1)(q-1)}{q_0-1}. \ee
Re-arranging the above expression, we have (noting that $q=q_0^{a_0}$ and $q_0 =
p^{R/aa_0}$) \ben \label{ER} p-1  \geq \frac{(q^a-1)(q_0-1)}{(q_0^a-1)(q-1)} \geq
q_0^{a(a_0-1)} \cdot \frac{q_0-1}{q-1} > \frac{q_0^{a(a_0-1) -(a_0-1)}}{2}. \een Since
$2(p-1) < p^2$, it follows that \be \frac{R}{aa_0}(a-1)(a_0-1) = \frac{R}{aa_0} (a(a_0-1)
-(a_0 -1)) < 2. \ee

Suppose $a_0 \geq 2$. Then since $\frac{R}{aa_0} \geq 1$, we have $a-1 \leq
\frac{R}{aa_0} (a-1)(a_0-1) <2$ and hence $a=2$. Thus the inequality in (\ref{ER})
becomes $p-1 \geq \frac{q+1}{q_0+1}$ and hence $q+1 \leq (p-1)(q_0+1) = pq_0 - q_0 + p
-1 < pq_0 \leq q$, which is a contradiction. So $a_0 =1$. \epf

We now consider the examples from Examples 2.6 - 2.9 in \cite{GPPS} (see also
Tables~\ref{tab:no3} - \ref{tab:no7}) of subgroups $A_1$ of $GL(a,q)$ with order
divisible by a $\ppd$ $r$ of $q^a-1$, where $a \geq 2$ and $q=q_0=p^{R/a}$. Here, $L
\leq A_1/(A_1 \cap Z) \leq Aut(L)$ for some non-abelian simple group $L$.

\vspace{3cm}
\begin{table}[ht]
\begin{center}
\begin{tabular}{|c|c|}
\hline (1) & $A_n \leq A_1 \leq S_n \times Z$; \ $a =r-1$.\\   & Furthermore, $n = a+1$
if $p \nmid n$ or $n =a+2$ if $p \mid n$. \\ \hline
\end{tabular}
\caption{Example 2.6 (a) of \cite{GPPS}} \label{tab:no3}
\end{center}
\end{table}

\vspace{5cm}
\begin{center}

(Turn over for Tables~\ref{tab:no4} - \ref{tab:no7}.)

\end{center}

\newpage

\begin{landscape}
\begin{table}
\begin{center}
\begin{tabular}{|c|c|c|c|}
\hline
(2) & $L = A_7$& $a=4$ & $r=5$ and $q_0 =q =7$.\\

(3) & $L = A_7$& $a=4$ & $r=5$, $p \equiv 1,2,4$ $($mod $7)$ and
$q_0=q=p$. \\

(4) & $L = A_6$& $a=4$ & $r=5$, $p \geq 7$, $p \equiv \pm 2$ (mod 5)
and $q_0=q=p$. \\

(5) & $L= A_5$& $a=2$& $r=5$, $p \equiv 4$ (mod 5) and $q_0 =q =p$. \\

(6) & $L= A_5$& $a=2$& $r=5$, $p \equiv \pm 2$ (mod 5) and $q_0 =q
=p^2$. \\

(7) & $L= A_5$& $a=2$& $r=3$, $p=5$ or $p \equiv 11, 14$ (mod 15)
and $q_0 =q =p$. \\

(8) & $L= A_5$& $a=4$& $r=5$, $p \geq 7$, $p \equiv \pm 2$ (mod 5)
and $q_0 =q =p$.\\

(9) & $L=A_7$& $a=3$& $r=7$, $p=5$ and $q_0 =q=25$. \\

(10) & $L=A_7$& $a=6$& $r=7$, $p \equiv 1$ (mod 6) and $q_0 =q=p$. \\

(11) & $L=A_7$& $a=6$& $r=7$, $p \equiv 1, 7$ (mod 24) and $q_0 =q=p$. \\
\hline
\end{tabular}
\caption{Example 2.6 (b) of \cite{GPPS}} \label{tab:no4}
\end{center}
\end{table}

\begin{table}
\begin{center}
\begin{tabular}{|c|c|c|c|}
\hline
(12) & $L=M_{11}$& $a=10$& $r=11$, $p \neq 11$ and $q_0 =q=p$. \\

(13) & $L=M_{12}$& $a=10$& $r=11$, $p \equiv 1, 3$ (mod 8) and $q_0
=q=p$. \\

(14) & $L=M_{22}$& $a=10$& $r=11$, $p=7$ or $p \equiv 1, 2, 4$ (mod
7) and $q_0 =q=p$. \\

(15) & $L=M_{23}$& $a=22$& $r=23$, $p \neq 2, 3$ and $q_0 =q=p$. \\

(16) & $L=J_{2}$& $a=6$& $r=7$, $p=5$ or $p \equiv \pm 1$ (mod 5)
and $q_0 =q=p$. \\

(17) & $L=J_{3}$& $a=18$& $r=19$, $p \equiv 1$ or $4$ (mod 15) and
$q_0 =q=p$. \\

(18) & $L=Ru$& $a=28$& $r=29$, $p \equiv 1$ (mod 4) and $q_0 =q=p$.
\\

(19) & $L=Suz$& $a=12$& $r=13$, $p \equiv 1$ (mod 6) and $q_0 =q=p$.
\\ \hline
\end{tabular}
\caption{Example 2.7 of \cite{GPPS}} \label{tab:no5}
\end{center}
\end{table}

\begin{table}
\begin{center}
\begin{tabular}{|c|c|c|c|}
\hline
(20) & $L=G_{2}(4)$& $a=12$& $r=13$, $p >2$ and $q_0 =q=p$. \\

(21) & $L=PSU(4,2)$& $a=4$& $r=5$, $p \equiv 1$ (mod 6) and $q_0
=q=p$. \\

(22) & $L=PSU(4,3)$& $a=6$& $r=7$, $p \equiv 1$ (mod 6) and $q_0
=q=p$. \\

(23) & $L=PSL(3,4)$& $a=6$& $r=7$, $p \equiv 1$ (mod 6) and $q_0
=q=p$. \\
\hline
\end{tabular}
\caption{Example 2.8 of \cite{GPPS}} \label{tab:no6}
\end{center}
\end{table}

\begin{table}
\begin{center}
\begin{tabular}{|c|c|c|c|}
\hline
(24) & $L=PSL(n,s)$ where $n \geq 3$& $a = \frac{s^n-1}{s-1} -1$& $r =
\frac{s^n-1}{s-1}$, $n$ is prime and \\ & & & $s$ is a prime power such that $gcd(s,p)
=1$.\\ \hline

(25) & $L=PSU(n,s)$ where $n \geq 3$& $a = \frac{s^n+1}{s+1} -1$& $r =
\frac{s^n+1}{s+1}$, $n$ is prime and \\ & & & $s$ is a prime power such that $gcd(s,p)
=1$. \\ \hline

(26) & $L=PSp(2n,s)$& $a = \frac{1}{2}(s^n-1)$& $r = \frac{1}{2}(s^n + 1)$, $n =2^b \geq
2$ and \\ & & &  $s$ is a odd prime power such that $gcd(s,p) =1$. \\ \hline

(27) & $L=PSL(2,s)$ where $s \geq 7$& $a = s$& $r = s+1$, $s =2^c$,
\\ & & & $c=2^{c'}$ and $gcd(s,p) =1$. \\ \hline

(28) & $L=PSL(2,s)$ where $s \geq 7$& $a = s-1$& $r = s$, \\ & & & $s$ is prime and
$gcd(s,p) =1$.
\\ \hline

(29) & $L=PSL(2,s)$ where $s \geq 7$& $a = \frac{1}{2}(s-1)$& $r = s$, \\ & & & $s$ is
prime and $gcd(s,p) =1$. \\ \hline

(30) & $L=PSL(2,s)$ where $s \geq 7$& $a = \frac{1}{2}(s-1)$& $r = \frac{1}{2}(s+1)$, and
\\ & & & $s$ is a odd prime power such that $gcd(s,p) =1$. \\ \hline
\end{tabular}
\caption{Example 2.9 of \cite{GPPS}} \label{tab:no7}
\end{center}
\end{table}
\end{landscape}

\bl \label{lem:CB} Let $a\geq 2$ and $A_1 = A_0 \cap GL(a,q)$. Then case $\mathbf{(B)}$
does not hold. \el \bpf Suppose $a \geq 2$ and $A_1 \leq GL(a,q)$ belongs to one of the
families of examples in case $\mathbf{(B)}$. We consider all the examples Examples 2.6 -
2.9 in \cite{GPPS}.

In Example $2.6 (a)$ of \cite{GPPS} (see Table~$\ref{tab:no3}$), $A_n \leq A_1 \leq S_n \times Z$,
$r=a+1$, $n \geq 5$ and $n=a+1$ if $p \nmid n$, or $n=a+2$ if $p \mid n$. Now $A_1$ contains an
element of order $\frac{p^R-1}{p-1}$, which is a multiple of $r$. However, since $n =r$ or $r+1$
and $r$ is prime, any element of $S_n$ of order divisible by $r$ is an $r$-cycle and hence has
order equal to $r$. Thus an element of $S_n \times Z$ of order divisible by $r$ has order at most
$r(q-1)$. Hence $r(q-1) \geq \frac{p^R-1}{p-1}$, so $(a+1)(p-1) = r(p-1) \geq \frac{p^R-1}{q-1} =
\frac{q^a-1}{q-1}$. It follows that $a=2$ and so $n \leq 4$, which contradicts $n \geq 5$.

In Example $2.6 (b)$ of \cite{GPPS} (see Table~$\ref{tab:no4}$), $L =A_n$ and $(n,r)$ is $(7,7)$,
$(7,5)$, $(6,5)$, $(5,5)$ or $(5,3)$. In $Aut(L)$, an element of order a multiple of $r$ has order
at most $\delta \cdot r \leq 2r$ and hence an element of $A_1$ of order $\frac{p^R-1}{p-1}$ has
order at most $(q-1) \delta \cdot r$. Thus \ben \label{eq:N2} \delta \cdot r(p-1) \geq
\frac{p^R-1}{q-1} =\frac{q^a-1}{q-1}. \een More details on values of $a$, $p$ and $q$ are given in
Table~\ref{tab:no4}. First suppose that $a=6$. Then $r \leq p =q$, so $R=6$ and (\ref{eq:N2})
gives $2p^2 > \frac{p^R-1}{q-1} =\frac{p^6-1}{p-1}$, which is not true. Next, suppose that $a=3$.
Then $r=7$, $q=p^2 =25$ and so $R=6$ and $8r = 2r(p-1) \geq \frac{5^6-1}{24}$, which is not the
case. Now suppose that $a=4$. Then $r=5$ and $q=p \geq 7$. So $10 (p-1) \geq \frac{p^R-1}{q-1} =
\frac{p^4-1}{p-1}$, which is again not true. Hence $a=2$, $L=A_5$ and as $R \geq 3$, we must have
$R \geq 4$. Hence $q \geq p^2$, so line (6) of Table~\ref{tab:no4} holds. Thus $q=p^2$, $r=5$, and
(\ref{eq:N2}) fails since $\delta =1$ in this case.

In Example $2.7$ of \cite{GPPS} (see Table~$\ref{tab:no5}$), one of the lines of
Table~\ref{tab:no5} holds. In particular $q=p$ and so $A_1$ contains an element of order
$\frac{p^a-1}{p-1} > p^{a-1}$. However (see the \textsc{Atlas} \cite{Atlas}) an element of $A_1$
of order a multiple of $r$ has order at most $2(p-1)r$ and this is a contradiction for all lines
of Table~\ref{tab:no5}.

In Example $2.8$ of \cite{GPPS} (see Table~$\ref{tab:no6}$), one of the lines of
Table~\ref{tab:no6} holds. In particular $q=p$ and so $A_1$ contains an element of order
$\frac{p^a-1}{p-1}$. This quantity is at most $(p-1) \delta_L$ where $\delta_L$ is the largest
order for an element of $Aut(L)$ of order divisible by $r$, which is (see \textsc{Atlas}
\cite{Atlas}): 13, 10, 28 and 21 in lines (20), (21), (22) and (23) respectively. This gives a
contradiction in all cases.

In Example $2.9$ of \cite{GPPS} (see Table~$\ref{tab:no7}$), one of the lines of
Table~\ref{tab:no7} holds. Following the notation used in \cite{Aron}, we let $m(K)$ denote the
maximum of the orders of the elements of a finite group $K$. Now $\omb^k$ has order modulo scalars
at least $\frac{p^R-1}{k} \cdot \frac{1}{q-1} \geq \frac{p^R-1}{p-1} \cdot \frac{1}{q-1}$ and
hence, \ben \label{sharp} q^{a-2}  < \frac{q^a-1}{q-1} \cdot \frac{1}{p-1} \leq m(Aut(L)). \een We
shall make use of the inequality in (\ref{sharp}) to go through lines (24) to (30) of
Table~\ref{tab:no7}.

\textit{Line $(24)$ of Table~$\ref{tab:no7}$}: Here $L=PSL(n,s)$ ($n \geq 3$) and
$a=\frac{s^n-1}{s-1} -1$, where gcd$(s,p)=1$ and $s$ is a prime power. Thus $s = x^c$ where $x$ is
a prime and $x \neq p$. Since $n \geq 3$, we have $a \geq \frac{s^3-1}{s-1} -1 = s^2 + s \geq 6$.
Now $m(Aut(L)) \leq m(GL(n,s)) \cdot 2c = 2c(s^n-1) \leq s(s^n-1)$. Using inequality
(\ref{sharp}), we have $q^{a-2} <  \frac{q^a-1}{q-1} \cdot \frac{1}{p-1} \leq s(s^n-1)
=s(s-1)(a+1) < s^2(a+1) < a(a+1)$. Thus $q^{a-2} < a(a+1)$, which is impossible for $q \geq 5$ and
$a \geq 6$.

\textit{Line $(25)$ of Table~$\ref{tab:no7}$}: Here, $L=PSU(n,s)$ ($n \geq 3$; $n$ is prime) and
$a=\frac{s^n+1}{s+1} -1$ where $gcd(s,p)=1$. Again we let $s =x^c$ where $x$ is a prime, $x \neq
p$ and we have $a \geq \frac{s^3+1}{s+1} -1 = s^2-s \geq s$. Note that $m(Aut(L)) \leq
m(GL(n,s^2)) \cdot 2c = 2c (s^{2n} - 1) < 2c (s^n+1)^2$. Since $s^n+1 =(a+1)(s+1)$, we have
$m(Aut(L)) < 2c (s+1)^2(a+1)^2 \leq 2a(a+1)^4$ (since $c \leq s \leq a$).

Using inequality (\ref{sharp}), we have \be q^{a-2} \ < \ \frac{q^a-1}{q-1} \cdot \frac{1}{p-1} &
< & 2a(a+1)^4 \ee and hence, since $q \geq p \geq 5$, this implies that $a \leq 9$. Also, as $r =
\frac{s^n+1}{s+1} = a+1$ is prime, $n$ is odd and $a \leq 9$, it follows that $(r,a,s,n)$ is
either $(3,2,2,3)$ or $(7,6,3,3)$. However in the first case, $L = PSU(3,2)$ is not a nonabelian
simple group. Hence $(r,a,s,n) = (7,6,3,3)$ and $L = PSU(3,3)$. Using the \textsc{Atlas}
\cite{Atlas}, we know that $m(Aut(PSU(3,3)))=12$ and by inequality (\ref{sharp}),
$\frac{q^a-1}{q-1} \cdot \frac{1}{p-1} \leq 12$. However, for $a =6$ and $q \geq p \geq 5$, this
inequality does not hold.

\textit{Line $(26)$ of Table~$\ref{tab:no7}$}: Here, $L=PSp(2n,s)$ and $a=\frac{1}{2}(s^n-1)$
where $gcd(s,p)=1$. Also, $s$ is an odd prime power and $n =2^b \geq 2$. Let $s=x^c$. Now
$m(Aut(\Gamma)) \leq c(s^{2n}-1) < s(s^{2n} -1) $. Since $s^n =2a+1$ and by inequality
(\ref{sharp}), we have \be q^{a-2} \ < \ \frac{q^a-1}{q-1} \cdot \frac{1}{p-1} & < & s(s^{2n} -1)
\ < \ s^{3n} \ = \ (2a+1)^3. \ee Since $q \geq p \geq 5$, this implies that $a \leq 6$, and since
$r = \frac{1}{2} (s^n+1) =a+1$ is prime, $n=2^b \geq 2$, and $s$ is odd, it follows that
$(r,a,s,n)$ is $(5,4,3,2)$. Thus $L=PSp(4,3)$ and using the \textsc{Atlas} \cite{Atlas},
$m(Aut(PSp(4,3))) = 12$, which is not greater than $ \frac{q^a-1}{q-1} \cdot \frac{1}{p-1}$ for
$a=4$ and $q \geq p \geq 5$.

\textit{Lines $(27) - (30)$ of Table~$\ref{tab:no7}$}: In all cases, $L =PSL(2,s)$ with $s \geq
7$, $s \neq p$ and $a \geq \frac{s-1}{2} \geq 3$. Then $m(Aut(L)) =s+1$ (see \cite{Aron}) and
hence inequality (\ref{sharp}) yields \be q^{a-2} < s+1 \leq 2a+2. \ee Since $q \geq5$ and $a \geq
\frac{s-1}{2} \geq 3$, this implies that $a =3 = \frac{s-1}{2}$ and $q \leq 7$. Thus line (29) or
(30) holds with $L = PSL(2,7)$, and $p \neq s =7$. So $q=p=5$ and $a=R=3$. However the only $\ppd$
of $p^R-1 = 5^3-1$ is 31 whereas by Table~\ref{tab:no7}, $r \leq 7$.

We have exhausted all possibilities for $A_1$ and so Lemma~\ref{lem:CB} is proved. \epf

We formalise the proof of Theorem~\ref{thm:no2}.

\vs \noindent \textit{Proof of Theorem~\ref{thm:no2}.} \ Suppose $\Gamma =\GP =Cay(V,S)$
is a generalised Paley graph where $k \geq 2$ and $k \mid (p-1)$. Then by
Lemma~\ref{lem:conns}, $\Gamma$ is connected. Let $A =Aut(\Gamma)$. If $k=2$, then $A = T
\rtimes \brac{\omb^2, \al}$ (see for example \cite{Peisert2001}).

So suppose that $k \geq 3$. Since $k \mid (p-1)$ and $\Gamma$ is connected, then by
Lemma~\ref{lem:newstuff}, we have $A = T \rtimes A_0 \leq AGL(R,p)$ and $\brac{W, \al}
\leq A_0$ where $W=\brac{\omb^k}$. If $R=2$, then by Lemma~\ref{lem:pre} and
Lemma~\ref{lem:R2}, we have $A = T \rtimes \brac{W, \al}$ and so the theorem holds.

Now suppose $R \geq 3$, and let $a$ be minimal such that $A_0 \leq \Gamma L(a,p^{R/a})$.
By Lemmas~\ref{lem:CaseA} and \ref{lem:CB}, $a=1$, and by Lemma~\ref{lem:mc}, $A_0
=\brac{W, \al}$. \qedd

\section{Cyclotomic association scheme} \label{sec:cyc}
In this section, we briefly mention the relationship between a generalised Paley graph
and the class of graphs associated with a symmetric cyclotomic association scheme. We
start by defining a symmetric association scheme.

\bdeff \label{def:AS} {\rm \cite[Chapter 2, p. 43]{bcn} } \rm{A \textit{symmetric $k$-class
association scheme} is a pair $(V,\M{R})$ such that \bnum
\item the set $\M{R} =\{R_0, R_1, \ldots, R_k\}$ is a partition of $V^2 = V \times V$
($R_0, R_1, \ldots, R_k$ are called \textit{classes} of the association scheme $(V,\M{R})$);
\item $R_0=\{(x,x) \mid x\in V\}$;
\item $R_i=R_i^T$ (that is, $(x,y) \in R_i$ $\imply$ $(y,x) \in R_i$) for all $i \in
\{0,1,\ldots, k\}$;
\item there are numbers $p_{ij}^h$ (called the \textit{intersection numbers} of the
scheme) such that for any pair $(x,y) \in R_h$, the number of $z \in V$ with $(x,z) \in
R_i$ and $(z,y) \in R_j$ equals $p_{ij}^h$. \enum} \edeff


Equivalently, we may see each class of a symmetric association scheme $(V,\M{R})$ as an undirected
graph with vertex set $V$, that is, $\Gamma_i =(V, E_i)$ where $E_i =\{\{x,y\} \mid (x,y) \in
R_i\}$ (see for example \cite[Chapter 12]{Godsil93}). Thus from now on we shall identify the $k$
classes of a symmetric association scheme $(V,\M{R})$ with a set of undirected graphs $\Gamma_1,
\ldots, \Gamma_k$. Note that the graphs $\Gamma_i$ need not be isomorphic to one another. We say
that an association scheme is \textit{primitive} if each of the graphs $\Gamma_i$ is connected,
otherwise it is called \textit{imprimitive}.

Let $V=\F{p^R}$ and suppose $S =\brac{\om^k} \subseteq V^* = \F{p^R}^*$ where $k \mid
(p^R-1)$. Define the \textit{cyclotomic association scheme} (or \textit{cyclotomic
scheme}) on $V$ by $\M{R} = \{R_i \mid 0 \leq i \leq k \}$ where for $1 \leq i \leq k$,
$R_i=\{(x,y) \mid y-x \in S \om^i\}$. We shall denote a cyclotomic scheme by Cyc$(V,
\M{R})$. A cyclotomic scheme is symmetric when $\frac{q-1}{k}$ is even, or when $q$ is
even. If $k=2$, then the cyclotomic schemes are just the association schemes of the
Paley graphs (see \cite[p. 66]{bcn}). Note that in a symmetric cyclotomic scheme Cyc$(V,
\M{R})$, the classes $\Gamma_i$ of the scheme (identified as undirected graphs on $V$)
are all isomorphic to one another.

It is easily seen that a generalised Paley graph $\Gamma = \GP$ is isomorphic to each of
the graphs $\Gamma_i$ (where $1 \leq i \leq k$) of a $k$-class symmetric cyclotomic
scheme Cyc$(V, \M{R})$ where $V =\F{p^R}$. By Theorem~\ref{thm:add}, they are connected
if and only if $k$ is not a multiple of $\frac{p^R-1}{p^r-1}$ for any proper divisor $r$
of $R$. This divisibility condition is therefore a necessary and sufficient condition
for Cyc$(V, \M{R})$ to be primitive.

Suppose now Cyc$(V, \M{R})$ is a primitive $k$-class symmetric cyclotomic scheme where
$V =\F{p^R}$. Then the set of $k$ undirected graphs $\Gamma_1, \ldots, \Gamma_k$
associated with Cyc$(V, \M{R})$ are all connected, and are isomorphic to a connected
generalised Paley graph $\GP$. By Theorem~\ref{thm:A}, the graphs $\Gamma_i$ are Hamming
graphs if and only if there exists a divisor $a$ of $R$ with $1 \leq a < R$ such that $k
=\frac{a(p^R-1)}{R(p^a-1)}$. Thus, in the language of association scheme, we may restate
Theorems \ref{thm:add} and \ref{thm:A} as follow (note that a proper divisor $a$ of $R$
is one satisfying $1 \leq a <R$).

\bt \label{thm:CYC} Suppose {\rm Cyc}$(V, \M{R})$ is a $k$-class symmetric cyclotomic
scheme where $V =\F{p^R}$. Let $\Gamma_1, \ldots, \Gamma_k$ be the set of $k$ undirected
graphs associated with {\rm Cyc}$(V, \M{R})$. Then$:$ \bnum
\item[$1.$] {\rm Cyc}$(V, \M{R})$ is primitive if and only if $k$ is not a multiple of
$\frac{p^R-1}{p^a-1}$ for every proper divisor $a$ of $R$.
\item[$2.$] Each $\Gamma_i$ is a Hamming graph if and only if $k=\frac{a(p^R-1)}{R(p^a-1)}$
for some proper divisor $a$ of $R$. \enum \et
\chapter{Methods and background theory} \label{c3}

In this chapter, we develop the general theory in our investigation of edge-transitive
homogeneous factorisations of complete graphs. In particular, we show that in an
edge-transitive homogeneous factorisation $(M,G,V,\M{E})$ of a complete graph, the group
$G$ is necessarily $2$-homogeneous. Equipped with this result, we then spell out the
approach and the methods we will use to characterise the factor graphs arising in such
factorisations.

\section{Some basics} \label{c3:strategy}
We shall start with the following result which is fundamental to our study of
edge-transitive homogeneous factorisations of complete graphs.

\bl \label{G_2homo} Suppose $(M,G,V,\mathcal{E})$ is a homogeneous factorisation of
$K_n$ such that $M$ is edge-transitive on each factor. Then $G$ is a $2$-homogeneous
permutation group on $V$. \el \bpf Let $\{x,y\}$ and $\{u,v\}$ be two unordered pairs
from the vertex set $V$ of $K_n$. If $\{x,y\}$ and $\{u,v\}$ are from the same factor,
say $\Gamma_i$, then since each factor is $M$-edge-transitive, there exists $m \in M
\leq G$ such that $\{x,y\}^{m} = \{u,v\}$. If $\{x,y\}$ and $\{u,v\}$ are not of the
same factor, then there exists $g \in G$ that maps $\{x,y\}$ to a pair $\{x',y'\}$ of
the same factor as $\{u,v\}$. We then apply some $m \in M \leq G$ to reach $\{u,v\}$.
Thus $G$ is 2-homogeneous. \epf

All finite 2-homogeneous permutation groups have been classified. They consist of all finite
2-transitive permutation groups (up to permutation isomorphism) listed in Theorem~\ref{2_trxclass}
and all finite 2-homogeneous but not 2-transitive permutation groups described by Kantor (see
Theorem~\ref{2hom}). As we shall see, it will not be necessary to use details of the latter
classification.

In light of Lemma~\ref{G_2homo}, we can now take a finite 2-homogeneous permutation group
$G$, find all its possible normal subgroups $M$ that are not 2-homogeneous, and for each
$M$, determine the structure of the $M$-edge-transitive (or $M$-arc-transitive) factors.
We shall adopt this general approach to classify all edge-transitive homogeneous
factorisations of complete graphs.

The next result tells us that if $G$ is a 2-transitive permutation group on $V$ with a
normal subgroup $M$ of even order, then all $M$-orbitals in $V$ are self-paired. In
particular, we show that given an edge-transitive homogeneous factorisation
$(M,G,V,\M{E})$ where the groups $G$ and $M$ have such properties, the factor graphs are
always $M$-arc-transitive.

\bp \label{selfpaired} Let $(M,G,V,\M{E})$ be an edge-transitive homogeneous
factorisation of $K_n =(V,E)$ of index $k$ with $|V|=n$. Also, let $\mathcal{E}= \{E_1,
\ldots, E_k \}$ with corresponding factors $\Gamma_i =(V,E_i)$. Suppose $G$ is a
$2$-transitive permutation group on $V$ and $M$ is a transitive normal subgroup of $G$
such that $|M|$ is even. Then the following hold. \bnum
\item[$1.$] All nontrivial $M$-orbitals in $V$ are self-paired and the number of such
$M$-orbitals equals $k$. Let them be $\Or_1, \ldots, \Or_k$.
\item[$2.$] The partition $\mathcal{E}= \{E_1, \ldots, E_k \}$ may be labelled such that
$E_i =\{ \{x,y\} \mid (x,y) \in \Or_i \}$ for each $i \in \{ 1, \ldots, k\}$.
Furthermore, the factors $\Gamma_i$ are all $M$-arc-transitive of valency
$\frac{n-1}{k}$. \enum \ep \bpf (1) \ By \cite[Theorem 16.5]{wie}, a finite permutation
group has a nontrivial self-paired orbital if and only if it is of even order. Thus the
group $M$ has at least one nontrivial self-paired orbital in $V$. Since $G$ is
2-transitive and $M$ is normal in $G$, it follows by Lemma~\ref{gin} that $G$ leaves the
set of nontrivial $M$-orbitals invariant and permutes the set of all nontrivial
$M$-orbitals transitively. Thus all nontrivial $M$-orbitals are self-paired.

Since $M$ acts edge-transitively on each $\Gamma_i$ and all nontrivial $M$-orbitals are
self-paired, it follows by Lemma~\ref{lem:orbgraph} that each $E_i$ corresponds to a nontrivial
self-paired $M$-orbital $\Or_i$ in $V$ and $E_i =\{ \{x,y\} \mid (x,y) \in \Or_i \}$. It follows
immediately that the number of nontrivial self-paired $M$-orbitals equals $|\M{E}| =k$.

(2) \ The first statement follows from the proof of part (1). The part about $\Gamma_i$ being
$M$-arc-transitive follows from Lemma~\ref{lem:orbgraph} (1). Also, since all $\Gamma_i$ ($1 \leq
i \leq k$) are vertex-transitive and isomorphic to one another, it follows that a vertex in
$\Gamma_i$ (for any $i \in \{1,\ldots, k\}$) has the same valency as all other vertices in
$\Gamma_i$, as well as $\Gamma_j$ where $j \neq i$ and $j \in \{1,\ldots, k\}$. Now every vertex
in $\Gamma_i$ is also a vertex in $K_n$ (of valency $n-1$). Consider a vertex $v$ in $K_n$. Then
the set of vertices adjacent to $v$ in $K_n$ is a (disjoint) union of the sets of vertices
adjacent to $v$ in $\Gamma_i$ for all $i$. As the vertices in $\Gamma_i$ (for all $i$) have equal
valency, it follows that each $\Gamma_i$ is regular of valency $(n-1)/k$. \epf

In what follows, we will use the list of finite 2-transitive permutation groups found in
Theorem~\ref{2_trxclass}. In the case of finite 2-homogeneous but not 2-transitive
permutation groups, we use the fact that they are primitive groups of affine type. In
our treatment of this case in Section~\ref{c3:affineGP}, we show that it is sufficient
to only consider the 2-transitive affine case (see Proposition~\ref{fact2}). We first
look at the case when $G$ is an almost simple 2-transitive permutation group, that is,
the list of groups with minimal normal subgroups $N$ listed in Theorem~\ref{2_trxclass}
(1).

\section{Almost simple case} \label{c3:AS}
Let $G$ be a 2-homogeneous group on $V$ that is not of affine type. Then by
Theorems~\ref{2_trxclass} and \ref{2hom}, $G$ has a normal subgroup $N$ which is a nonabelian
simple group, and $N$ is listed in Theorem~\ref{2_trxclass} (1). Moreover, $N$ is also
2-transitive, except in case 1(g) where $N=PSL(2,8)$ (see \cite{Cameron81} or \cite[p.
245-253]{DM_b96}). (We note that $N=PSL(2,8)$ is transitive on $V$ and $|PSL(2,8)|=504$ is even.)

Suppose that $(M,G,V,\M{E})$ is an edge-transitive homogeneous factorisation of
$K_n=(V,E)$ of index $k >1$, and let $\M{E} =\{E_1, \ldots, E_k\}$. Since the group $M$
is a normal subgroup of $G$ and acts transitively on $V$, it must contain $N$. Thus
$|M|$ is even, and by Lemma~\ref{selfpaired}, $M$ has $k$ nontrivial self-paired
orbitals $\Or_1, \ldots, \Or_k$ in $V$ and each $E_i =\{\{x,y\} \mid (x,y) \in \Or_i\}$.
Furthermore, since $k >1$, $M$ is not 2-transitive and so $N=PSL(2,8) \leq M$ acts on
$28$ points (hence we have $|V|=28$ and $K_{28} =(V,E)$). The next result shows that an
edge-transitive homogeneous factorisation $(M,G,V,\mathcal{E})$ of $K_{28}$ where $G$ is
almost simple has index $k=3$.

\bp \label{psl28} Let $(M,G,V,\mathcal{E})$ be an edge-transitive homogeneous
factorisation of $K_n$ of index $k$. Suppose $G$ is an almost simple $2$-transitive
permutation group on $V$ with minimal normal subgroup $N$. Then $M = N = PSL(2,8)$,
$G=Aut(N)=P\Gamma L(2,8)$, $k=3$ and $n=28$. Furthermore, each factor is an
$M$-arc-transitive graph of valency $9$. \ep \bpf By the remarks above, $N =PSL(2,8)$
and $|V| =n=28$. Since $N=PSL(2,8)$ is a maximal subgroup of $Aut(N)$, we have $M=N$ and
$Aut(N)=G$. Now the point stabiliser $N_{v}$ has three orbits in $V \setminus \{v \}$ of
equal length 9 (see for example, \cite[p. 24]{Cameron_b99}). Since there is a natural
one-to-one correspondence between the orbits of $N_{v}$ in $V \setminus \{v\}$ and the
$N$-orbits in $V^{2} \setminus \{(v, v) \mid v \in V\}$, we have exactly three nontrivial
$N$-orbitals in $V$. Also, as $|N|$ is even, then by Proposition~\ref{selfpaired} (1),
all three nontrivial $N$-orbitals are self-paired and so $k=3$. Thus by
Proposition~\ref{selfpaired} (2), each factor is $N$-arc-transitive of valency $9$. \epf

\brk \label{rk:psl28} \rm{We used \textsc{Magma} (see Appendix \textbf{A2}) to compute
the full automorphism group of the above arc-transitive homogeneous factor $\Gamma$ and
found that $|Aut(\Gamma)| = 504 =|PSL(2,8)|$ (note that in \textsc{Magma}, the primitive
group $PSL(2,8)$ on 28 points is written as \verb+PrimitiveGroup(28,2)+). Thus it
follows that $Aut(\Gamma) =PSL(2,8)$. Now no subgroup of $Aut(\Gamma)$ acts regularly on
$V$ and so by Lemma~\ref{lem:cayreg}, the factor $\Gamma$ is non-Cayley.} \erk

\section{Affine case} \label{c3:affineGP}
In this section, we will undertake a general study of the case where $G$ is an affine 2-homogeneous
permutation group. In this case, $G$ is either a 2-homogeneous but not 2-transitive permutation
group (see Theorem~\ref{2hom}), or $G$ belongs to the list of all finite 2-transitive permutation
groups (up to permutation isomorphism) in Theorem~\ref{2_trxclass} (2). In both cases, $G$ contains
a unique minimal normal subgroup $N$ which is regular and is an elementary abelian $p$-group for
some prime $p$, so $N=\Z_p^R$.

Let $(M,G,V, \M{E})$ be an edge-transitive homogeneous factorisation of $K_n=(V,E)$ where $G$ is
an affine 2-homogeneous permutation group with a unique minimal normal subgroup $N = \Z_p^R$.
Since $M$ is a transitive normal subgroup of $G$, it must contain $N$. Furthermore, as $N$ is
regular on $V$, we may identify the vertex set $V$ with a finite vector space (see
Section~\ref{c2:perm_gp}). Let $q$ be maximal such that $G$ preserves on $V$ the structure of a
vector space over a field $\F{q}$ of order $q$. Then $V = \F{q}^a =V(a,q)$ where $q^a =p^R$, $a
\geq 1$. Finally, we identify the regular normal subgroup $N$ with the translation group $T$ which
acts regularly and additively on $V=V(a,q)$ (see Section~\ref{c2:perm_gp}).

Thus we can denote the groups $G$, $M$ and $N$ by: \[ G:= T \rtimes G_0, \hspace{0.7cm}
M:= T \rtimes M_0 \hspace{0.7cm} {\rm and} \hspace{0.7cm} N:=T, \] where $G_0, M_0 \leq
\Gamma L(a,q)$ (stabilisers of the zero-vector in $G$ and $M$ respectively). So from now
on, whenever we have an edge-transitive homogeneous factorisation $(M,G,V,\M{E})$ of
$K_n=(V,E)$ such that $G$ is an affine 2-homogeneous permutation group, we shall adopt
the above notations for $G$ and $M$.

\bl \label{affcayfac} Let $(M,G,V,\M{E})$ be an edge-transitive homogeneous factorisation of
$K_n=(V,E)$ of index $k$ with factors $\{\Gamma_1, \ldots, \Gamma_k \}$ $($and so by definition,
$G$ acts transitively on $\M{E}= \{E\Gamma_1, \ldots, E\Gamma_k\})$. Suppose that $G =T \rtimes
G_0$ is an affine $2$-homogeneous permutation group on $V=V(a,q)$. Then the following hold. \bnum
\item[$1.$] The factors are Cayley graphs of $V$, that is, $\Gamma_i =Cay(V, S_i)$ where $S_i =-S_i$
and $E\Gamma_i =\{ \{u,v\} \mid v-u \in S_i \}$ for all $i \in \{1, \ldots, k\}$. Moreover, the
set $\M{S} := \{S_1, \ldots, S_k\}$ forms a partition of $V \setminus \{\0\}$.
\item[$2.$] The group $G_0$ acts transitively on $\mathcal{S}:=\{S_1, \ldots, S_k\}$. \enum
\el \bpf (1) \ The $M$-edge-transitive homogeneous factors admit the translation group
$T$ as a regular subgroup of automorphisms. So by Lemma~\ref{lem:cayreg}, the undirected
factors are Cayley graphs of $V$. Hence $\Gamma_i = Cay(V,S_i)$, where $E\Gamma_i =\{
\{u,v\} \mid v-u \in S_i \}$ and $S_i =-S_i$. Observe that $S_i =\{v \mid \{\0, v\} \in
E\Gamma_i\}$, and since $\M{E} =\{E\Gamma_1, \ldots, E\Gamma_k\}$ forms a partition of
$E$ of $K_n$, it follows that $\M{S} =\{S_1, \ldots, S_k\}$ forms a partition of $V
\setminus \{\0\}$.

(2) \ From part (1), we know that $\Gamma_i = Cay(V,S_i)$, and by definition of
$\Gamma_i$, $T$ fixes setwise each $E\Gamma_i$. Suppose $g  \in G$. Then $g = t_v
\sigma$ with $t_v \in T$ and $\sigma \in G_0$ (since $G=T \rtimes G_0$). Suppose also
that $E\Gamma_i^g =E\Gamma_j$. Then $E\Gamma_i^{\sigma} =E\Gamma_j$ (since $t_v$ fixes
$E\Gamma_i$ setwise). Let $s_i \in S_i$, then $\{\0,s_i\} \in E\Gamma_i$. Since
$E\Gamma_i^{\sigma} = E\Gamma_j$, we have $\{\0,s_i\}^{\sigma} = \{\0, s_i^{\sigma}\}
\in E\Gamma_j$. This implies that $s_i^{\sigma} \in S_j$. So $S_i^{\sigma} \subseteq
S_j$. Similarly if $s_j \in S_j$, then $\{\0,s_j\} \in E\Gamma_j$. Furthermore
$\{\0,s_j\}=e^{\sigma}$ for some edge $e \in E\Gamma_i$. Thus
$e=\{\0,s_j\}^{\sigma^{-1}} = \{\0^{\sigma^{-1}},s_j^{\sigma^{-1}}\} =
\{\0,s_j^{\sigma^{-1}}\} \in E\Gamma_i$ and so $s_j^{\sigma^{-1}} \in S_i$. Hence
$S_j^{\sigma^{-1}} \subseteq S_i \imply S_j \subseteq S_i^{\sigma}$ and we have $S_j
=S_i^{\sigma}$. Thus $G_0$ acts on $\mathcal{S}$. Since $G$ is transitive on $\M{E}$,
for all $i,j \leq k$, there is an element $g = t_v \sigma \in G$ such that
$E\Gamma_i^{g} =E\Gamma_j$, and hence such that $E\Gamma_i^{\sigma} =E\Gamma_j$. Thus
$G_0$ is transitive on $\M{S}$. \epf

By Lemma~\ref{affcayfac} (1), when $G$ is an affine 2-homogeneous permutation group, all
$M$-edge-transitive factors are Cayley graphs. The next result shows that for an
edge-transitive homogeneous factorisation $(M,G,V,\M{E})$ of $K_n=(V,E)$ such that $G$
is affine 2-homogeneous, we may, without loss of generality, assume $G$ to be affine
2-transitive and $M$ to be arc-transitive on each homogeneous factor.

\bp \label{fact2} Let $(M,G,V,\M{E})$ be an edge-transitive homogeneous factorisation of
$K_n=(V,E)$ of index $k$ with factors $\Gamma_i$ $(1 \leq i \leq k)$. Suppose $G=T \rtimes G_0$ is
an affine $2$-homogeneous permutation group on $V=V(a,q)$ and $M =T \rtimes M_0$ is
edge-transitive but not arc-transitive on $\Gamma_i =Cay(V,S_i)$. $($Note that by
Lemma~$\ref{affcayfac}$ $(1)$, $\Gamma_i$ is a Cayley graph.$)$ Then the following hold. \bnum
\item[$1.$] $|V|$ is odd and for each $i$, $S_i=D_i \cup (-D_i)$ where $D_i =s^{M_0}$ for some
$s \in S_i$, and $D_i \neq -D_i$.
\item[$2.$] The map $\phi:V \mapp V$, where $\phi: v \mapp -v $ for all $v \in V$, is such that
$\phi \in \bigcap_{i} (Aut(\Gamma_i))$, $\phi$ centralises $G$, $\phi$ fixes $\0$
$($zero-vector in $V$$)$, and $\brac{G, \phi}$ is an affine $2$-transitive permutation
group. Furthermore, $|\brac{M, \phi}|$ is even and $\brac{M, \phi}$ is arc-transitive on
each $\Gamma_i$, with $\brac{M_0,\phi}$ acting transitively on $S_i$.
\item[$3.$] $(\brac{M,\phi}, \brac{G,\phi}, V, \M{E})$ is an arc-transitive homogeneous
factorisation of $K_n$ $($where $\brac{M,\phi}$ is arc-transitive on each factor$)$.
\enum \ep \bpf (1) \ The result follows immediately from Lemma~\ref{fact1} (2).

(2) \ The map $\phi$ of $V$ as defined above is an automorphism of the additive group of
$V$. So $\phi \in GL(a,q) =Aut(V)$. Suppose $\{x,y\} \in E\Gamma_i$ for some $i$. Then
$y-x \in S_i$. Now $\{x,y\}^{\phi} = \{x^{\phi},y^{\phi}\} = \{-x,-y\}$, and since $S_i
=-S_i$, we have $-(y-x) \in S_i$ and so $\{-x,-y\} \in E\Gamma_i$. As $i$ is arbitrarily
chosen, $\phi \in Aut(\Gamma_i)$ for all $i$. Since $\phi$ acts as $-I$ on $V$ (where
$I$ is the identity element of $GL(a,q)$), $\phi$ centralises $G$.

Observe that $\phi$ fixes the zero-vector in $V$ and $\brac{G,\phi} \leq A\Gamma
L(a,q)$. As $G$ is an affine 2-homogeneous permutation group, it follows that
$\brac{G,\phi}$ is also an affine 2-homogeneous permutation group on $V$ where the
stabiliser of the zero-vector is $\brac{G_0,\phi}$. Since $|\phi|=2$ divides $|\brac{G,
\phi}|$, it follows by Kantor's result on 2-homogeneous but not 2-transitive permutation
groups \cite{Kantor69} that $\brac{G, \phi}$ is 2-transitive.

As shown earlier, $\phi \in Aut(\Gamma_i)$ for all $i$ and so we have $\brac{M,\phi} \leq
Aut(\Gamma_i)$. Moreover, $|\brac{M,\phi}|$ is even since $\phi$ has order 2. Now by
part (1), $S_i =D_i \cup -D_i$  where $D_i =s^{M_0}$ for some $s \in S_i$, and as $\phi:
D \mapp -D$, we have that $\brac{M_0,\phi}$ is transitive on $S_i$ for each $i$. Thus
$\brac{M,\phi}$ is arc-transitive on $\Gamma_i$ by Lemma~\ref{fact1} (1).

(3) \ Since $\brac{M,\phi}$ is normal in $\brac{G,\phi}$, part (3) follows from part (2).
\epf

Thus whenever $(M,G,V,\M{E})$ is an edge-transitive homogeneous factorisation of
$K_{n}=(V,E)$ with $G$ affine 2-homogeneous, then by Proposition~\ref{fact2}, we may
extend both $G$ and $M$ by $\phi$ (defined in Proposition~\ref{fact2}) such that
$\brac{G, \phi}$ is an affine 2-transitive permutation group on $V$, $\brac{M, \phi}$ is
arc-transitive on $\Gamma_i =Cay(V,S_i)$ and $|\brac{M, \phi}|$ even.

So from now on, whenever $G$ is an affine group on $V=V(a,q)$, we will without loss of
generality, assume $(M,G,V,\M{E})$ to be an arc-transitive homogeneous factorisation of
$K_n=(V,E)$ of index $k$ with factors $\Gamma_i=Cay(V,S_i)$ ($1 \leq i \leq k$).
Moreover, we also assume $G$ to be an affine 2-transitive permutation group, and $M$ a
transitive normal subgroup of even order. In particular, $M_0$ is transitive on $S_i$
for all $i$ and $|M_0|$ is even if $q$ odd.

In the next chapter, we will look at the list of affine 2-transitive permutation groups
$G = T \rtimes G_0$ in Theorem~\ref{2_trxclass} (2), and determine all their possible
transitive normal subgroups $M = T \rtimes M_0$ that are not 2-homogeneous. In
particular, since (from the above assumptions) the homogeneous factors are
$M$-arc-transitive Cayley graphs $\Gamma_i =Cay(V,S_i)$, we will determine the
$M_0$-orbits in $V^*$ (which are also elements of the set $\M{S} = \{ S_i \mid 1 \leq i
\leq k \}$) and use them to characterise\footnote{Note that since the group $G=T \rtimes
G_0$ induces isomorphisms between the factors $\Gamma_i =Cay(V,S_i)$ (or rather, it is
$G_0$ that induces the isomorphisms; since by Lemma~\ref{affcayfac} (2), $G_0$ acts
transitively on the set $\M{S}=\{S_i \mid 1 \leq i \leq k \}$), we only need to consider
ONE typical $S_i \in \M{S}$ to describe all the factors.} the homogeneous factors.
Finally, we note that the number of $M_0$-orbits in $V^*$ is precisely the index $k$ of
the homogeneous factorisation.

\section{Examples: two families of arc-transitive homogeneous factorisations} \label{c3:examples}
We first construct two families of arc-transitive homogeneous factorisations
$(M,G,V,\M{E})$ of $K_n =(V,E)$ where the group $G$ is an affine 2-transitive
permutation group. In the examples below, we will use the notations found in the
beginning of Section~\ref{c6:main}. For convenience, we shall repeat some of them here.

Let $V=\F{q}$ be a finite field of order $q=p^R$ and $V^*= V \setminus \{\0\}$. For a
fixed primitive element $\om$ in $V$, we let $\omb$ be the corresponding scalar
multiplication $\omb: x \mapp x\om$ for $x \in V$. Let $\brac{\omb^i}$ be the
multiplicative subgroup of $V^*$ generated by the element $\omb^i$ where $i \geq 1$ and
$i \mid (q-1)$. Since for each $\omb^i$, there is a corresponding $\om^i \in V$, we
shall use $\brac{\om^i}$ to denote the set $\{1,\om^i,\om^{2i}, \ldots, \om^{((q-1)/i
)-i} \} \subseteq V^*$. Also, we let $\al$ denote the Frobenius automorphism of
$\F{p^R}$, that is, $\al: x \mapp x^p$. Finally, note also that $\Gamma L(1,p^R)
=\brac{\omb, \al}$ and hence $A\Gamma L(1,p^R) = T \rtimes \brac{\omb, \al}$, where $T$
is the translation group acting additively on $V =\F{p^R}$.

Recall in Chapter~\ref{c6} where we mentioned that the class of generalised Paley graphs
arise as factors of arc-transitive homogeneous factorisations of complete graphs. We
shall prove this assertion here. Again for convenience, we shall give the definition of a
generalised Paley graph and define what we mean by a generalised Paley partition of a
complete graph.

\bdeff \label{defn:gpaley} \rm{\textbf{(Generalised Paley graph and partition)} \ Let
$V$ and $V^*$ be as defined above. Let $k \geq 2$ be an integer which divides $q-1$ such
that if $p$ is odd, then $\frac{q-1}{k}$ is even. Let $\om$ be a fixed primitive element
in $V$. The graph $\GPq$ is the Cayley graph $Cay(V,S)$ with connecting set $S =
\brac{\om^k} \varsubsetneq V^*$ and is called a \textit{generalised Paley graph} with
respect to $V$.

The \textit{generalised Paley partition} of the complete graph $K_q=(V,E)$ is the
partition $\M{E}_{GP}(q,k) =\{E_1, \ldots, E_k\}$ of the edge set $E$ into $k$ parts
such that $E_i =\{ \{u,v\} \mid v-u \in \om^{i-1}S \}$ for $1 \leq i \leq k$. (The claim
that $\M{E}_{GP}(q,k)$ forms a partition of $E$ is proved in
Proposition~\ref{prop:gpaley}.)} \edeff

\bp \label{prop:gpaley} Let $\Gamma = \GP =Cay(V,S)$ be a generalised Paley graph as
defined in Definition~$\ref{defn:gpaley}$. Let $G =T \rtimes \brac{\omb}$ and $M = T
\rtimes \brac{\omb^k}$. Then the following hold. \bnum
\item $\Gamma$ is an undirected, $M$-arc-transitive graph of valency $\frac{p^R-1}{k}$.
\item $(M,G,V,\M{E}_{GP}(p^R,k))$ is an arc-transitive homogeneous factorisation of
$K_{p^R}=(V,E)$ of index $k$. Furthermore, each factor $\Gamma_i :=(V,E_i)$ is isomorphic
to $\Gamma$. \enum \ep \bpf \ (1) By Lemma~\ref{fel}, $\Gamma =\GP$ admits $M = T
\rtimes \brac{\omb^k} \leq AGL(1,p^R)$ as an arc-transitive subgroup of automorphisms.
To show that $\Gamma$ is undirected, it suffices to show that $-1 =
\om^{\frac{p^R-1}{2}}$ is an element of $S =\brac{\om^k}$. Since $\frac{p^R-1}{k}$ is
even, it follows that $\frac{p^R-1}{2k}$ is an integer. Now $(\om^k)^{\frac{p^R-1}{2k}}
= \om^{\frac{p^R-1}{2}}$ has order 2 and hence is equal to $-1$. Thus $-1 \in S$.
Finally, the valency of $\Gamma$ equals $|S| =\frac{p^R-1}{k}$.

(2) \ It is easy to see that the sets $\om^{i-1} S$, where $1 \leq i \leq k$, are
precisely the $\brac{\omb^k}$-orbits in $V^*$.  Since each $E_i =\{ \{u,v\} \mid v-u \in
\om^{i-1}S \}$ corresponds to $\om^{i-1}S$ (an $\brac{\omb^k}$-orbit) and since
$\om^{i-1}S \neq \om^{j-1}S$ for $i \neq j$, it follows that $E_i \cap E_j = \emptyset$.
Now take an edge $\{u,v\} \in E$. Then $v-u \in \om^{i-1} S$ for some $i$ and it follows
that $\{u,v\} \in E_i$. Hence $\M{E}_{GP}(p^R,k)$ forms a partition of the edge set $E$
of $K_n =(V,E)$.

Note that the $\Gamma_i =(V,E_i)$ are the Cayley graphs $Cay(V, \om^{i-1} S)$ for $1 \leq
i \leq k$. Since $\brac{\omb^k} \normlt \brac{\omb}$ and $\brac{\omb}$ is transitive on
$V^*$ (by Lemma~\ref{gpaley}), it follows by Lemma~\ref{lem:DD} that $\brac{\omb}$
permutes the the sets $\om^{i-1} S$ transitively. Thus for each $i$, there exists an
element of $\brac{\omb}$ which induces an isomorphism between $\Gamma =\Gamma_1
=Cay(V,S)$ and $\Gamma_i =Cay(V, \om^{i-1}S)$. Also, since each $\Gamma_i$ admits $M$ as
a subgroup of automorphisms and $\brac{\omb^k}$ acts transitively on $\om^{i-1} S$ for
each $i$, it follows (by Lemma~\ref{fact1} (1)) that $M$ is arc-transitive on $\Gamma_i$.

Now since $\brac{\omb^k} \normlt \brac{\omb}$, we have $M \normlt G$. Also, $\brac{\omb}$
acts transitively on the sets $\om^{i-1} S$. Now suppose $\sigma \in \brac{\omb}$ maps
$\om^{i-1}S$ to $\om^{j-1}S$ for some $1 \leq i,j \leq k$ and $i \neq j$. Then $\{u,v\}
\in E_i \imply v-u \in \om^{i-1}S$ and so $(v-u)^{\sigma} = v^{\sigma} -u^{\sigma} \in
\om^{j-1}S$. Thus $\{u^{\sigma}, u^{\sigma}\} \in E_j$ and it follows that $E_i^{\sigma}
=E_j$. Since $T$ fixes each $E_i$ setwise, it follows that $G=T \rtimes \brac{\omb}$
acts on $\mathcal{E}_{GP}(p^R,k)$. As $\brac{\omb}$ is transitive on the sets $\om^{i-1}
S$, we have $G$ acting transitively on $\mathcal{E}_{GP}(p^R,k)$. Finally, as $M$ fixes
each $E_i$ set-wise and $G$ is transitive on $\M{E}_{GP}(p^R,k)$, the result on $(M,G,V,
\M{E}_{GP}(p^R,k))$ being a homogeneous factorisation with $M$-arc-transitive factors
follows. \epf


\bdeff \label{defn:tgpaley} \rm{\textbf{(Twisted generalised Paley graph and partition)}
\ Again, let $V$ and $V^*$ be as previously defined and $\om$ be a primitive element in
$V$. Let $R \geq 1$ be even and $p \equiv 3$ (mod 4). Let $h \geq 1$ be an odd integer
such that $2h \mid (p-1)$. Let $S_0 =\brac{\om^{4h}} = \{1, \om^{4h}, \om^{8h}, \ldots,
\om^{q -1 -4h} \} \varsubsetneq V^*$ (observe that $R$ even, $p \equiv 3$ (mod $4$), and
$2h \mid (p-1)$ imply that $8h \mid (q-1)$). The \textit{twisted generalised Paley
graph}, denoted as $TGPaley(q,\frac{q-1}{2h})$, with respect to $V$ is the Cayley graph
$Cay(V, S_0 \cup \om^{3h} S_0)$.

The \textit{twisted generalised Paley partition} of the complete graph $K_q=(V,E)$ is the
partition $\M{E}_{TGP}(q,2h) =\{E_1, \ldots, E_{2h}\}$ of the edge set $E$ into $2h$
parts such that $E_i =\{ \{u,v\} \mid v-u \in \om^{2(i-1)}(S_0 \cup \om^{3h} S_0) \}$ for
$1 \leq i \leq 2h$. }\edeff

\bp \label{prop:tgpaley} Let $\Gamma = TGPaley(p^R,\frac{p^R-1}{2h}) =Cay(V,S_0 \cup
\om^{3h} S_0)$ be a twisted generalised Paley graph as defined in
Definition~$\ref{defn:tgpaley}$. Let $G =T \rtimes \brac{\omb^2, \omb \al}$ and $M = T
\rtimes \brac{\omb^{4h}, \omb^h \al}$. Then the following hold. \bnum
\item $\Gamma$ is an undirected $M$-arc-transitive graph of valency $\frac{p^R-1}{2h}$.
\item $(M,G,V,\M{E}_{TGP}(p^R,2h))$ is an arc-transitive homogeneous factorisation of
$K_{p^R}=(V,E)$ of index $2h$. Furthermore, each factor $\Gamma_i :=(V,E_i)$ is
isomorphic to $\Gamma$. \enum \ep \bpf (1) \ Clearly multiplication by $\om^{4h}$ leaves
both $S_0$ and $\om^{3h}S_0$ invariant. Also, $\omb^h \al$ maps $S_0$ to $(\om^h
S_0)^{\al} =\om^{ph} S_0 = \om^{3h} S_0$ (since $p \equiv 3$ (mod $4$)), and maps
$\om^{3h} S_0$ to $(\om^{4h} S_0)^{\al} = S_0^{\al}=S_0$. Thus, setting $M_0
=\brac{\omb^{4h}, \omb^h \al}$, $\Gamma$ admits $M = T \rtimes \brac{\omb^{4h}, \omb^h
\al}$ as a subgroup of automorphisms. Then since $\brac{\omb^{4h}}$ acts regularly on
$S_0$ and $\om^{3h} S_0$, with $\omb^h \al$ interchanging $S_0$ and $\om^{3h} S_0$, we
have that $M_0$ is transitive on $S_0 \cup \om^{3h} S_0$. Thus by Lemma~\ref{fact1} (1),
$\Gamma$ is $M$-arc-transitive. To show that $\Gamma$ is undirected, it suffices to show
that $-1 =\om^{\frac{p^R-1}{2}} \in S_0 =\brac{\om^{4h}}$. Observe that the conditions
$R$ even, $p \equiv 3$ (mod $4$) and $2h \mid (p-1)$ imply that $8h \mid (p^R-1)$. Thus
$\frac{p^R-1}{8h}$ is an integer and so $(\om^{4h})^{\frac{p^R-1}{8h}} =-1 \in S_0$.
Finally, the valency of $\Gamma$ equals $|S_0 \cup \om^{3h} S_0| =\frac{p^R-1}{2h}$
(since $|S_0|= |\om^{3h} S_0|=\frac{p^R-1}{4h}$).

(2) \ We first show that the sets $\om^{2(i-1)} (S_0 \cup \om^{3h} S_0)$ are $M_0$-orbits
in $V^*$ for all $i=1, \ldots, 2h$. It is easy to see that $\omb^{4h}$ acts transitively
on $\om^{2(i-1)}S_0$ and $\om^{2(i-1) +3h} S_0$. So it suffices for us to show that
$\omb^h \al$ interchanges $\om^{2(i-1)}S_0$ and $\om^{2(i-1) +3h} S_0$. Thus, using the
fact that $p \equiv 3$ (mod 4), we have (recall that $\omb$ acts by multiplication by
$\om$, and $\al: x \mapp x^p$) \be (\om^{2(i-1)}S_0)^{\omb^h \al} &=&
(\om^{2(i-1)})^{\omb^h \al} S_0 \\
            &=& \om^{2p(i-1)+ hp} S_0 \\
            &=& \om^{2p(i-1)+3h} S_0. \ee
Note that since $2h \mid (p-1)$, we have $4h \mid (2p-2)$. So, $2p \equiv 2$ (mod $4h$).
In light of this, it follows that $(\om^{2(i-1)} S_0)^{\omb^h \al} = \om^{2(i-1)+3h}
S_0$. Using a similar argument, we obtain $(\om^{2(i-1) +3h}S_0)^{\omb^h \al} =
\om^{2(i-1)}S_0$ and thus the sets $\om^{2(i-1)} (S_0 \cup \om^{3h} S_0)$ are
$M_0$-orbits in $V^*$. Clearly their union is the whole set $V^*$.

Now by definition, a pair $\{u,v\} \in E_i$ if and only if $v-u$ lies in the $M_0$-orbit
$\om^{2(i-1)}(S_0 \cup \om^{3h}S_0)$. It follows that the graph $\Gamma_i =(V,E_i)$ is
the Cayley graph $Cay(V,\om^{2(i-1)}(S_0 \cup \om^{3h}S_0))$, admitting $M$ as an
arc-transitive subgroup of automorphisms, for $1 \leq i \leq 2h$; and that
$\M{E}_{TGP}(q,2h)$ is a partition of the edge set of $K_{p^R} =(V,E)$.

We will prove in Lemma~\ref{newcon} that $M_0$ is a normal subgroup of $\brac{\omb^2,
\omb \al}$ and that the latter group is transitive on $V^*$. Hence $M \normlt G$ and
$\brac{\omb^2, \omb \al}$ permutes the $M_0$-orbits in $V^*$ transitively. It follows
that $\M{E}_{TGP}(q,2h)$ is $G$-invariant and $G$ permutes the parts (in
$\M{E}_{TGP}(q,2h)$) transitively. (The details of proof are similar to those in the
last paragraph of the proof of Proposition~\ref{prop:gpaley} and are omitted.)

Thus $(M,G,V,\M{E}_{TGP}(p^R,2h))$ is an arc-transitive homogeneous factorisation of
$K_{p^R}$ of index $2h$. Also, as $\brac{\omb^2, \omb \al}$ permutes the $\Gamma_i$
transitively, each $\Gamma_i$ is isomorphic to $\Gamma=\Gamma_1$. \epf

\brk \label{rk:tgpaley} \rm{Note that for $h=1$, the twisted generalised Paley graphs
form the class of self-complementary, arc-transitive graphs, called the $\M{P}^*$-graphs,
which are constructed and studied by Peisert in \cite{Peisert2001}. } \erk

Although the constructions for the generalised Paley graphs and the twisted generalised
Paley graphs look similar at first glance, the two class of graphs are in fact very
different. We will, later in Chapter~\ref{c5}, show that in general (except for one small
case) $\GP \ncong TGPaley(p^R,\frac{p^R-1}{2h})$. Furthermore, we note that in both of
the above cases (the generalised Paley partition and the twisted generalised Paley
partition), the groups $M$ and $G$ are contained in the one-dimensional affine group
$A\Gamma L(1,p^R)$. In fact, the factors that arise from these two partitions are the
only two infinite families of examples known so far that arise from arc-transitive
homogeneous factorisations, where $G \leq A\Gamma L(1,p^R)$.
\chapter{The affine $2$-transitive case} \label{c4}

In this chapter, we will make a case-by-case study of the list of finite 2-transitive
permutation groups $G$ found in Theorem~\ref{2_trxclass} (2) and determine all their
transitive normal subgroups $M$. Then using the information we know about $M$ and $G$,
we determine the structure of the $M$-arc-transitive factors.

Recall from Proposition~\ref{fact2} and the discussion following it, that if $G$ is an
affine $2$-homogeneous group on $V=V(a,q)$ (in an edge-transitive homogeneous
factorisation $(M,G,V,\M{E})$), we may assume $G$ to be 2-transitive, $M$ to have even
order (and hence $M_0$ to have even order if $q$ is odd), and $M$ to act
arc-transitively on each homogeneous factor. Furthermore, we can denote the groups $G$
and $M$ by $G= T \rtimes G_0$ and $M= T \rtimes M_0$ respectively, where $G_0, M_0 \leq
\Gamma L(a,q)$, $M_0 \normlt G_0$ and $T$ is the translation group acting (regularly) on
$V$ by addition. Finally, $G_0$ is transitive on the set of nonzero vectors $V^* =V
\setminus \{\0\}$ while $M_0$ is not transitive on $V^*$. Note that we use the notation
of Theorem~\ref{2_trxclass} (2), and so $q^a =p^R =|V|$.

\section{Cases $2(b)$ - $2(f)$} \label{c4:BtoF}
Here, for cases $2(b)$ to $2(d)$ of Theorem~\ref{2_trxclass}, $G_0$ contains $SL(a,q)$,
$Sp(a,q)$ or $G_2(q)'$ respectively; while for cases $2(e)$ and $2(f)$, $G_0 \cong A_6$,
$A_7$ or $G_0 =SL(2,13)$ respectively. In all these cases $V=V(a,q) =\F{q}^a$ with $a
\geq 2$. We first prove some basic results.

\bl \label{b2f_2} Let $G_0$ be as in cases $2(b)$ to $2(f)$ of
Theorem~$\ref{2_trxclass}$. Define $L$, a normal subgroup of $G_0$, as one of the
following$:$ \bnum
\item[$1.$] $L := SL(a,q)$, $Sp(a,q)$ or $G_2(q)'$ with $a$ and $q$ of $V=V(a,q)$ satisfying:
\bnum
\item[$(a)$] $a \geq 3$ or $q>3$ if $L=SL(a,q)$,
\item[$(b)$] $a$ even and $a \geq 4$ if $L=Sp(a,q)'$,
\item[$(c)$] $q$ is even and $a=6$ if $L=G_2(q)'$,
\enum
\item[$2.$] $L:= A_6$ or $A_7$ and $a=4$, $q=2$ $($that is, $V=V(4,2)$$)$ or,
\item[$3.$] $L:= SL(2,13)$ and $a=6$, $q=3$ $($that is, $V=V(6,3)$$)$.
\enum Suppose $M_0 \normleq G_0$ such that $M_0 \nleq Z(GL(a,q))$. Then $M_0 \geq L$. \el
\bpf Consider case $(1)$ of the above Lemma. Let $Z:= Z(GL(a,q))$. Suppose $\varphi :
\Gamma L(a,q) \mapp P\Gamma L(a,q)$ is the canonical homomorphism (with appropriate
restriction of $\varphi$ to $Z \circ \Gamma Sp(a,q)$ and $Z \times Aut(G_2(q))$ for
cases $2(c)$ and $2(d)$ of Theorem~\ref{2_trxclass} respectively). Moreover, for any $H
\leq \Gamma L (a,q)$, denote $\varphi(H) = HZ/Z$ by $\oline{H}$. Then under $\varphi$,
we have $\oline{L} \leq \oline{G_0} \leq P\Gamma L(a,q)$ where the simple group
$\oline{L}$ is either (i) $PSL(a,q) = (SL(a,q)Z)/Z$, (ii) $PSp(a,q)' = (Sp(a,q)'Z)/Z$ or
(iii) $G_2(q)'$. In particular, if

\vs

\begin{tabular}{l l}
(A) & $\oline{L}=PSL(a,q)$, then $\oline{L} \leq \oline{G_0} \leq P\Gamma L(a,q)$, \\
(B) & $\oline{L}=PSp(a,q)'$, then $\oline{L} \leq \oline{G_0} \leq P\Gamma Sp(a,q)$, \\
(C) & $\oline{L}=G_2(q)'$,  then $\oline{L} \leq \oline{G_0} \leq Aut(G_2(q))$.
\end{tabular}

\vs

For all the above three cases, $\oline{G_0}$ contains $\oline{L}$ as the unique minimal
normal subgroup. Furthermore since $M_0 \normleq G_0$, it follows that $\oline{M_0}
\normleq \oline{G_0}$, where $\oline{M_0}:= M_0Z/Z \neq 1$ since $M_0 \nleq Z$. Hence
$\oline{M_0} \geq \oline{L}$. Then by the Correspondence Theorem applied to the
homomorphism $\varphi$, we have $M_0 Z \geq LZ$. Now $M_0 Z$ contains and normalises $L$
and $M_0$. By the Second Isomorphism Theorem, $M_0Z/M_0 \cong Z/(Z \cap M_0)$ is cyclic.
In particular, the subgroup $M_0 L/M_0 \normleq M_0Z/M_0$ is also cyclic and hence
abelian. Now $M_0 L /M_0 \cong L/(M_0 \cap L)$ and so $L/(M_0 \cap L)$ is abelian.
However, we know that $L$ is a perfect group, i.e. $L'= L$. Hence $M_0 \cap L = L$,
implying that $M_0 \geq L$.

Next consider case (2). Here $G_0 \cong L = A_6$ or $A_7$. Since $A_6$ and $A_7$ are
nonabelian simple groups, we must have $M_0 =L$.

Finally, consider case (3). Here $G_0 = L =SL(2,13)$. The only proper nontrivial normal
subgroup of $G_0$ is $G_0 \cap Z \cong \Z_2$. Since $M_0 \nleq Z$, it follows that $M_0
=L =SL(2,13)$. \epf

Note that the groups $L$ in cases (1), (2) and (3) of Lemma~\ref{b2f_2} are transitive on
$V^*$. (Cases (1) and (2) are well known. For case (3), we refer the reader to \cite[p.
244] {DM_b96}.) Observe also that there are cases involving small values of $a$ and $q$
are not covered in Lemma~\ref{b2f_2}: $SL(2,2)$ and $SL(2,3)$. We will consider them in
the following lemma.

\bl \label{b2f_3small} Let $G_0$ be as in cases $2(b)$ of Theorem~$\ref{2_trxclass}$
with $G_0 \normgeq L$ where $L = SL(2,2)$ or $SL(2,3)$. Suppose $M_0$ is normal in $G_0$
such that $M_0$ is not transitive on the non-zero vectors $V^*$. Then $M_0 \leq
Z(GL(a,q))$ where $(a,q) =(2,2)$ or $(2,3)$. \el \bpf If $L=SL(2,2) \cong S_3$, then $G_0
= L$ and the only proper nontrivial normal subgroup of $SL(2,2)$ has order 3 and acts
regularly on the set of non-zero vectors in $V=V(2,2)$. Thus it follows that $M_0 =1$.

If $L=SL(2,3)$, then there are 2 possibilities for $G_0$, namely $G_0 =GL(2,3)$ or $G_0 =
SL(2,3)$. Every nontrivial normal subgroup $M_0$ of $G_0$ contains $Z=Z(GL(2,3)) \cong
\F{3}^*$, and if $M_0 \nleqslant Z$, then $M_0$ contains a subgroup $Q \cong Q_8$ of
order 8. However, $Q$ is transitive on $V^*$, and $M_0$ is intransitive on $V^*$. Hence
$M_0 \leq Z$. \epf

\bp \label{b2fmain} Let $(M,G,V,\M{E})$ be an arc-transitive homogeneous factorisation of
$K_{q^a}=(V,E)$ of index $k$ with factors $\Gamma_i$ $($where $1 \leq i \leq k)$. Suppose
$G= T \rtimes G_0$ is an affine $2$-transitive permutation group on $V=V(a,q)$ where
$G_0$ is as listed in Theorem~$\ref{2_trxclass}$ $2(b)$ - $2(f)$. Then $M_0 \leq
Z=Z(GL(a,q))$, $\frac{q^a-1}{q-1} \mid k$, $\M{E} =\M{E}_{GP}(q^a,k)$, and each
$\Gamma_i$ is isomorphic to the generalised Paley graph $GPaley(q^a, \frac{q^a-1}{k})$
with $\frac{q^a-1}{k}$ even if $q$ is odd. \ep \bpf By Lemma~\ref{affcayfac} (1),
$\Gamma_i=Cay(V,S_i)$ with $S_i =-S_i$, for $1 \leq i \leq k$. Since $S_i =-S_i$, we
have that $|S_i| = \frac{q^a-1}{k}$ is even if $q$ is odd. Note also that $M_0$ is
intransitive on $V^*$, but is transitive on $S_i$ since $\Gamma$ is $M$-arc-transitive
(see also Lemma~\ref{fact1} (1)). If $G_0$ is as in Lemma~\ref{b2f_2}, then the group
$L$ as defined in that lemma is transitive on $V^*$, and hence $M_0$ cannot contain $L$;
it follows from Lemma~\ref{b2f_2} that $M_0 \leq Z$. Similarly, if $G_0$ is one of the
groups of Lemma~\ref{b2f_3small}, then by that lemma, $M_0 \leq Z$ also. Thus in all
cases, $M_0 \leq Z$.

Now we may identify $V$ with $\F{q^a}$. In this case, $M_0$ is contained in the multiplicative
group of $\F{q^a}^*=V^*$ generated by $\omb$ for a primitive element $\om$ of $\F{q^a}$. (See
Section~\ref{c3:examples} for explanation of notations $\omb$ and $\om$ used here.) Moreover, $M_0$
acts regularly on $S_i$ for each $i$ (since $M_0$, as a multiplicative subgroup of $V^*$, which
acts semiregularly on $V^*$ where $V=\F{q^a}$). Thus it follows that $|M_0| = |S_i| =
\frac{q^a-1}{k}$ and $M_0 =\brac{\omb^k}$. Since $|M_0|$ divides $|Z| =q-1$, we have
$\frac{q^a-1}{k} \mid (q-1)$ and so $\frac{q^a-1}{q-1} \mid k$.

Finally, as $V$ is identified with a finite field $\F{q^a}$, and elements of the set
$\{S_i \mid 1 \leq i \leq k\}$ are orbits of $M_0$ in $V^*$ where $M_0 =\brac{\omb^k}$,
it follows from Proposition~\ref{prop:gpaley} (2) that  $\M{E} =\M{E}_{GP}(q^a,k)$ and
$\Gamma_i \cong GPaley(q^a, \frac{q^a-1}{k})$. \epf

\brk \label{rk:2bf} \rm{Here each $S_i$ is contained in some 1-space of $V$, where we
identify $V$ as the vector space $V(a,q)$, and since $a \geq 2$, the factors $\Gamma_i$
are disconnected (since $\brac{S_i} \neq V$).} \erk

\section{Case $2(g)$} \label{c4:SL}
We first look at the case where $\Gamma L(2,q) \geq G_0 \normgeq SL(2,5)$ and
$V=V(2,q)$, with $q=9,11,19,29$ or $59$. Again we use $V^*$ to denote the set of nonzero
vectors in $V$ and $Z =Z(GL(2,q))$. Clearly, $G_0 \leq N_{\Gamma L(2,q)}(SL(2,5))$. Let
$\varphi: \Gamma L(2,q) \mapp P\Gamma L(2,q)$ be the canonical homomorphism. Then $A_5
\cong PSL(2,5) \normleq \oline{G_0} \leq \oline{N}$ where $\varphi(SL(2,5)) =PSL(2,5)$,
$\oline{G_0}:=\varphi(G_0)$ and $\oline{N}:= \varphi(N_{\Gamma L(2,q)}(SL(2,5))) =
N_{P\Gamma L(2,q)}(PSL(2,5))$. Finally we note that $PSL(2,5) \cong A_5$ is transitive
on the set $P_1(V)$ of 1-spaces of $V$.

\bl \label{sl25} Suppose $P\Gamma L(2,q) \geq \oline{N} \normgeq PSL(2,5)$ where
$q=9,11,19,29$ or $59$. Then $\oline{G_0} = PSL(2,5)$ if $q=11,19,29$ or $59$. If $q=9$,
then $\oline{G_0} =PGL(2,5)$ or $\oline{G_0} =PSL(2,5)$. In all cases, $M_0 \leq Z$ or
$M_0 \geq SL(2,5)$. \el \bpf From \cite[p. 417 (Ex. 7)]{Suzuki1} (which gives the
classification of the maximal subgroups of $PSL(2,q)$ where $q$ is a prime power), we
know that $A_5 \cong PSL(2,5)$ is contained and maximal in $PSL(2,q)$ for $q=9,11,19,29$
or $59$, and there are two conjugacy classes of subgroups $A_5$ in $PSL(2,q)$ that are
interchanged by $PGL(2,q)$. It follows that $A_5$ is self-normalising in $PGL(2,q)$.
Thus $\oline{N} =\oline{G_0}=A_5$ for $q=11,19,29$ or $59$. For $q=9$, we have $P\Gamma
L(2,9) = PGL(2,9) \rtimes \Z_2$ and $\oline{N}= PSL(2,5) \rtimes \Z_2 \cong S_5$, and so
$\oline{G_0} =PSL(2,5)$ or $PGL(2,5)$. Now $M_0$ is normal in $G_0$ and $M_0$ is not
transitive on $V^*$. So under the canonical homomoprhism $\varphi$, we have $\oline{M_0}
\normlt \oline{G_0}$ where $\oline{M_0} =\varphi(M_0)$. Since $PSL(2,5)$ is simple, it
follows that $\oline{M_0} =1$ or $\oline{M_0} \geq PSL(2,5)$, and hence $M_0 \leq Z$
(when $\oline{M_0}=1$) or $M_0 \geq SL(2,5)$ (when $\oline{M_0} \geq PSL(2,5)$). \epf

\brk \label{rk:sl25} \rm{Note that for $q=11$, $SL(2,5)$ is transitive on $V^*$. Thus
for this case, $M_0 \leq Z=Z(GL(2,11))$. } \erk

We now look at the case when $\Gamma L(2,q) \geq G_0 \normgeq SL(2,3)$ where $V=V(2,q)$
and $q=5,7,11$ or $23$. As before, we let $\varphi$ be the canonical homomorphism with
$\varphi(G_0) = \oline{G_0}$ and $\varphi(N_{\Gamma L(2,q)}(SL(2,3))) = \oline{N}$.
Since $q=5,7,11$ and $23$ are primes, we have $PGL(2,q) \geq \oline{N} \geq \oline{G_0}
\normgeq PSL(2,3) \cong A_4$. Note also that $PSL(2,3)$ is transitive on the set of
1-spaces $P_1(V)$ of $V$ only for $q=5$ or $11$. Finally, we let $Z =Z(GL(2,q))$ where
$q =5,7,11$ or $23$.

\bl \label{sl23} Suppose $PGL(2,q) \geq \oline{N} \normgeq PSL(2,3)$ where $q=5,7,11$ or
$23$. Then $\oline{G_0}=PGL(2,3) \cong S_4$ or $\oline{G_0} =PSL(2,3) \cong A_4$ $($only
for $q=5$ or $11)$. Also, $M_0 \leq Z$, $Q_8 \leq M_0 \leq Z \circ Q_8$ $($where $M_0 Z/Z
\cong \Z_2 \times \Z_2)$, or $M_0 \geq SL(2,3)$. \el \bpf By \cite[Theorem
6.26(ii)]{Suzuki1}, it follows that $\oline{N} =PGL(2,3) \cong S_4$, and so we have
either $\oline{G_0} =PGL(2,3)$, or $q=5, 11$ and $\oline{G_0}=PSL(2,3)$ (since
$PSL(2,3)$ is transitive on the $P_1(V)$ only for $q=5$ or $11$). Now $PGL(2,3)$ has
normal subgroups $1$, $\Z_2 \times \Z_2$, $PSL(2,3)$ and $PGL(2,3)$. Thus $\oline{M_0}$
is one of these subgroups. Arguing as in the proof of Lemma~\ref{sl25}, we have $M_0
\leq Z$, or $Q_8 \leq M_0 \leq Z \circ Q_8$ (such that $M_0 Z/Z \cong \Z_2 \times \Z_2$),
or $M_0 \geq SL(2,3)$. \epf

\brk \label{rk:sl23} \rm{Note that for $q=5$, $SL(2,3)$ is transitive on $V^*$. In this
case, we can only have $M_0 \leq Z=Z(GL(2,5))$ or $Q_8 \leq M_0 \leq Z \circ Q_8$. } \erk

We summarise the results obtained in Lemmas~\ref{sl25} and \ref{sl23} and list the
possibilities for $\oline{M_0}$ and $M_0$ (where $M_0 \nleqslant Z$) in
Tables~\ref{tab:TK2} and \ref{tab:TK1}. In addition, using \textsc{Magma}, we compute
(see Appendix \textbf{A3}) the number of $\oline{M_0}$-orbits on the 1-spaces of $V$.

\vspace{0.5cm}

\begin{table}[ht]
\begin{center}
\begin{tabular}{c|c|c|c|c}
\hline
$q$ & $M_0$     & $\oline{M_0}$ & No. of $\oline{M_0}$-orbits  & Length of each $\oline{M_0}$-orbit   \\
    &           &               & in $P_1(V)$                  & in $P_1(V)$                  \\
\hline \hline
$9$, $11$, \  &                     &   &  &       \\
$19$, $29$ \  & $\supseteq SL(2,5)$ &  $PSL(2,5)$ &  $1$  &  $|P_1(V)|$  \\
or $59$    \  &                     &   &  &       \\
\hline
\end{tabular}
\vspace{1mm} \caption{$\oline{M_0}$-orbits on the 1-spaces of $V=V(2,q)$ where $M_0$ is
as in Lemma~\ref{sl25}.} \label{tab:TK2}
\end{center}
\end{table}

\vspace{0.5cm}

\begin{table}[ht]
\begin{center}
\begin{tabular}{c|c|c|c|c}
\hline
$q$ & $M_0$     & $\oline{M_0}$ & No. of $\oline{M_0}$-orbits  & Length of each $\oline{M_0}$-orbit   \\
    &           & & in $P_1(V)$                  & in $P_1(V)$                  \\
\hline \hline
   $5$  &                 &                     &  $3$  &  $2$      \\
   $7$  & $Q_8 \leq M_0 \leq Z \circ Q_8$ &  $\Z_2 \times \Z_2$ &  $2$  &  $4$      \\
   $11$ &                 &                     &  $3$  &  $4$      \\
   $23$ &                 &                     &  $6$  &  $4$       \\
\hline
 $5$  &                     &                       &  $1$  &  $6$       \\
 $7$  & $\supseteq SL(2,3)$ &  $PSL(2,3)$           &  $2$  &  $4$       \\
 $11$ &                     &                       &  $1$  &  $12$       \\
 $23$ &                     &                       &  $2$  &  $12$       \\
\hline
         &      &               &       &                   \\
 $5,7$,  & $\supseteq SL(2,3)$ & $PGL(2,3)$    &  $1$  &  $|P_1(V)|$ \\
 $11,23$ &                     &               &       &                   \\
\hline
\end{tabular}
\vspace{1mm} \caption{$\oline{M_0}$-orbits on the 1-spaces of $V=V(2,q)$ where $M_0$ is
as in Lemma~\ref{sl23}.} \label{tab:TK1}
\end{center}
\end{table}

Although it is clear that the results in Tables~\ref{tab:TK2} and \ref{tab:TK1} do not tell us
much about the number of $M_0$-orbits in $V^*$, they are useful in enabling us to see (almost
directly) if the resulting $M$-arc-transitive homogeneous factors are connected.

\bl \label{lem:Chi}  Let $(M,G,V,\M{E})$ be an arc-transitive homogeneous factorisation
of $K_{q^2}=(V,E)$ of index $k$ with factors $\Gamma_i$ $($where $1 \leq i \leq k)$.
Suppose $G= T \rtimes G_0$ is an affine $2$-transitive permutation group on $V=V(2,q)$,
where $G_0$ is as listed in Theorem~$\ref{2_trxclass}$ $2(g)$. Suppose also that $M_0
\nleqslant Z=Z(GL(2,q))$. Then the $M$-arc-transitive factors $\Gamma_i =Cay(V,S_i)$ are
all connected. \el \bpf The assertion that $\Gamma_i$ is a Cayley graph $Cay(V,S_i)$
follows from Lemma~\ref{affcayfac}. Also, since $\Gamma_i$ is $M$-arc-transitive (where
$M=T \rtimes M_0$), we have that $M_0$ is transitive on $S_i$ (see Lemma~\ref{fact1}
(1)). Since $M_0 \nleqslant Z$, then by Lemmas~\ref{sl25} and \ref{sl23}, $M_0$ is one
of the groups listed in Tables~\ref{tab:TK2} or \ref{tab:TK1}. From Tables~\ref{tab:TK2}
and \ref{tab:TK1}, we see that each orbit of $M_0$ in $P_1(V)$ has length at least 2. In
particular, each $M_0$-orbit in $V^*$, and hence each set $S_i$, must contain elements
from at least two different 1-spaces. As dim$(V) =2$, it follows that each $S_i$ spans
$V$ and so each $\Gamma_i$ is connected. \epf

Again using \textsc{Magma}, we are able to construct \textit{explicitly} all the possibilities for
$M_0$ (where $M_0 \nleqslant Z$, and $M_0$ are as \textit{briefly described} in
Tables~\ref{tab:TK2} and \ref{tab:TK1}) and compute the number of $M_0$-orbits in $V^*$ (Appendix
\textbf{A4}). In doing so, we can immediately determine, for each $M$ and $G$, the index $k$ of
the corresponding arc-transitive homogeneous factorisation $(M,G,V,\M{E})$ where $V=V(2,q)$. The
results are given in Tables~\ref{t25} and \ref{t23}.

Also, for each group $M = T \rtimes M_0$ (where $M_0$ is as in Tables~\ref{t25} and \ref{t23}), we
construct its corresponding $M$-arc-transitive factors $\Gamma_i$ and compute $Aut(\Gamma_i)$
using \textsc{Magma}. A point to note is that as all the $M$-orbital graphs are pairwise
isomorphic, we only need to construct one of them (see Appendix \textbf{A4} and also the footnote
at the end of Section~\ref{c3:affineGP}). The $M$-arc-transitive graphs constructed are denoted by
$G(q^2,k)$, where $k$ is the index of the corresponding arc-transitive homogeneous factorisation
(and also $k$ is the number of $M_0$-orbits in $V^*$).

\subsection{On the graphs $G(q^2,k)$} \label{c4:Gqk}

We shall now give a formal definition\footnote{Although we define $G(q^2,k)$ as a Cayley graph on
$V=V(2,q)$, for computational purposes (using \textsc{Magma}), it is easier to construct it as an
$M$-orbital graph.} for the graphs $G(q^2,k)$ given in Tables~\ref{t25} and \ref{t23}.

\bdeff \label{def:Gqk} \rm{\textbf{(The graph $G(q^2,k)$)} \ \ Let $M =T \rtimes M_0 \leq
A\Gamma L(2,q)$ be an affine permutation group on $V=V(2,q)$ such that $M_0$ is one of
the groups listed in Tables~\ref{t25} and \ref{t23}. Let $v$ be a fixed element of $V^*
=V \setminus \{\0\}$ and let $S = v^{M_0}$. Then $G(q^2,k)$ is defined to be the Cayley
graph $Cay(V,S)$ where $k$ is the number of $M_0$-orbits in $V^*$. } \edeff

\brk \label{rk:Gqk} \rm{Note that each of the groups $M_0$ in Tables~\ref{t25} and \ref{t23} is
such that $M_0 \cap Z(GL(2,q)) \supseteq \{\pm I \}$, where $I$ is the identity matrix. It follows
that if $v \in S$, then $-v \in S$. Thus $S=-S$ and $G(q^2,k)$ is undirected. Also since $G_0$ (as
listed in Theorem~\ref{2_trxclass} $2(g)$) permutes transitively the connecting sets $\{S_i \mid 1
\leq k \}$ of the $M$-arc-transitive homogeneous factors $\Gamma_i =Cay(V,S_i)$, the ``fixed" $v
\in V^*$ in Definition~\ref{def:Gqk} can be arbitrarily chosen.} \erk

In the remainder of this small section, we make some simple observations about the graphs
$G(q^2,k)$, and comment on the information given in Tables~\ref{t25} and \ref{t23}. Furthermore,
we briefly describe how some of the results in Tables~\ref{t25} and \ref{t23} are obtained using
\textsc{Magma} (see Appendix \textbf{A4}\footnote{Note that not all the codes used for performing
our computations were written in Appendix \textbf{A4}. Only function codes and examples that
sufficiently illustrate how the rest of results may be obtained were written.}). Our comments are
listed in ``point-form" below. (In what follow, we have $Z=Z(GL(2,q))$ where $q$ is one of the
values given in Tables~\ref{t25} and \ref{t23}.)

\vspace{0.5cm}

\bnum

\item As noted before, each $G(q^2,k)$ is an orbital graph for the permutation
group $M =T \rtimes M_0$ on $V=V(2,q)$, where $M_0$ is as given in Tables~\ref{t25} and
\ref{t23}. Furthermore by Lemma~\ref{lem:Chi}, all $G(q^2,k)$ are connected.


\item Observe from Tables~\ref{t25} and \ref{t23} that different groups $M_0$ can give
rise to isomorphic $M$-arc-transitive factors $\Gamma_i \cong G(q^2,k)$. (The isomorphisms are
checked using \textsc{Magma}.)


\item It is relatively clear from Tables~\ref{t25} and \ref{t23} that different values
of $q$ and $k$ yield non-isomorphic graphs $G(q^2,k)$, that is, $G(q_1^2, k_1) \ncong
G(q_2^2, k_2)$ whenever $q_1 \neq q_2$ or $k_1 \neq k_2$.


\item We do not include the case where $q =11$ and $M_0 \supseteq SL(2,5)$ in Table~\ref{t25}
(see Remark~\ref{rk:sl25}), and the case where $q=5$ and $M_0 \supseteq SL(2,3)$ in
Table~\ref{t23} (see Remark~\ref{rk:sl23}), the reason being that in both cases, $M_0$ is
transitive on $V^*$ and hence there are no homogeneous factorisations.


\item Note that when $k=2$, we have one of the following: \bnum

\item $q=7$, $SL(2,3) \leq M_0 \leq Z \circ SL(2,3)$ or $M_0 = Z \circ
Q_8$ and $\Gamma_i \cong G(7^2,2)$ (line (3) of Table~\ref{t23}),

\item $q=23$, $M_0 = Z \circ SL(2,3)$ and $\Gamma_i \cong G(23^2,2)$ (line (11) of Table~\ref{t23}) or

\item $q=9$, $SL(2,5) \leq M_0 < (Z \circ SL(2,5)) \cdot \Z_2$ and $\Gamma_i \cong G(9^2,2)$ (line (1) of Table~\ref{t25}).
\enum

In each of the cases above, the homogeneous factors correspond to some arc-transitive
self-complementary graphs. In fact, they are the three exceptional graphs studied in
\cite{Peisert2001}, where they are denoted as \ \textbf{(a)} $G(7^2)$, \ \textbf{(b)}
$G(23^2)$ \ and \ \textbf{(c)} $G(9^2)$ \ respectively (see also Theorem~\ref{MainPei}
and Remark~\ref{rk:pei} in Section~\ref{c1:lit_rev}).

However, in \cite[Lemma 6.6 and 6.7]{Peisert2001}, it was shown that $G(7^2,2)$ (or $G(7^2)$) and
$G(9^2,2)$ (or $G(9^2)$) are isomorphic (using our notation) to the twisted generalised Paley
graphs $TGPaley(49,24)$ and $TGPaley(81,40)$ respectively. (See also Section~\ref{c3:examples} for
details on twisted generalised Paley graphs.) Only the graph $G(23^2,2)$ (or $G(23^2)$) is ``new"
in the sense that it is neither a Paley graph nor a $\mathcal{P}^*$-graph \cite[Lemma
6.8]{Peisert2001}.


\item It was determined using \textsc{Magma} whether or not each of the graphs $G(q^2,k)$
graph is isomorphic to a generalised Paley graph (\textbf{GPaley}) or a twisted
generalised Paley graph (\textbf{TGPaley}). The results are shown under the ``Remarks"
column of Tables~\ref{t25} and \ref{t23}.


\item In line (1) of Table~\ref{t23}, we verified using \textsc{Magma} that
the $\Gamma_i$ are generalised Paley graphs; the fact that the factors are also Hamming
graphs follows from Theorem~\ref{thm:A}. Thus in this case, $Aut(\Gamma_i) = S_5 \ wr \
S_2$. In lines (4) and (5) of Table~\ref{t23}, we again used \textsc{Magma} to verify
that the factors are generalised Paley graphs. It follows that by Lemma~\ref{fel}, the
factors $\Gamma_i$ admit subgroups of automorphisms isomorphic to $T \rtimes
\brac{\omb^{15}, \al}$ (line (4)) and $T \rtimes \brac{\omb^{3}, \al}$ (line (5)). As
the order of the full automorphism groups (computed using \textsc{Magma}) in both cases
equal $|T \rtimes \brac{\omb^{15}, \al}|$ and $|T \rtimes \brac{\omb^{3}, \al}|$, the
results follow. Finally, we note that $Aut(\Gamma_i)$ in each of the remaining lines of
Tables~\ref{t25} and \ref{t23} has the same order as its largest $M=T \rtimes M_0$. As
each $Aut(\Gamma_i)$ contains $M$ as a subgroup of automorphisms ($\Gamma_i$ are
$M$-orbital graphs), it follows that each $Aut(\Gamma_i)$, in the remaining lines of
Tables~\ref{t25} and \ref{t23}, is exactly its largest $M$.


\item By first noting that each $M_0$ in Tables~\ref{t25} and \ref{t23} contains
either $SL(2,5)$, $Q_8$ or $SL(2,3)$ (see Lemmas~\ref{sl25} and \ref{sl23}), the structure of
$M_0$, and hence $M = T \rtimes M_0$, can be easily determined using \textsc{Magma}. In almost all
cases, namely lines (1) to (4) of Table~\ref{t25} and all of Table~\ref{t23} except the second part
of line (6) and line (10), we only need to test for isomorphisms to see if $M_0$ contains
$SL(2,5)$, $Q_8$ or $SL(2,3)$, and compute $Z(M_0)$ (whereby we are able to determine $M_0 \cap
Z$). It turns out that in all these cases, we either have \bnum \item $M_0 < (Z \circ SL(2,5))
\cdot \Z_2$ (see line (1) of Table~\ref{t25}; note that by $X \cdot Y$, we mean an
\textit{extension} of $X$ by $Y$), or \item $M_0 \leq Z \circ B$, where $B= SL(2,5)$, $Q_8$ or
$SL(2,3)$. \enum For the second part of line (6) of Table~\ref{t23}, \textsc{Magma} was used to
show that $M_0$ is isomorphic to $GL(2,3)$. In line (10) of Table~\ref{t23}, we considered the
quotient group of $M_0$ by $SL(2,3)$ ($SL(2,3)$ is normal in $M_0$) and verified using
\textsc{Magma} that it is isomorphic to $\Z_2$; thus $M_0 \cong SL(2,3) \cdot \Z_2$, an extension
of $SL(2,3)$ by $\Z_2$.


\item Finally, for each of the factors $\Gamma_i$ in lines (7) and (8) of Table~\ref{t23}, we
were able to construct (using \textsc{Magma}, see Appendix \textbf{A4: Example A4-4(iv)}) a
$K$-orbital graph $\Sigma$ such that $\Sigma \cong \Gamma_i$, $K \leq A\Gamma L(1,23^2)$ and $|K|
=|Aut(\Gamma_i)| =|M|$, where $M =T \rtimes Q_8$ (line (7)) or $M =T \rtimes (Z \circ Q_8)$ (line
(8)). Thus in both of these cases, $Aut(\Gamma_i)$ is isomorphic to a subgroup of a one-dimensional
affine group $A\Gamma L(1,23^2)$.

We further note that even though $Aut(\Gamma_i) \lnsim A\Gamma L(1,23^2)$, the factors $\Gamma_i$
are neither isomorphic to the generalised Paley graphs (verified using \textsc{Magma}) nor the
twisted generalised Paley graphs. They cannot be isomorphic to the twisted generalised Paley graphs
$TGPaley(23^2, \frac{528}{2h})$, since if they were, then by definition (see
Definition~\ref{defn:tgpaley}) the index $k =2h$ must divide $p-1 = 23-1 =22$. A quick check from
Table~\ref{t23} shows that this is not the case ($k =66$ (line (7)) and $k=6$ (line (8)) do not
divide $22$).

\enum

\vspace{7cm}
\begin{center}

(Turn over for Tables~\ref{t25} and \ref{t23}.)

\end{center}

\newpage

\begin{landscape}
\begin{table}
\begin{center}
{\small
\begin{tabular}
{|l|c|c|c|c|c|c|c|c|c|c|} \hline
               & {\small $G(q^2,k)$ }      &      &         &       &         &       & {\small  No. of }           &     $|S_i|=$            &  $|Y|$,           &                   \\
               & {\small $\cong \Gamma_i$ }&  $q$ & $|V^*|$ & $M_0$ &$|M_0|$  & $|M|$ & {\small $M_0$-orbits }      &   $\frac{|V^*|}{k}$     &  where            & {\small Remarks } \\
               &                           &      &         &       &         &       & {\small  in $V^*$$ = k$ }   &                         &  $Y < Z $         &                   \\
\hline \hline
               &                           &      &         & $SL(2,5)$      & 120 & 9720 &   &   & & {\small $\mathbf{TGPaley}$ } \\
$\mathbf{(1)}$ & {\small $G(9^2,2)$ }      & 9    & 80      & $Y \circ SL(2,5)$      & 240 & 19440&  2 & 40 & $4$ & {\small (see Lemma 6.7,} \\
               &                           &      &         & $(Y \circ SL(2,5)) \cdot \Z_2$      & 480 & 38880 &  &    & & {\small \cite{Peisert2001}) } \\
\hline
               &                           &      &         & $SL(2,5)$      & 120 & 43320 & & & &                               {\small see }             \\
$\mathbf{(2)}$ & {\small $G(19^2,3)$ }     & 19   & 360     & &                      & & 3 & 120  & $6$ & {\small Lemma~\ref{except} }\\
               &                           &      &         & $Y \circ SL(2,5)$ & 360 &129960 & &  & & \\
\hline
               &                           &      &         & $SL(2,5)$ & 120 & 100920 && &                         &       {\small see }               \\
$\mathbf{(3)}$ & {\small $G(29^2,7)$ }     & 29   & 840     & &     & & 7 & 120 & $4$ & {\small Lemma~\ref{except} } \\
               &                           &      &         & $Y \circ SL(2,5)$ & 240 & 201840 &  &  & &\\
\hline
&            & & & & & & & &                                                        & {\small see }                   \\
$\mathbf{(4)}$ & {\small $G(59^2,29)$ }    & 59   & 3480    & $SL(2,5)$ & 120 & 417720 & 29 & 120 & $1$ & {\small Lemma~\ref{except} } \\
\hline

\end{tabular}

\vsl \caption{$SL(2,5) \leq M_0 \normlt G_0 \leq \Gamma L(2,q)$, where $q =9, 19, 29$ or
$59$. If $q =19, 29$ or $59$, then $Aut(\Gamma_i) = T \rtimes (Y \circ SL(2,5))$, where
$Y <Z=Z(GL(2,q))$. If $q=9$, then $Aut(\Gamma_i) = T \rtimes ((Y \circ (SL(2,5)) \cdot
\Z_2)$, where $Y < Z = Z(GL(2,9))$ (the ``$\Z_2$" is due to the group of Frobenius
automorphisms in $\Gamma L(2,9)$).} \label{t25} }
\end{center}
\end{table}
\end{landscape}

\newpage

\begin{landscape}
\begin{table}
{\small
\begin{tabular}{|l|c|c|c|c|c|c|c|c|c|c|}
\hline
             & {\small $G(q^2,k) $ }        &                       &                       &                           &                   &                  & {\small No. of }           & {\small $|S_i|$ }     &                            &   \\
             & {\small $\cong \Gamma_i$}    &   {\small $q$ }       & {\small $|V^*|$ }     &  {\small $M_0$ }          & {\small $|M_0|$ } & {\small $|M|$ }  & {\small $M_0$-orbits}      & {\small $=|V^*|/k$}   & {\small $Aut(\Gamma_i)$}   & {\small Remarks }   \\
             &                              &                       &                       &                           &                   &                  & {\small in $V^*$ $ = k$}   &                       &                            &                     \\

\hline  \hline
$\textbf{(1)}$ & {\small $G(5^2,3)$ }     &  $5$      &  $24$     &  $Q_8$            &  $8$  & $200$     &  $3$  & $8$   &  $S_5 \ wr \ S_2$         &  $\mathbf{GPaley}$    \\
             &                          &           &           &  $Z \circ Q_8$    & $16$  & $400$     &       &       &  (in product action)      &  $H(5,2)$             \\
\hline
$\textbf{(2)}$ & {\small $G(7^2,6)$ }     & 7         & 48        &  $Q_8$            & 8     & 392       & 6     & 8     & $T \rtimes (Q_8 \cdot \Z_3)$   & {\small $\mathbf{TGPaley}$} \\
\hline
             &                          &           &           & $Z \circ Q_8$     & 24    & 1176      &       &       &                                     & {\small $\mathbf{TGPaley}$ } \\
$\textbf{(3)}$ & {\small $G(7^2,2)$ }     & 7         & 48        & $SL(2,3)$         & 24    & 1176      & 2     & 24    & $T \rtimes (Z \circ SL(2,3))$  & {\small (see also } \\
             &                          &           &           & $Z \circ SL(2,3)$ & 72    & 3528      &       &       &                                     &  {\small \cite[Lemma 6.6]{Peisert2001}) } \\
\hline
$\textbf{(4)}$ & {\small $G(11^2,15)$ }   & 11        & 120       & $Q_8$             & 8     & 968       & 15    & 8     & $T \rtimes \brac{\omb^{15}, \al}$  & {\small $\mathbf{GPaley}$ } \\
\hline
$\textbf{(5)}$ & {\small $G(11^2,3)$ }    & 11        & 120       & $Z \circ Q_8$     & 40    & 4840      & 3     & 40    & $T \rtimes \brac{\omb^{3}, \al}$   & {\small $\mathbf{GPaley}$ } \\
\hline
             &                          &           &           & $ SL(2,3)$        & 24    & 2904      &       &       &                                          &                                   \\
$\textbf{(6)}$ & {\small $G(11^2,5)$ }    & 11        & 120       &                   &       &           & 5     & 24    & $T \rtimes GL(2,3)$              & {\small see Lemma~\ref{except} }  \\
             &                          &           &           & $ GL(2,3)$        & 48    & 5808      &       &       &                                          &                                   \\
\hline

$\textbf{(7)}$ & {\small $G(23^2,66)$ }  & 23 & 528 & $Q_8$ & 8 & 4232 & 66 & 8 & $T \rtimes Q_8$ & {\small $Aut(\Gamma_i) \lnsim A\Gamma L(1,23^2)$ } \\

\hline

$\textbf{(8)}$ & {\small $G(23^2,6)$ }   & 23 & 528 & $Z \circ Q_8$ & 88 & 46552& 6 & 88 & $T \rtimes (Z \circ Q_8)$  & {\small $Aut(\Gamma_i) \lnsim A\Gamma L(1,23^2)$ } \\

\hline

$\textbf{(9)}$ & {\small $G(23^2,22)$ }  & 23 & 528 & $SL(2,3)$ & 24 & 12696 &  22 & 24 & $T \rtimes SL(2,3)$         & {\small see Lemma~\ref{except} } \\

\hline

$\textbf{(10)}$ & {\small $G(23^2,11)$ } & 23 & 528 & $ SL(2,3) \cdot \Z_2 $  & 48 & 25392 & 11 & 48 & $T \rtimes (SL(2,3) \cdot \Z_2) $& {\small see Lemma~\ref{except}  } \\
\hline

& &   & &     &&&&&& {\small see Lemma~\ref{except} } \\
$\textbf{(11)}$ & {\small $G(23^2,2)$ } &23 & 528 & $Z \circ SL(2,3)$ & 264 & 139656  & 2 & 264 & $T \rtimes (Z \circ SL(2,3))$         & {\small (see also  } \\
& &   & &     &&&&&& {\small \cite[Lemma 6.8]{Peisert2001}) } \\ \hline

\end{tabular} }

\vspace{2mm} \caption{$M_0 \supseteq SL(2,3)$ or $M_0 \supseteq Q_8$, \ $M_0 \normlt G_0
\leq GL(2,q)$ and $T$ is the translation subgroup of $AGL(2,q)$, where $q =5,7,11$ or
$23$} \label{t23}
\end{table}
\end{landscape}

\newpage

For the rest of the graphs $G(q^2,k)$ which are not isomorphic to the generalised Paley
graphs, the twisted generalised Paley graphs, or the graphs in line (7) and (8) of
Table~\ref{t23} (see ``Remark" columns of Tables~\ref{t25} and \ref{t23}), we are able to
show that their full automorphism groups do not contain any arc-transitive subgroups
which are contained in the one-dimensional affine group. Thus they will not occur in the
one-dimensional case studied in Chapter~\ref{c5}.

\bl \label{except} Let $G(q^2,k)$ be such that $(q,k) = (19,3)$, $(29,7)$ or $(59,29)$ $($see
Table~$\ref{t25}$$)$ or $(q,k)= (11,5)$, $(23,22)$, $(23,11)$ or $(23,2)$ $($see
Table~$\ref{t23}$$)$. Then $Aut(G(q^2,k))$ does not contain an arc-transitive subgroup isomorphic
to a subgroup of $A\Gamma L(1,q^2)$. \el \bpf Let $A:=Aut(G(q^2,k))$ and $A_0$ be its point
stabiliser. Now suppose there exists a subgroup $H \leq A_0$ such that $H \lesssim \Gamma
L(1,q^2)$ $\cong \Z_{q^2-1} \rtimes \Z_2$ and $H$ is transitive on $S_i$. Then since $S_i =-S_i$
is such that $|S_i| =\frac{q^2-1}{k}$ is even, say $|S_i| =2t$, it follows that $H$ must contain a
cyclic subgroup of order at least $t$.

Consider the cases from Table~\ref{t25}, that is, $(q,k) = (19,3)$, $(29,7)$ or $(59,29)$. In all
of these cases, $t= |S_i|/2 = 60$ and so we have one of the following (with $Z=Z(GL(2,q))$; see
also Table~\ref{t25}). \bnum

\item $(q,k)=(19,3)$ and $\Z_{60} \leq  H  \leq A_0 = \Z_6 \circ SL(2,5) < Z \circ SL(2,5)$.

\item $(q,k)=(29,7)$ and $\Z_{60} \leq  H  \leq A_0 = \Z_4 \circ SL(2,5) < Z \circ SL(2,5)$.

\item $(q,k)=(59,29)$ and $\Z_{60} \leq  H  \leq A_0 = SL(2,5)$.

\enum Note that the largest cyclic subgroup in $SL(2,5)$ has order $10$ and $Z  \cap
SL(2,5) \cong \Z_2$. Thus it follows that in all of the above cases, the order of a
cyclic subgroup in $H$ is at most $30$, contradicting the fact that $H \geq \Z_{60}$.

Next, we consider the cases from Table~\ref{t23}, that is, $(q,k)= (11,5)$, $(23,22)$,
$(23,11)$ or $(23,2)$. In these cases, $t = 12,12,24,132$ respectively, and one of the
following holds (also see Table~\ref{t23}): \bnum

\item $(q,k)=(11,5)$ and $\Z_{12} \leq  H  \leq A_0  =GL(2,3)$.

\item $(q,k)=(23,22)$ and $\Z_{12} \leq  H  \leq A_0  = SL(2,3)$.

\item $(q,k)=(23,11)$ and $\Z_{24} \leq  H  \leq A_0  = SL(2,3) \cdot \Z_2$.

\item $(q,k)=(23,2)$ and $\Z_{132} \leq  H  \leq A_0  = Z \circ SL(2,3)$.

\enum Note that largest cyclic subgroup in $SL(2,3)$ has order $6$, and in the four
cases above, one can easily check using \textsc{Magma} that the order of a cyclic
subgroup in $A_0$ is at most $8$ for cases (1), (2) and (3), and at most $66$ for case
(4). It follows that all the above cases contradict the fact that $H \geq \Z_{t}$, and
so the lemma is proved. \epf

We are now ready to prove the following result for case $2(g)$ of
Theorem~\ref{2_trxclass}.

\bp \label{SLL} Let $(M,G,V,\M{E})$ be an arc-transitive homogeneous factorisation of
$K_{q^2}=(V,E)$ of index $k$ with factors $\Gamma_i$ $($$1 \leq i \leq k$$)$. Suppose $G
=T \rtimes G_0$ is an affine $2$-transitive permutation on $V=V(2,q)$ with $G_0$ as
listed in Theorem~$\ref{2_trxclass}$ $2(g)$. Then $\frac{q^2-1}{k}$ is even if $q$ is
odd and one of the following holds. \bnum
\item[$1.$] $M_0 \leq Z=Z(GL(2,q))$, $(q+1) \mid k$, $\M{E} = \M{E}_{GP}(q^2,k)$, and
each $\Gamma_i$ is a $($disconnected$)$ generalised Paley graph $GPaley(q^2,
\frac{q^2-1}{k})$.

\item[$2.$] $M_0$ is as in Tables~$\ref{t25}$ and $\ref{t23}$, the $\Gamma_i$ are connected, and
one of the following holds. \bnum \item[$(a)$] $\Gamma_i \cong G(q^2,k)$ and $(q,k)=(11,5)$,
$(23,22)$, $(23,11)$, $(23,2)$, $(19,3)$, $(29,7)$ or $(59,29)$ $($the graph $G(q^2,k)$ is as in
Definition~$\ref{def:Gqk}$$)$.
\item[$(b)$] $\Gamma_i \cong G(q^2,k)$ where $Aut(\Gamma_i) \lnsim A\Gamma L(1,23^2)$ and
$(q,k)=(23,66)$ or $(23,6)$.
\item[$(c)$] $\Gamma_i \cong GPaley(q^2, \frac{q^2-1}{k})$ and $(q,k) = (5,3)$, $(11,15)$ or
$(11,3)$.
\item[$(d)$] $\Gamma_i \cong TGPaley(q^2, \frac{q^2-1}{k})$ and $(q,k) = (9,2)$, $(7,6)$ or $(7,2)$.
\enum \enum \ep \bpf By Lemma~\ref{affcayfac} (1), $\Gamma_i=Cay(V,S_i)$ for all $1
\leq i \leq k$. Also each $\Gamma_i$ is undirected and so $S_i =-S_i$. Hence $|S_i| =
\frac{q^a-1}{k}$ is even if $q$ is odd.

Now $M_0$ is one of the groups given by Lemmas~\ref{sl25} and \ref{sl23}. Suppose first
that $M_0 \leq Z$. Then as in the proof of Proposition~\ref{b2fmain} and also by
Remark~\ref{rk:2bf}, the factors $\Gamma_i$ are disconnected generalised Paley graphs,
$\M{E}$ is the corresponding generalised Paley partition, and so (1) holds. From now on,
we will assume that $M_0 \nleqslant Z$.

Then $M_0$ is one of the groups listed in Tables~\ref{t25} and \ref{t23}, with
corresponding graphs $G(q^2,k)$ as in Defintion~\ref{def:Gqk}. Furthermore, by
Lemma~\ref{except} as well as the results shown under the ``Remark" columns of
Tables~\ref{t25} and \ref{t23}, the homogeneous factors $\Gamma_i$ are as listed in
$2(a)$ - $2(d)$ above. Finally, by Lemma~\ref{lem:Chi}, all the factors $\Gamma_i$ in
$(2)$ are connected. \epf

\brk \label{rk:partition} \rm{For a finite set $V$, we say that two partitions $\M{E}$ and $\M{L}$
of $V^{\{2\}}$ are \textit{equivalent} if there exists $g \in$ Sym$(V)$ such that $\M{E}^g =\M{L}$
and $E \in \M{E} \imply E^{g} \in \M{L}$. Now, even though the factors in $2(c)$ and $2(d)$ of
Proposition~\ref{SLL} are generalised Paley graphs and twisted generalised Paley graphs
respectively, we have not decided if their corresponding partitions $\M{E}$ are equivalent to the
generalised Paley partitions or the twisted generalised Paley partitions. (Note that we first have
to identify $V=V(2,q)$ with the finite field $\F{q^2}$.) Although we are primarily interested in
identifying the structure of the factors that arise from an arc-transitive homogeneous
factorisation of a complete graph, it would be interesting to know if isomorphic homogeneous
factors (from two homogeneous factorisations $(M,G,V,\M{E})$ and $(M',G',V',\M{E}')$ such that
there is a bijection between $V$ and $V'$) can give rise to non-equivalent partitions.} \erk

\section{Case $2(h)$} \label{c4:extraspecial}
Here we are looking at the scenario where $GL(4,3) \geq G_0 \normgeq \mathbf{E}$ and
$\mathbf{E}$ is an extraspecial subgroup of order $2^5$. Also $V=V(4,3)$ and $V^*$
denotes the set of nonzero vectors in $V$. Now $G_0 \leq N_{GL(4,3)}(\mathbf{E})$, and
under the canonical homomorphism $\varphi: GL(4,3) \mapp PGL(4,3)$, we have $PGL(4,3)
\geq \oline{N} \geq \oline{G_0} \geq \oline{\mathbf{E}}$ where $\oline{N}:= \varphi
(N_{GL(4,3)}(\mathbf{E})) = N_{PGL(4,3)}(\oline{\mathbf{E}})$, $\oline{\mathbf{E}}:=
\mathbf{E}Z/Z \cong \Z_2^4$ and $Z=Z(GL(4,3))$.

\bl \label{esp1} $\oline{N} \cong \Z_2^4 \rtimes S_5$ and $\oline{G_0} \gtrsim \Z_2^4 \rtimes
\Z_5$. \el \bpf The group $\oline{\mathbf{E}}$ has 5 orbits of length 8 in the set $P_1(V)$ of
1-spaces in $V =\F{3}^4$. By Theorem~\ref{2_trxclass} $2(h)$, $\oline{G_0}$, and hence also
$\oline{N}$, act transitively on $P_1(V)$, and so $\oline{G_0}$ and $\oline{N}$ are transitive on
the 5 orbits of $\oline{\mathbf{E}}$ in $P_1(V)$. Since $\oline{G_0}$ contains the elementary
abelian subgroup $\oline{\mathbf{E}}$ of order 16, it follows that $2^4 \cdot 5$ divides
$|\oline{G_0}|$. Now $\oline{\mathbf{E}} \leq PSL(4,3)$ (checked using \textsc{Magma}) and
$|PSL(4,3)| = 2^7 \cdot 3^6 \cdot 5 \cdot 13$. It follows that if $K$ is a maximal subgroup of
$PSL(4,3)$ containing $\oline{N} \cap PSL(4,3)$, then $|PSL(4,3):K|$ is not divisible by 5 and $K$
contains a subgroup isomorphic to $\Z_2^4$. From the \Atlas \cite{Atlas}, we find that $K=PSU(4,2)
\rtimes \Z_2 \cong PGU(4,2)$. Then, by checking through the list of maximal subgroups of
$PGU(4,2)$ in the \Atlas \cite{Atlas}, only the subgroup isomorphic to $\Z_2^4 \rtimes S_5$
contains $\oline{\mathbf{E}} \cong \Z_2^4$. Since $\Z_2^4$ is obviously normal in $\Z_2^4 \rtimes
S_5$, we have that $\oline{N} \cap PSL(4,3)\cong \Z_2^4 \rtimes S_5$. In fact, this is the full
normaliser of $\oline{\mathbf{E}} \cong \Z_2^4$ in $PGL(4,3)$ (see Remark~\ref{esp1_magma} below)
and so is equal to $\oline{N}$. Finally since $2^4.5 \mid |\oline{G_0}|$ and $\oline{G_0} \leq
\oline{N} \cong \Z_2^4 \rtimes S_5$, it follows that $\oline{G_0}$ must contain a subgroup
isomorphic to $\Z_2^4 \rtimes \Z_5$. \epf

\brk \label{esp1_magma} \rm{We could use \Magm to arrive at the above result and this is
the first way the result was obtained. First compute the normaliser $\oline{N}$, of the
elementary abelian group $\oline{\mathbf{E}}$, of order $2^4$ in $PGL(4,3)$. Now
$\oline{N}$ has order $|\oline{N}| = 1920$. Clearly $\oline{N} \normgt
\oline{\mathbf{E}}$. Upon inspecting the maximal subgroups (there are 5 conjugacy
classes of maximal subgroups) of $\oline{N}$, only one (isomorphic to $S_5$) does not
contain $\oline{\mathbf{E}}$ as a subgroup. For such a subgroup $S_5$, we have $S_5 \cap
\oline{\mathbf{E}} =1$. Since $|\oline{\mathbf{E}}.S_5| =1920 = |\oline{N}|$, the result
follows.} \erk

\bc \label{esp2} The full pre-image of $\oline{N} \cong \Z_2^4 \rtimes S_5$ in $GL(4,3)$
is $N_{GL(4,3)}(\mathbf{E}) = \mathbf{E}H$ where $H >Z$ and $H/Z \cong S_5$. \ec \bpf
From Lemma~\ref{esp1}, $\oline{N} \cong \Z_2^4 \rtimes S_5$. The pre-image of
$\oline{\mathbf{E}} \cong \Z_2^4$ is $\mathbf{E}$. Let $H$ be the full pre-image of
$S_5$ with respect to $\varphi$. It follows that, since $\oline{\mathbf{E}} \cap S_5
=1$, we have $\mathbf{E} \cap H =Z$. This implies that $H > Z$ and $H/Z \cong S_5$. Thus
the full pre-image of $\oline{N} = \Z_2^4 \rtimes S_5$ is $N_{GL(4,3)}(\mathbf{E}) =
\mathbf{E}H$. \epf

Since $\mathbf{E} < G_0 \leq N_{GL(4,3)}(\mathbf{E}) =\mathbf{E}H$
(Corollary~\ref{esp2}), we have $G_0 =\mathbf{E}L$ where $L \leq H$, and $\oline{G_0} =
\oline{\mathbf{E}} \rtimes \oline{L}$ where $\oline{L} = LZ/Z$ (note that $L$ does not
necessarily contain the scalars $Z$). By Lemma~\ref{esp1}, $G_0$ (and $\oline{G_0}$)
have order divisible by 5, so $|L|$ and $|\oline{L}|$ are divisible by 5. Now $\oline{L}
= LZ/Z \lesssim S_5$ and since $|\oline{L}|$ is divisible by 5, it follows that
$\oline{L}$ is isomorphic to one of $\Z_5$, $D_{10}$, $F_{20}$ (the Frobenius group,
$\Z_5 \rtimes \Z_4$), $A_5$ or $S_5$. Using \textsc{Magma}, we can easily check that all
these (possible) subgroups $\oline{G_0}$ are transitive on the set $P_1(V)$ of 1-spaces.

Now the pre-image of $\oline{G_0}$ contains $\mathbf{E}$ which also contains the scalars
$Z$. Moreover, since $\oline{\mathbf{E}}$ has 5 orbits of length $8$ in $P_1(V)$, it
follows that $\mathbf{E}$ then has $5$ orbits of length 16 in $V^*$. Also, since $L$ (or
$\oline{L}$) permutes transitively the set of $\oline{\mathbf{E}}$-orbits in $P_1(V)$,
it follows that $G_0= \mathbf{E}L$ (the pre-image of $\oline{G_0}$ under $\varphi$) is
transitive on $V^*$ (note that the extraspecial subgroup $\mathbf{E}$ of $G_0$ contains
$Z$, and so even if $L$ does not contain the scalars $Z$, $G_0$ is \textit{still}
transitive on $V^*$ \footnote{In fact, using \textsc{Magma}, we can show that there
exist a group ($\mathbf{E} \rtimes \Z_5$) whereby $L$ does not contain the scalars.} ).
Now $M_0 \normlt G_0$ and $M_0$ is not transitive on $V^*$. Under the canonical
homomorphism $\varphi$, we have $\oline{M_0} \normlt \oline{G_0}$ where $\varphi(M_0) =
\oline{M_0}$. The following result determines the possibilities for $\oline{M_0}$.

\bl \label{esp4} Suppose $\mathbf{E} \leq G_0 \leq N_{GL(4,3)}(E)$. Then $G_0
=\mathbf{E}L$ where $L \leq H$ and $\oline{L} = LZ/Z = \Z_5$, $D_{10}$, $F_{20}$, $A_5$
or $S_5$. Further, $\oline{M_0}=1$ or $\oline{M_0}=\oline{\mathbf{E}}$, and $M_0 =Z$ or
$M_0 =\mathbf{E}$ respectively. \el \bpf By the discussion in the above paragraph, we
have that $G_0 =\mathbf{E}L$ where $L \leq H$, and the possibilities for $\oline{L}$ are
as listed above.

Suppose $K \normlt \oline{G_0}$ and $5$ divides $|K|$. Then $K$ must contain
$\oline{\mathbf{E}} \rtimes \Z_5$ (otherwise we have $K \cap \oline{\mathbf{E}} =1$ and
so $[K,\oline{\mathbf{E}}]=1$ whereas elements of order 5 do not centralise
$\oline{\mathbf{E}}$). Now if 5 divides $|\oline{M_0}|$, then by the previous argument,
we have $\oline{M_0} \supseteq \oline{\mathbf{E}} \rtimes \Z_5$. However, as mentioned
earlier, the pre-image of $\oline{\mathbf{E}} \rtimes \Z_5$ is transitive on $V^*$,
which is a contradiction since $M_0$ is intransitive on $V^*$. Hence $5 \nmid
|\oline{M_0}|$ and so we have $\oline{M_0} \leq \oline{\mathbf{E}}$. It follows that
$\oline{M_0} =1$ or $\oline{\mathbf{E}}$.

If $\oline{M_0}=1$, then the full pre-image of $\oline{M_0}$ is $Z$ (since we require
that $|M_0|$ is even, so we must have $M_0 =Z$). If $\oline{M_0}=\oline{\mathbf{E}}$,
then $M_0$ is either $\mathbf{E}$ or $M_0 \cong \oline{\mathbf{E}} =\Z_2^4 <\mathbf{E}$.
The latter case is not possible since $\mathbf{E}$ does not contain an elementary abelian
subgroup of order $16$. So $M_0 =\mathbf{E}$. \epf

Now we are ready to state our main result for this case. (Identification of the
isomorphism type of the factors $\Gamma_i$ in the case where $M_0 = \mathbf{E}$ will be
done by referring to the main result in Chapter~\ref{c7}: Theorem~\ref{hamm}.)

\bp \label{espmain} Let $(M,G,V,\M{E})$ be an arc-transitive homogeneous factorisation of
$K_{81}=(V,E)$ of index $k$ with factors $\Gamma_i$. Suppose $G=T \rtimes G_0$ is an
affine $2$-transitive permutation group on $V$ such that $G_0$ is as in case $2(h)$ of
Theorem~$\ref{2_trxclass}$ with $V=V(4,3)$. Then $\Gamma_i$ are Cayley graphs
$Cay(V,S_i)$ for all $1 \leq i \leq k$ and $M = T \rtimes M_0$ where either $M_0 =
Z=Z(GL(4,3))$ or $M_0 =\mathbf{E}$. Furthermore, we have one of the following. \bnum
\item[$1.$] $M_0 = Z$, $k=40$, $\M{E} =\M{E}_{GP}(81,40)$, and each $\Gamma_i$ is a
$($disconnected$)$ generalised Paley graph $GPaley(3^4,2)$, or
\item[$2.$] $M_0 =\mathbf{E}$, $k =5$ and each homogeneous factor $\Gamma_i$ is isomorphic to the
Hamming graph $H(9,2)$ $($see Definition~$\ref{def:Hamming}$$)$. \enum \ep \bpf By
Lemma~\ref{affcayfac} (1), $\Gamma_i=Cay(V,S_i)$ for all $1 \leq i \leq k$. Also, the
possibilities for $M_0$ are given in Lemma~\ref{esp4}. Suppose $M_0 =Z$, then as in the
proof of Proposition~\ref{b2fmain} and also by Remark~\ref{rk:2bf}, the factors
$\Gamma_i$ are disconnected generalised Paley graphs. Moreover, as $M_0 =Z \cong \Z_2$
acts regularly on $S_i$ for each $i$, it follows that $|S_i| =2$ and so $k=40$. Thus (1)
holds.

Now suppose $M_0=\mathbf{E}$. Then since $\mathbf{E}$ has 5 orbits of equal length 16 in
$V^*$ (with $G_0$ leaving the orbits of $\mathbf{E}$ invariant and permuting them
transitively; see also Lemma~\ref{lem:DD}), it follows that the homogeneous
factorisation has index $k=5$. Furthermore, by Theorem~\ref{hamm}, we have that
$\Gamma_i \cong H(9,2)$. \epf

\brk \label{rk:extra} \rm{Note that for part (1) of the above Proposition, we can easily
work out the structure of the connected components of the $\Gamma_i$. As the $\Gamma_i$
are generalised Paley graphs, we can identify $V$ as the finite field $\F{3^4}$. Here
$M_0 \cong \Z_2$ (seen as a multiplicative subgroup of $\F{3^4}^*$) acts semiregularly
by field multiplication on $V^* =\F{3^4}^*$ and regularly on $S_i$. Suppose $\Gamma_1
=Cay(V,S_1)$ is such that $1 \in S_1$. Then $S_1 =\{\pm1\}$ and it follows that the
$\F{3}$-span of $S_1$ equals $\F{3}$. Thus the connected components of $\Gamma_1$, and
hence of all the $\Gamma_i$, are isomorphic to $K_3$. } \erk

\section{Summing up the cases} \label{c4:SUM}
Finally, we sum up the results from Sections~\ref{c4:BtoF} to \ref{c4:extraspecial} in
the following theorem.

\bt \label{thm:mainc4} Let $(M,G,V,\M{E})$ be an arc-transitive homogeneous
factorisation of $K_{q^a}=(V,E)$ of index $k$ with factors $\Gamma_i$. Suppose $G=T
\rtimes G_0$ is an affine $2$-transitive permutation group on $V=V(a,q)$ such that $G_0$
is in one of the cases $2(b) - 2(h)$ of Theorem~$\ref{2_trxclass}$. Then $\Gamma_i
=Cay(V,S_i)$, where $|S_i| =\frac{q^a-1}{k}$ and $|S_i|$ is even if $q$ is odd, and $M =
T \rtimes M_0$ is such that precisely one of the following holds. \bnum

\item $M_0 \leq Z=Z(GL(a,q))$, $(q+1) \mid k$, $\M{E} =\M{E}_{GP}(q^a,k)$, and each
$\Gamma_i$ is a $($disconnected$)$ generalised Paley graph $GPaley(q^a,
\frac{q^a-1}{k})$.

\item $M_0$ is as in Tables~$\ref{t25}$ and $\ref{t23}$, the $\Gamma_i$ are connected, and
one of the following holds. \bnum \item $\Gamma_i \cong G(q^2,k)$ and $(q,k)=(11,5)$, $(23,22)$,
$(23,11)$, $(23,2)$, $(19,3)$, $(29,7)$ or $(59,29)$.
\item $\Gamma_i \cong G(q^2,k)$ where $Aut(\Gamma_i) \lnsim A\Gamma L(1,23^2)$ and
$(q,k)=(23,66)$ or $(23,6)$.
\item $\Gamma_i \cong GPaley(q^2, \frac{q^2-1}{k})$ and $(q,k) = (5,3)$, $(11,15)$ or
$(11,3)$.
\item $\Gamma_i \cong TGPaley(q^2, \frac{q^2-1}{k})$ and $(q,k) = (9,2)$, $(7,6)$ or $(7,2)$.
\enum

\item $M_0 =\mathbf{E}$ and $k =5$ and the homogeneous factors $\Gamma_i$ are all
isomorphic to the Hamming graph $H(9,2)$.

\enum \et \bpf By Lemma~\ref{affcayfac} (1), $\Gamma_i =Cay(V,S_i)$ and $|S_i| =\frac{q^a-1}{k}$
for all $1 \leq i \leq k$. Furthermore, since we require each $\Gamma_i$ to be undirected, we have
$S_i =-S_i$ and so $|S_i|$ is even if $q$ is odd. The fact that one of $(1)$ to $(3)$ holds now
follows from Propositions~\ref{b2fmain}, \ref{SLL} and \ref{espmain}. Finally to see that the
factors arising from $(M,G,V, \M{E})$ belong to exactly one of the above cases, we note the
following.

The graphs/factors from case (1) are not connected, they are not isomorphic to those connected
ones found in cases (2) and (3). Also, as mentioned in Section~\ref{c4:Gqk}, the graphs for
different parts within case (2) are not isomorphic as they correspond to different choices of
$(q,k)$. Finally, the factors isomorphic to $H(9,2)$ in case (3) are not isomorphic to those found
in (2); since for case (3), $(q,k) =(9,5)$ but there is no graph in (2) which has the same
parameters. \epf

\chapter{Hamming graphs as arc-transitive homogeneous factors} \label{c7}

We saw in Section~\ref{c4:Gqk} that there is a homogeneous factorisation of $K_{25}$ of index $3$
in which each of the factors is a Hamming graph $H(5,2)$. In Chapter~\ref{c3}, we showed that the
generalised Paley graphs, of which some are isomorphic to the Hamming graphs (see
Theorem~\ref{thm:A}), also arise as factors of arc-transitive homogeneous factorisations of
complete graphs. However, of these two examples, only the former one occurs with factors arising
from a homogeneous factorisation $(M,G,V,\M{E})$ where $G$ is not a one-dimensional affine group
(even though $H(5,2)$ is isomorphic to the generalised Paley graph $GPaley(25, 8)$). Such a
situation also arises for the complete graph $K_{81}$.

Recall from Section~\ref{c4:extraspecial} that there exists an arc-transitive homogeneous
factorisation $(M,G,V,\M{E})$ of $K_{81}=(V,E)$ of index $5$ where $M = T \rtimes \mathbf{E}$, and
$\mathbf{E}$ is an extraspecial subgroup of $GL(4,3)$ of order $2^5$. In particular, in
Proposition~\ref{espmain}, it is claimed (by using Theorem~\ref{hamm} in this chapter) that the
homogeneous factors that arise in this case are isomorphic to the Hamming graphs $H(9,2)$. The aim
of this chapter is to show that the claim is true.

\section{Some representation theory} \label{c7:rep}
As in Section~\ref{c4:extraspecial}, let $M = T \rtimes \mathbf{E}$ be a subgroup of $AGL(4,3)$
acting on $V=V(4,3)$, where $\mathbf{E}$ is an extraspecial subgroup of $GL(4,3)$ of order $2^5$.
Let $\Gamma_i=Cay(V,S_i)$ (for $1 \leq i \leq 5$) be the $M$-arc-transitive factors that arise from
the arc-transitive homogeneous factorisation of $K_{81}$ of index $5$ in
Section~\ref{c4:extraspecial}. We will show that the $\Gamma_i$ are all isomorphic to the Hamming
graph $H(9,2)$.

It is well known that the Hamming graph $H(9,2)$ is isomorphic to a Cartesian product of
complete graphs $K_9 \ \square \ K_9$ (see also Definition~\ref{def:Hamming}). By a
Cartesian product $\Gamma \ \square \ \Gamma'$ of graphs $\Gamma =(V,E)$ and
$\Gamma'=(V', E')$, we mean the graph with vertex set $V \times V'$ and edge set $\{
\{(x,y), (z,y)\} \mid \{x,z\} \in E \}$ $\cup$ $\{ \{(i,j), (i,k)\} \mid \{j,k \} \in E'
\}$. The following result is taken from \cite{HMP2001}.

\bp \label{prop:prod} {\rm \cite[Proposition 3.3]{HMP2001}} \ If $\Gamma =Cay(V,S)$ and
$\Gamma'=Cay(V', S')$. Then $\Gamma \ \square \ \Gamma'$ is the Cayley graph on $V \times V'$ with
generating set $(S \times 1_{V'}) \cup (1_V \times S')$. \ep

Note that the group $M = T \rtimes \mathbf{E}$ contains $L:= T \rtimes Q_8$ as a subgroup.
Moreover, since $Q_8 \normlt \mathbf{E}$, it follows that $L$ is normal in $M$. We first show (in
Lemma~\ref{irred_rep}) that $Q_8$ (as a subgroup of $GL(4,3)$ acting on $V =V(4,3)$) does not have
an irreducible $\F{3}$-representation of degree 4. To do that, we need to recall some basic
concepts from representation theory (for more details of which, see for example \cite{Huppert_b98}
and \cite{Isaacs_b76}).

For any arbitrary field $\F{}$ and any group $G$, an $\F{} G$-module $U$ is said to be
\textit{irreducible} if it is non-zero and it has no $\F{} G$-submodules apart from $\{\0\}$ and
$U$ (otherwise it is called \textit{reducible}). A representation of $G$ over $\F{}$ is a
homomorphism $\Psi$ from $G$ to $GL(n,\F{})$, that is, $\Psi: G \mapp GL(n,\F{})$ (we often refer
to $\Psi$ as an $\F{} G$-representation or an $\F{}$-representation of the group $G$). The
\textit{degree} of the representation $\Psi$ is the integer $n$. A representation $\Psi$ is said
to be \textit{irreducible} if the corresponding $\F{} G$-module is irreducible. In other words, a
representation $\Psi$ is irreducible if its corresponding $\F{} G$-module $U$ is non-zero and it
has no $\F{} G$-submodules apart from $\{\0\}$ and $U$ that is invariant under $\Psi(G)$. The
\textit{$\FF$-character} or \textit{character} $\chi$ of $G$ afforded by the $\FF
G$-representation $\Psi$ is the function given by $\chi(g) := tr(\Psi(g))$ where $tr(\Psi(g))$ is
the trace of the matrix $\Psi(g) \in GL(n,\FF)$ for all $g \in G$. If the $\FF G$-representation
$\Psi$ is irreducible, then we call the corresponding $\FF$-character afforded by $\Psi$ an
\textit{irreducible character} of $G$ over $\FF$.

Suppose now $\FF \subseteq \EE$ and $\chi$ is an $\EE$-character of a group $G$ afforded
by some $\EE G$-representation. We write $\FF(\chi)$ to denote the subfield of $\EE$
generated by $\FF$ and the character values $\chi(g)$ for $g \in G$. In fact,
$\FF(\chi)$ is the smallest extension of $\FF$ containing all $\chi (g)$, where $g \in
G$. Furthermore, if $\FF \subseteq \EE$ and $\Psi$ is an $\FF$-representation of $G$,
then we may view $\Psi$ as an $\EE$-representation of $G$. As such, we denote it by
$\Psi^\EE$.

Let $\Psi$ be an $\FF$-representation of a group $G$. Then $\Psi$ is \textit{absolutely
irreducible} if $\Psi^\EE$ is irreducible for every field $\EE \supseteq \FF$. The field $\FF$ is
a \textit{splitting field} for $G$ if every irreducible $\FF$-representation of $G$ is absolutely
irreducible.

\bl {\rm \cite[Corollary 9.15]{Isaacs_b76}} \label{split} \ Let $G$ have exponent $m$ and assume
the polynomial $x^m -1$ splits into linear factors in the field $\FF$. If $\FF$ has prime
characteristic, then it is a splitting field for $G$. $($Note that the \textit{exponent} of $G$ is
the least positive integer $m$ such that $g^m=1$ for all $g \in G$.$)$ \el

\bl \label{splita} The finite field $\F{9}$ is a splitting field for $Q_8$. \el \bpf The group
$Q_8$ has exponent $4$ and the finite field $\F{9}$ has characteristic $3$. Since $x^4 -1 =
(x-1)(x+1)(x-\al^2) (x+\al^2)$ in $\F{9}$ where $\al$ is a primitive element of $\F{9}$, the
result then follows from Lemma~\ref{split} \epf

\bl \label{splitb} {\rm \cite[Lemma 34.3 (1)]{Ashbacher_b86}} \ Let $G$ be a finite group and $\FF$
be a splitting field for $G$ such that the characteristic does not divide $|G|$. Then the number
of equivalence classes of irreducible $\FF G$-representations is equal to the number of conjugacy
classes of $G$. \el

By Lemma~\ref{splita}, $\F{9}$ is a splitting field for $Q_8$. Since the characteristic of $\F{9}$
does not divide $|Q_8|=8$, by Lemma~\ref{splitb}, the number of non-isomorphic irreducible
$\F{9}Q_8$-representations is equal to the number of conjugacy classes of $Q_8$. Now $Q_8
=\brac{a,b \mid a^4 =1, \ b^2=a^2, \ b^{-1}ab=a^{-1}}$ has five conjugacy classes. Thus we have
exactly five non-isomorphic irreducible representations of $Q_8$ over $\F{9}$. Let $\Psi_1,
\Psi_2, \ldots, \Psi_5$ be the five non-isomorphic irreducible $\F{9}$-representations of $Q_8$
and $\chi_1, \chi_2, \ldots, \chi_5$ be their respective irreducible characters. Furthermore (using
the fact that $\F{9}$ is a splitting field for $Q_8$ with characteristic not dividing $|Q_8|$)
from \cite[Lemma 34.2 and 34.4]{Ashbacher_b86}, we have $|Q_8| = 8 = \sum_{i=1}^{5}
\chi_i(1_{Q_8})^2$ where $\chi_i(1_{Q_8})$ is the degree of the irreducible $\F{9}
Q_8$-representation $\Psi_i$. Thus it is easy to see that $Q_8$ has four non-isomorphic
(irreducible) representations of degree 1 and one (unique) irreducible representation of degree 2
over $\F{9}$.  The four representations of degree 1 are:
\[\begin{array}{ll} \Psi_1 : \left\{ \begin{array}{l}
                        a \mapsto 1\\ b \mapsto 1
                        ,\end{array} \right.&
\Psi_2 : \left\{ \begin{array}{l}
                        a \mapsto 1\\ b \mapsto -1
                        ,\end{array} \right.\\ \\
\Psi_3 : \left\{ \begin{array}{l}
                        a \mapsto -1\\ b \mapsto 1
                        ,\end{array} \right.&
\Psi_4 : \left\{ \begin{array}{l}
                        a \mapsto -1\\ b \mapsto -1
                        .\end{array} \right.
                        \end{array} \]
We note that $\Psi_i(Q_8) \subseteq \F{3}$ for each $i =1,\ldots,4$. Thus $\F{3}(\chi_i)
= \F{3}$ for each character $\chi_i$ afforded by $\Psi_i$.

The irreducible $\F{9} Q_8$-representation of degree 2 is $\Psi_5 :Q_8 \mapp GL(2,9)$
where
\[ \Psi_5 : \left\{ \begin{array}{l}
                        a \mapsto \left( \begin{array}{cc}
                        1 & 1\\
                        1 & -1 \end{array} \right) \\ \\
                        b \mapsto \left( \begin{array}{cc}
                        0 & -1\\
                        1 & 0 \end{array} \right)
                        .\end{array} \right.\]
Note that $\Psi_5(Q_8) \subseteq GL(2,3)$. Thus $\F{3}(\chi_5) =\F{3}$ for the character
$\chi_5$ afforded by $\Psi_5$.

\section{$H(9,2)$ as arc-transitive factors} \label{c7:h92}
We are now ready to show that $Q_8$ does not have an irreducible $\F{3} Q_8$-representation of
degree $4$. To do that, we need the concept of an irreducible constituent of a possibly reducible
representation. Briefly (see \cite[p. 37]{Ashbacher_b86} or \cite[p. 146 - 147]{Isaacs_b76} for
more details), an \textit{irreducible constituent} of an $\FF G$-representation $\Psi$ (with
corresponding $\FF G$-module $U$) is the restriction $\Psi_i = \Psi|_{U_i/U_{i-1}}$ of $\Psi$ to
the composition factor $U_i/U_{i-1}$ of the composition series $\0 = U_0 \leq U_1 \leq \ldots \leq
U_n =U$ (and each composition factor $U_i/U_{i-1}$ is an irreducible $\FF G$-module). Clearly the
irreducible constituents of the $\FF G$-representation $\Psi$ of degree $m$ are irreducible $\FF
G$-representations of degrees $m' \leq m$.

\bl {\rm \cite[Corollary 9.23]{Isaacs_b76}} \label{vital} \ Let $\FF \subseteq \EE$ be fields of
prime characteristic. Let $\Psi$ be an irreducible $\EE$-representation of $G$ which affords the
character $\chi$. Let $\eta$ be an irreducible $\FF$-representation such that $\Psi$ is a
constituent of $\eta^\EE$. Then $deg(\eta) = | \FF(\chi):\FF| \ deg(\Psi)$. In particular, if
$\FF(\chi) = \FF$, then $\Psi$ is similar to $\eta^\EE$. \el

\bl \label{irred_rep} Each irreducible $\F{3} Q_8$-representation has degree at most
$2$.\el \bpf Let $m>2$. Suppose there exists a degree $m$ irreducible $\F{3}
Q_8$-representation: $\eta :Q_8 \mapp GL(m,3)$. Since $\F{3} \subset \F{9}$, we have
$\eta^{\F{9}} :Q_8 \mapp GL(m,9)$. Let $\Psi$ be an irreducible constituent of
$\eta^{\F{9}}$ which affords the character $\chi$. It is easy to see that the
irreducible constituent $\Psi$ of $\eta^{\F{9}}$ is one of the irreducible $\F{9}
Q_8$-representations $\Psi_1, \ldots, \Psi_5$ (as defined above), affording characters
$\chi_1, \ldots, \chi_5$ respectively. Applying Lemma~\ref{vital}, we have $m =
deg(\eta) = | \F{3}(\chi):\F{3}| deg(\Psi)$ where $\Psi = \Psi_i$ and $\chi =\chi_i$ for
$i \in \{1,2,3,4,5 \}$. Since $\F{3}(\chi) = \F{3} (\chi_i) = \F{3}$, we have
\be         m   &=& deg(\eta) \\
                &=& | \F{3}(\chi):\F{3}| deg(\Psi) \\
                &=& 1. deg(\Psi) \\
                &=& deg(\Psi) \leq 2.
\ee This is a contradiction since $deg (\eta) = m >2$. Thus the result follows. \epf

By Lemma~\ref{irred_rep}, the group $Q_8$ does not have an irreducible $\F{3}
Q_8$-representation of degree $4$. Now in $L =T \rtimes Q_8$, the group $Q_8$ acts
semiregularly (by conjugation) on $T \cong \Z_3^4$. So $Q_8$ has 10 orbits of length 8
in $V^*= V \setminus \{ \0 \}$. Let $\Or_1, \ldots, \Or_{10}$ be the orbits of $Q_8$ in
$V^*$. We will show that $\brac{\Or_i}=\F{3}^2 \cong\Z_3^2$ for some $i$ (note that by
$\brac{\Or_i}$, we mean the $\F{3}$-span of vectors in $\Or_i$).

\bl \label{orr} Suppose $\Or_1, \ldots, \Or_{10}$ are orbits of $Q_8$ in $V^*$ $($where
$V =V(4,3)=\F{3}^4$$)$. Then there exists an $\Or_i$ such that $\brac{\Or_i}=\F{3}^2
\cong \Z_3^2$, and moreover $\Or_i$ is the set of non-zero vectors in $\brac{\Or_i}$.
\el \bpf Since the characteristic of $\F{3}$ does not divide $|Q_8|$, then by Maschke's
Theorem \cite[Theorem 12.9]{Ashbacher_b86}, every $\F{3} Q_8$-module is semisimple, that
is, every $\F{3} Q_8$-module is the direct sum of irreducible submodules. So $V =
\bigoplus U_i$, where each $U_i$ is an irreducible submodule of $V$. By
Lemma~\ref{irred_rep}, we know that all irreducible $\F{3} Q_8$-representations have
degree at most $2$. Thus $dim(U_i) \leq 2$. Suppose now all $U_i$ have dimension 1. So
$V = \bigoplus_{i=1}^4 U_i$. Let $\varphi: Q_8 \mapp GL(4,3)$ be the corresponding
(completely reducible) $\F{3} Q_8$-representation with respect to $V$. Then $\varphi$
splits into a direct sum of irreducible $\F{3} Q_8$-representations of degree 1, and
with respect to an appropriate basis, $\varphi(Q_8) \subseteq \left \{ \left
[\begin{array}{cccc} a & 0 & 0 & 0
\\ 0 & b & 0 & 0 \\ 0 & 0 & c & 0 \\ 0 & 0 & 0 & d \end{array} \right]  :  a,b,c,d =
\pm 1 \right \} \cong \Z_2^4$. Furthermore, since $Q_8$ acts faithfully on $V$, we have
$Q_8 \lesssim \Z_2^4$. This is clearly not possible since $Q_8$ is non-abelian. Thus
there exists a $Q_8$-invariant irreducible $\F{3} Q_8$-submodule, say $U_k$, such that
$dim(U_k) =2$. Let $u \in U_k \setminus \{\0\}$. Then $u \in \Or_i$ for some $i$.
Observe that $\Or_i = u^{Q_8} \subset U_k$ and since $|\Or_i|=8 > |\Z_3|$, it follows
that $\brac{\Or_i} =U_k = \F{3}^2 \cong \Z_3^2$ and $U_k = \Or_i \cup \{\0\}$. \epf

\bl \label{key} Let $v \in \Or_i$ where $\Or_i$ is a $Q_8$-orbit in $V^*$. Let $S =
v^{M_0}=v^E$ $($where $v \in V^*$$)$ be an $M_0$-orbit in $V^*$. Then $S = \Or_i \cup
\Or_j$ where $\Or_j$ is an $Q_8$-orbit in $V^*$ and $\Or_i \neq \Or_j$. \el \bpf Since
$Q_8 \normlt \mathbf{E}$, the group $\mathbf{E}$ leaves the set $\{\Or_1, \ldots,
\Or_{10} \}$ of $Q_8$-orbits in $V^*$ invariant. Now $\Or_i^{\mathbf{E}} = S$ with $16=
|S| > |\Or_i| =8$. As the set of $Q_8$-orbits is invariant under $\mathbf{E}$, there
exists an $x \in \mathbf{E}$ such that $\Or_i^x = \Or_j$ where $\Or_j$ is a $Q_8$-orbit
and $\Or_j \neq \Or_i$. Since $\Or_i \cap \Or_j = \emptyset$ and $|\Or_i \cup \Or_j|
=16$, it follows that we have $S = \Or_i \cup \Or_i^x = \Or_i \cup \Or_j$. \epf

We are now ready to prove the following result.

\bt \label{hamm} Let $\Gamma_b =Cay(V,S_b)$ $($for $1 \leq b \leq 5)$ be the $M$-arc-transitive
homogeneous factors as in Section~$\ref{c4:extraspecial}$, where $V=V(4,3) \cong \Z_3^4$ and $M =T
\rtimes \mathbf{E}$. Then $\Gamma_b \cong H(9,2)$ for all $b$. \et \bpf Let $\Gamma_b =Cay(V,S_b)$
$($$1 \leq b \leq 5)$ be the $M$-arc-transitive homogeneous factors as in
Section~$\ref{c4:extraspecial}$, where $V=V(4,3) \cong \Z_3^4$ and $M = T \rtimes \mathbf{E}$
(recall that the homogeneous factorisation has index $5$; see also Proposition~\ref{espmain}). Let
$\Or_i$ be a $Q_8$-orbit in $V^*$ such that $U:= \brac{\Or_i}= \F{3}^2$ (such a $Q_8$-orbit exists
by Lemma~\ref{orr}). Since each $S_b$ is an $M_0$-orbit in $V^*$, it follows that by
Lemma~\ref{key}, we have $S_b = \Or_i \cup \Or_j$ for some $j$ and $b$, and $\Or_j = \Or_i^x$ for
some $x \in \mathbf{E}$. Furthermore, $\brac{\Or_j} =\brac{\Or_i^x} =\brac{\Or_i}^x =U^x$, and
$\Or_j$ is the set of all non-zero vectors in $U^x$. Since $\Or_j \neq \Or_i$, we have $U^x \cap U
= \{\0 \}$ and $U^x \cong \Z_3^2$. So $U \oplus U^x =V=\F{3}^4$. Without loss of generality we may
suppose $S_1 = \Or_1 \cup \Or_2$. Set $\Sigma_1 := Cay(\brac{\Or_1}, \Or_1)$ and $\Sigma_2 :=
Cay(\brac{\Or_2}, \Or_2)$ where $\brac{\Or_1} \cong \Z_3^2 \cong \brac{\Or_2}$. Since $|\Or_1| =
|\Or_2|=8$, it is easy to see that $\Sigma_1 \cong K_9 \cong \Sigma_2$ where $K_9$ is the complete
graph of order 9. From Proposition~\ref{prop:prod}, we have $\Sigma_1 \square \Sigma_2 =
Cay(\brac{\Or_1} \oplus \brac{\Or_2}, \Or_1 \cup \Or_2) = Cay(U \oplus U^x, \Or_1 \cup \Or_2) =
Cay(V, S_1)$. Thus the homogeneous factor $Cay(V,S_1)$ is the Cartesian product of 2 complete
graphs of order 9 and this is isomorphic to the Hamming graph $H(9,2)$. Finally, as all the
$M$-arc-transitive homogeneous factors $\Gamma_b =Cay(V,S_b)$ are pairwise isomorphic, it follows
that $\Gamma_b \cong H(9,2)$ for all $b$. \epf

\section{$2$-Designs} \label{des}
In this section, we will show that the homogeneous factorisation found in
Proposition~\ref{espmain} (2) (where the factors are Hamming graphs $H(9,2)$; see also
Theorem~\ref{hamm}) gives rise to an edge partition of $K_{81}$ into $90$ copies of $K_9$.
Interestingly, this partition also corresponds to a point-transitive 2-design.

\bdeff \label{design} \rm{A \textit{$2{-}(v,\kappa,1)$ design} $\M{D} =(V,\M{B})$ is a system
consisting of a finite set $V$ of $v$ points and a collection $\M{B}$ of $\kappa$-subsets of $V$
called \textit{blocks} such that any $2$-subset of $V$ is contained in exactly one block. We will
always assume that $2 < \kappa < v$. Furthermore, given a $2{-}(v,\kappa,1)$ design $\M{D}$, let
$b$ denote the number of blocks and $r$ the number of blocks containing a given point. Finally, the
\textit{order} of a $2{-}(v,\kappa,1)$ design is the integer $\eta:= r-1$.} \edeff

A permutation $g$ of the points of a $2{-}(v,\kappa,1)$ design $\M{D}=(V,\M{B})$ is called an
automorphism of $\M{D}$ if $B^g \in \M{B}$ whenever $B \in \M{B}$. The set of all automorphisms of
$\M{D}$ constitutes a group which is called the \textit{automorphism group of $\M{D}$} and is
denoted by $Aut(\M{D})$. Note that every automorphism of a design induces a permutation of the
block set. So elements of $Aut(\M{D})$ can also be seen as permutations on $\M{B}$. We say that $G
\leq Aut(\M{D})$ \textit{is point-transitive} if $G$ is transitive on the point set $V$.

\bt \label{2design} The arc-transitive homogeneous factorisation of $K_{81} =(V,E)$ into $5$
factors, each isomorphic to the Hamming graph $H(9,2)$ $($see Proposition~$\ref{espmain} (2)$$)$,
gives rise to an edge partition of $K_{81}$ into $90$ copies of $K_9$. This edge partition
corresponds to a $2{-}(81,9,1)$ design $\M{D} =(V,\M{B})$ with point set $V$ and blocks the $90$
copies of $K_9$. Furthermore, $\M{D}$ is a $G$-point, $2$-transitive affine plane isomorphic to
the exceptional nearfield plane of order $9$ where $G = T \rtimes G_0 \leq AGL(4,3)$, and $G_0$
has a normal extraspecial subgroup $\mathbf{E}$ of order $2^5$. \et

\brk \rm{For a definition of the \textit{exceptional nearfield plane}, see \cite[p. 72]{Kantor85}
and the references therein. We identify our design with this plane by applying Kantor's
classification in \cite{Kantor85} of point 2-transitive, $2{-}(v, \kappa,1)$ designs.} \erk

From Theorem~\ref{hamm}, we know that the Hamming graph $H(9,2) \cong K_9 \ \square \ K_9$ occurs
as a factor of an arc-transitive homogeneous factorisation $(M,G,V,\M{E})$ of $K_{81}$ of index
$5$. Furthermore, the group $M=T \rtimes \mathbf{E}$ (where $\mathbf{E}$ is the extraspecial
subgroup of order $2^5$) acts arc-transitively on each factor, with $G =T \rtimes G_0$ (where
$G_0$ contains $\mathbf{E}$ as a normal subgroup; case $2(h)$ of Theorem~\ref{2_trxclass})
permuting the 5 isomorphic factors $\Gamma_i$ transitively. We shall show that this factorisation
gives rise to a $2{-}(81,9,1)$ design.

Before we proceed to show that the factorisation of $K_{81}$ into 5 copies of Hamming
graphs $H(9,2)$ gives rise to a 2-design, we need the following result about maximal
cliques in $H(9,2)$. (A \textit{clique} $C$ of a graph $\Gamma$ is a complete subgraph
of $\Gamma$. We say that a clique $C$ in $\Gamma$ is \textit{maximal} if there is no
larger clique in $\Gamma$ containing $C$.)

\bl \label{maxK} Let $\Gamma_1, \Gamma_2$ be complete graphs on $9$ vertices with vertex sets
$\Omega_1$ and $\Omega_2$ respectively. Then the Hamming graph $H(9,2) \cong \Gamma_1 \ \square \
\Gamma_2$ contains $18$ maximal cliques $K_9$ such that each edge of $H(9,2)$ is in exactly one of
these maximal cliques. Moreover these maximal cliques $K_9$ are unique. \el \bpf Given any edge
$\{u,v\}$ in $H(9,2)$, there are exactly 7 vertices adjacent to both $u$ and $v$, and these 9
vertices form a clique. Therefore there is exactly one $9$-clique (or $K_9$) on each edge and no
larger cliques. Furthermore, each 9-clique has 36 edges and since $H(9,2)$ has 648 edges, it
follows that there are $\frac{648}{36} =18$ $K_9$ in $H(9,2)$. \epf

We shall now prove Theorem~\ref{2design}.

\vspace{2mm}

\noindent \textit{Proof of Theorem~\ref{2design}} \ \ By Proposition~\ref{espmain} (2), there
exists an arc-transitive homogeneous factorisation of $K_{81}=(V,E)$ of index 5 such that each
factor $\Gamma_i$ (where $1 \leq i \leq 5$) is isomorphic to the Hamming graph $H(9,2)$. By
Lemma~\ref{maxK}, each $\Gamma_i \cong H(9,2)$ contains 18 edge-disjoint copies of (unique)
maximal cliques $K_9$ and each edge of $\Gamma_i$ is contained in exactly one of them. As there
are 5 factors of $K_{81}$ isomorphic to $H(9,2)$, we have an edge partition of $K_{81}$ into 90
copies of edge-disjoint $K_9$'s which are maximal cliques of some $\Gamma_i$. Let $\M{B}$ be the
set of vertex sets of the 90 maximal cliques $K_9$. Then each $B \in \M{B}$ has size 9, and each
edge of $K_{81}$ (a 2-subset of $V$) lies in precisely one block in $\M{B}$. Thus the
factorisation given in Proposition~\ref{espmain} (2) gives rise to a $2-(81,9,1)$ design
$\M{D}=(V,\M{B})$.

We recall that $G =T \rtimes G_0$ where $G_0 \normgt \mathbf{E}$ permutes the 5 factors
$\Gamma_i \cong H(9,2)$ transitively. So elements in $G$ induce isomorphisms between
each pair of factors. Suppose $g \in G$ maps $\Gamma_i$ to $\Gamma_j$ (where $i$ and $j$
are not necessarily distinct). Then it follows that $g$ will map a maximal clique $K_9$
(belonging to $\M{B}$) in $\Gamma_i$ to a (not necessarily distinct) maximal clique
$K_9$ (also belonging to $\M{B}$) in $\Gamma_j$. Thus $G$ leaves the block set $\M{B}$
invariant and so $G \leq Aut(\M{D})$. Moreover, $G$ is 2-transitive on the point set $V$
of $\M{D}$.

Finally, $|\M{B}| = b =90$ and $\kappa = |V(K_9)|= 9$ (where $V(K_9)$ is the vertex set of $K_9$).
Then by \cite[p. 57]{Dem}, which says that $rv =b \kappa$ (where $r$ is as in
Definition~\ref{design}), we have $r =(b \kappa) /v = 10$. Note that for the above $2-(81,9,1)$
design $\M{D}$, the order of $\M{D}$ is $\eta = \kappa = r-1=9$, and by \cite[2.2.6]{Dem}, $\M{D}$
is an affine plane. (Simply, an \textit{affine plane} is any $2{-}(\kappa^2, \kappa,1)$ design.) By
Kantor's classification of point 2-transitive $2{-}(v, \kappa,1)$ designs in \cite{Kantor85}, it
follows that $\M{D}$ is the ``exceptional near field plane". \qedd


\section{Final remarks} \label{c7:finale}
In this section, we shall make a few comments and observations on the results obtained in
this chapter.

\vspace{0.7cm}

\noindent (1) \ Note that by Proposition~\ref{prop:gpaley}, the generalised Paley graphs $\Gamma_i
= GPaley(81,16)$ ($i = 1, \ldots, 5$) are factors of an arc-transitive homogeneous factorisation
$(M,G,V,\M{E}_{GP}(81,5))$ of $K_{81}$ of index $5$, where $M = T \rtimes \brac{\omb^5} \normlt
G=T \rtimes \brac{\omb} < A\Gamma L(1,81)$. By Theorem~\ref{thm:A}, it follows that $GPaley(81,16)
\cong H(9,2)$. Thus the factors arising from the arc-transitive homogeneous factorisation in
Section~\ref{c4:extraspecial}, which are isomorphic to the Hamming graphs $H(9,2)$ by
Theorem~\ref{hamm}, are in fact isomorphic to the generalised Paley graphs $GPaley(81,16)$.
However, we have not decided if the partition $\M{E}$ in Section~\ref{c4:extraspecial} is
equivalent to the generalised Paley partition $\M{E}_{GP}(81,5)$. (See Remark~\ref{rk:partition}
for definition of equivalent partition.)

\vspace{0.7cm}

\noindent (2) \ The homogeneous factorisation of $K_{25}$ into three Hamming graphs $H(5,2)$ in
Section~\ref{c4:Gqk} also gives a similar edge partition to those described in the previous
section. In this case, the edge set $E(K_{25})$ (of $K_{25}$) is partitioned into $30=3 \times 10$
copies of $K_5$ and we have a $2{-}(25,5,1)$ design $\M{D}$. Furthermore, $\M{D}$ is a $G$-point,
2-transitive affine plane of order $5$ where $G = T \rtimes G_0 \leq AGL(2,5)$ and $G_0 \normgeq
SL(2,3)$. (The proof is analogous to that of Theorem~\ref{2design}.) Also, by looking up Kantor's
classification of point 2-transitive $2{-}(v, \kappa,1)$ designs in \cite{Kantor85}, we know that
$\M{D}$ is a Desarguesian affine plane.

\vspace{0.7cm}

\noindent (3) \ Similar observations (using an analogous approach to the proof of
Theorem~\ref{2design}) can also be made for homogeneous factorisations of $K_{p^R}$ into
generalised Paley graphs. By Theorem~\ref{thm:A}, the generalised Paley graphs $\Gamma=\GP$ such
that $k =\frac{a(p^R-1)}{R(p^a-1)}$ where $a \mid R$ and $a \neq R$, are isomorphic to the Hamming
graphs $H(p^a, \frac{R}{a})$. Now if $a = \frac{R}{2}$, then $k = \frac{1}{2}(p^{a}+1)$ and
$\Gamma \cong H(p^{a}, 2)$. Now the generalised Paley graphs are factors of arc-transitive
homogeneous factorisation $(M,G,V,\M{E}_{GP}(p^R,k))$ of $K_{p^R}$. Thus it follows that whenever
a complete graph $K_{p^R}$ is factorised into $\GP \cong H(p^{a}, 2)$ ($a=\frac{R}{2}$), then its
edge set is partitioned into $k \times 2p^{a} = \frac{1}{2}(p^{a}+1) \times 2p^{a} = p^{2a} +p^a$
copies of $K_{p^{a}}$ and we have a $2{-}(p^{2a},p^a,1)$ design $\M{D}$. Similar to the previous
observation, $\M{D}$ is also a $G$-point, 2-transitive affine plane of order $5$ where $G = T
\rtimes \brac{\omb} \leq A\Gamma L(1,p^{R})$. Finally, by Kantor's result in \cite{Kantor85},
$\M{D}$ is a Desarguesian affine plane.


\chapter{The one-dimensional affine case} \label{c5}

The only case remaining from the list in Theorem~\ref{2_trxclass} is case $2(a)$, which
we call the one-dimensional affine case. Here, $G = T \rtimes G_0$ is a $2$-transitive
subgroup of $A\Gamma L(1,p^R)$ acting on $V$ where $V$ is a finite field $\F{p^R}$ of
order $p^R$ for some prime $p$ and $R \geq 1$. As in the previous chapter, we need to
determine all possible transitive normal subgroups $M =T \rtimes M_0$ of $G$ such that
$M_0$ (which is also normal in $G_0$) is not transitive on $V^*$ and $|M_0|$ has even
order if $p$ is odd. Also as in the previous chapter, the problem of determining $G$ and
$M$ in $A\Gamma L(1,p^R)$ is equivalent to the problem of finding $G_0$ and $M_0$ in
$\Gamma L(1,p^R)$ such that $G_0$ is transitive on $V^*$, $M_0 \normlt G_0$ and $M_0$ is
not transitive on $V^* =V \setminus \{\0\}$.

In this chapter, we give generic constructions for both $G_0$ and $M_0$ in terms of a
fixed set of parameters and use them to construct arc-transitive homogeneous
factorisations $(M,G,V,\M{E})$ of $K_{p^R}=(V,E)$ where $V=\F{p^R}$ and $G \leq A\Gamma
L(1,p^R)$.

\section{Generic constructions for $M_0$ and $G_0$} \label{c5:std}
For the sake of completeness, we shall start off by repeating the explanations of the
notations used in Section~\ref{c6:main} (see also Section~\ref{c3:examples}).

Let $V=\F{p^R}$ and $V^* = V \setminus \{\0\}$ be as above. For a fixed primitive
element $\om$ in $V$, we let $\omb$ be the corresponding scalar multiplication $\omb: x
\mapp x\om$ for all $x \in V$. Then let $\brac{\omb^i}$ be the multiplicative subgroup
of $V^*$ generated by the element $\omb^i$ where $i \geq 1$ and $i \mid (q-1)$ (note
that $\omb^i: x \mapp x\om^i$, $x \in V$). Since for each $\omb^i$, there is a
corresponding $\om^i \in V$, we shall use $\brac{\om^i}$ to denote the set
$\{1,\om^i,\om^{2i}, \ldots, \om^{((q-1)/i )-i} \} \subseteq V^*$. Note that for $1, -1
\in V^*$, we will use $\pone$ and $\mone$ to denote the corresponding scalar
multiplications. Also, we let $\al$ denote the Frobenius automorphism of $\F{p^R}$, that
is, $\al: x \mapp x^p$. Finally, we note that the group $\Gamma L(1,p^R)$ is generated
by the scalar multiplication map $\omb$ and the Frobenius automorphism $\al$. Thus
$\Gamma L(1,p^R) =\brac{\omb, \al}$.

The following result from \cite{Foulser64} tells us that each subgroup $G_0$ of $\Gamma
L(1,p^R)$ can be expressed in the form $G_0 =\brac{\omb^d, \omb^e \al^s}$ where the
integers $d,e$ and $s$ can be chosen in a unique manner.

\bl \label{stdlem} $\cite[\mbox{Lemma $4.1$}]{Foulser64}$ Let $G_0$ of $\Gamma L(1,p^R)
=\brac{\omb,\al}$. Then there exist unique integers $d$, $e$ and $s$ such that $G_0
=\brac{\omb^d, \omb^e \al^s}$, and the following hold: \bnum
\item $d>0$ and $d \mid (p^R-1)$;
\item $s>0$ and $s \mid R$;
\item $0 \leq e < d$ and  $e(p^R-1)/(p^s-1) \equiv 0$ $(mod$ $d)$.
\enum \el

\bdeff \label{stddef} \rm{\textbf{(Standard Form)} If $G_0 =\brac{\omb^d,\omb^e \al^s}
\leq \Gamma L(1,p^R)$ and the integers $d,e$ and $s$ satisfy conditions (1)-(3) of
Lemma~\ref{stdlem}, then the representation $G_0 =\brac{\omb^d,\omb^e \al^s}$ is said to
be in \textit{standard form}.} \edeff

\brk \label{stdrk} \rm{In subsequent results (for example, see
Lemmas~\ref{gtran}-\ref{mintran}), we will always require a subgroup $G_0
=\brac{\omb^d,\omb^e \al^s}$ of $\Gamma L(1,p^R)$ to be in standard form. Otherwise, it
will cause ambiguities especially when certain restrictions or conditions are placed on
the 3-tuple $(d,e,s)$. For example, if $p=3$ then $G_0 =\brac{\omb^3,\al}$ is not
expressed in standard form since condition (1) of Lemma~\ref{stdlem} fails. The standard
form for this subgroup $G_0$ is $G_0 = \brac{\omb,\al}$ and by Lemma~\ref{stdlem}, we
know that this expression in standard form is unique.} \erk

To determine the 2-transitive subgroups $G=T \rtimes G_0$ of $A\Gamma L(1,p^R)$, we need
to know the conditions under which $G_0$ is transitive on $V^*$, that is, we need the
possible integers $d,e$ and $s$ such that $G_0 = \brac{\omb^d, \omb^e \al^s}$ is in
standard form and is transitive on the set of non-zero elements of $\F{p^R}$. The
transitivity criteria in Lemma~\ref{gtran} may be found in \cite[Section 3]{FK78}. We
provide a proof here as we need the details for determining the possibilities for $M_0$.

Let $G_0 =\brac{\omb^d, \omb^e \al^s}$ be in standard form and consider the orbits of
$H:= \brac{\omb^d}$ in $V^*$. Since $H$ acts semiregularly on $V^*$, it follows that $H$
has $d$ orbits of length $(p^R-1)/d$ in $V^*$. Let $\Omega:=\{\brac{\om^d},\om
\brac{\omega^d},\ldots, \omega^{d-1}\brac{\omega^d}\}$ be the set of $H$-orbits in
$V^*$. Now $\tau := \omb^e \alpha^s$ induces a permutation on the set $\Omega$ and $G_0$
is transitive on $V^*$ if and only if $\tau$ has one orbit in $\Omega$. Consider the
action of $\tau$ on $\Omega$. To determine the image of $\om^i \brac{\om^d}$ under
$\tau$, we simply need to find the ``coset" of $\brac{\om^d}$ containing
$(\om^i)^{\tau}$. Now $(\omega^i)^\tau = \om^{ip^s+ep^s}$ lies in $\om^j \brac{\om^d}$
where $ip^s+ep^s \equiv j$ (mod $d$) and $0 \leq j \leq d-1$. Further computations yield
$(\om^i)^{\tau^b} = \om^{ip^{bs}}. \om^{ep^s((p^{bs}-1)/(p^s-1))}$. In the proof below,
we shall use this observation for the case where $i=0$.

\bl \label{gtran} {\rm \textbf{(Transitivity)}} Suppose $G_0 = \brac{\omb^d, \omb^e
\al^s} \leq \Gamma L(1,p^R)$ is in standard form. Then $G_0$ is transitive on $V^*$ if
and only if either $d=1$ $($so $e=0$$)$, or both of the following hold$:$ \bnum
\item $e>0$, $d$ divides $e((p^{ds}-1)/(p^s-1))$, and
\item if $1< d' <d$, then $d$ does not divide $e((p^{d's}-1)/(p^s-1))$. \enum \el
\bpf Suppose first that $e=0$. Then $G_0 =\brac{ \omb^d, \al^s} = \brac{\omb^d} \rtimes
\brac{\al^s}$. Since $(\omb^d)^{\al^s} = \omb^{dp^s} \in \brac{\omb^d}$, it follows that
$\brac{\al^s}$ fixes each orbit of $\brac{\omb^d}$ in $V^*$ setwise. Thus, $G_0$ is
transitive on $V^*$ if and only if $d=1$, that is, $G_0 =\brac{\omb, \al^s}$.

Suppose now that $e \neq 0$. Then $G_0$ is transitive on $V^*$ if and only if $\brac{\tau}$ is
transitive on $\Omega$ where $\tau =\omb^e \al^s$ and $\Omega$ is as defined above. Using the
information above about the action of $\tau$ on $\Omega$, we have: \be \textrm{$\brac{\tau}$ is
transitive on $\Omega$} & \ifff & \textrm{(1) $\tau^d$ fixes $\brac{\om^d}$, and}\\ & &
\textrm{(2) if $1 \leq d' <d$, then $\tau^{d'}$ does not fix $\brac{\om^d}$}\\ & \ifff &
\textrm{(1) $\om^{ep^s((p^{ds}-1)/(p^s-1))} \in \brac{\om^d}$, and} \\ & & \textrm{(2) if $1 \leq
d'< d$, then $\om^{ep^s((p^{d's}-1)/(p^s-1))} \notin \brac{\om^d}$}\\ & \ifff & \textrm{(1) $d$
divides $ep^s((p^{ds}-1)/(p^s-1))$, and}\\ & & \textrm{(2) if $1 \leq d' <d$, then $d$ does not
divide $ep^s \left(\frac{p^{d's}-1}{p^s-1} \right)$.} \ee

Since $d \mid (p^R-1)$, it follows that $gcd(p,d)=1$. So \be \textrm{$\brac{\tau}$ is
transitive on $\Omega$} & \ifff & \textrm{(1) $d$ divides $e((p^{ds}-1)/(p^s-1))$,
and}\\ & & \textrm{(2) if $1 \leq d' <d$, then $d$ does not divide $e \left(
\frac{p^{d's}-1}{p^s-1} \right)$.} \ee \epf

Recall that, given a finite 2-transitive permutation group $G$ which permutes the
factors of some homogeneous factorisation of $K_{p^R}$, we need to determine the
possible normal subgroups $M$ of $G$ that are vertex transitive and arc transitive on
each factor. In particular, given a subgroup $G_0$ of $\Gamma L(1,p^R)$ which is
transitive on $V^*$, we want to know the possible normal subgroups $M_0$ of $G_0$ that
are not transitive on $V^*$. Let $M_0:= \brac{\omb^{d_1},\omb^{e_1} \al^{s_1}}$ be a
subgroup of $\Gamma L(1,p^R)$ expressed in standard form. Then the following results
tell us the conditions on the integers $d_1, e_1$ and $s_1$ under which the subgroup
$M_0$ is a normal subgroup of $G_0$. We also calculate the length of the $M_0$-orbits in
$V^*$ under the assumption that $M_0$ is normal in $G_0$.

\bl \label{contain} {\rm \textbf{(Containment)}} Suppose $M_0 =\brac{\omb^{d_1},
\omb^{e_1} \al^{s_1}}$ and $G_0 =\brac{\omb^d, \omb^e,\al}$ are subrgroups of $\Gamma
L(1,p^R)$ expressed in standard form. Then $M_0$ is a subgroup of $G_0$ if and only if
\bnum
\item $d \mid d_1$,
\item $s \mid s_1$,
\item and $d \mid (e \frac{(p^{s_1}-1)}{(p^s-1)} -e_1)$.
\enum \el \bpf Suppose $M_0 \leq G_0$. Then $\omb^{d_1}$ and $\omb^{e_1} \al^{s_1}$ are
elements of $G_0$. Let $B = \brac{\omb}$. Then $G_0 \cap B = \brac{\omb^d}$ contains $M_0
\cap B = \brac{\omb^{d_1}}$, so $|\brac{\omb^{d_1}}|$ divides $|\brac{\omb^{d}}|$ and
hence $d \mid d_1$. Also, we have $M_0/(M_0 \cap B) \cong M_0B/B \lesssim G_0B/B \cong
G_0/(G_0 \cap B)$. Since $M_0/(M_0 \cap B) = \brac{\brac{\omb^{d_1}} \omb^{e_1}
\al^{s_1}} \cong \Z_{R/s_1}$ and $G_0/(G_0 \cap B) = \brac{\brac{\omb^d} \omb^e \al^s}
\cong \Z_{R/s}$, it follows that $s \mid s_1$.

Given that $s$ divides $s_1$ (as shown above), we have $(\omb^e \al^s)^{s_1/s} \in G_0$.
Now by \cite[Lemma 2.1]{FK78}, for each $i \geq 1$, $(\omb^e \al^s)^i = \omb^J \al^{si}$
where $J \equiv e(p^{si}-1)/(p^s-1)$ (mod $p^R-1$), and hence
\[
(\omb^e \al^s)^{s_1/s} = \omb^{J'} \al^{s_1},
\]
where $J' = e(p^{s_1}-1)/(p^s-1)$. Writing
\[
(\omb^e \al^s)^{s_1/s} = \omb^{J'} \al^{s_1} = \omb^{J'} \omb^{-e_1} \omb^{e_1}
\al^{s_1},
\]
and since $(\omb^e \al^s)^{s_1/s}$, $\omb^{e_1} \al^{s_1} \in G_0$, we see that $\omb^{J' -e_1} \in
G_0$ and this is true if and only if $d \mid (J' -e_1)$.

Conversely, suppose the three conditions of the lemma are satisfied. Then since $d \mid
d_1$, there exists an integer $j$ such that $(\omb^d)^j = \omb^{d_1}$. Thus $\omb^{d_1}
\in G_0$. Now $s \mid s_1$ and $d \mid (e \frac{(p^{s_1}-1)}{(p^s-1)} -e_1)$. As above,
by \cite[Lemma 2.1]{FK78}, we get $(\omb^e \al^s)^{s_1/s} \in G_0$. However, we know that
\[
(\omb^e \al^s)^{s_1/s}= \omb^{J'} \al^{s_1} = \omb^{J' -e_1} \omb^{e_1} \al^{s_1} \in
G_0,
\]
where $J' - e_1 = (e \frac{(p^{s_1}-1)}{(p^s-1)} -e_1)$. Since $ d \mid (J' -e_1)$, we
have $\omb^{J' -e_1} \in G_0$, forcing $\omb^{e_1} \al^{s_1}$ to be in $G_0$. \epf

\bl \label{normal} {\rm \textbf{(Normality)}} Suppose $M_0 = \brac{\omb^{d_1},\omb^{e_1}
\al^{s_1}}$ is in standard form and is a subgroup of $G_0 = \brac{\omb^d, \omb^e \al^s}$ $($so
Lemma~$\ref{contain}$ holds$)$. Then $M_0$ is normal in $G_0$ if and only if \bnum

\item $d_1 \mid d(p^{s_1} -1)$ and

\item $d_1 \mid (e_1(p^s-1) + ep^s(p^{R-s_1} -1))$.

\enum \el \bpf Now $M_0$ is normal in $G_0$ if and only if $(\omb^{d_1})^g \in M_0$ and
$(\omb^{e_1} \al^{s_1})^g \in M_0$ for all $g \in G_0$. Since $\brac{\omb^{d_1}} \normlt
G_0$ whenever $M_0 = \brac{\omb^{d_1},\omb^{e_1} \al^{s_1}}$ is a subgroup of $G_0$, it
follows that $M_0$ is normal in $G_0$ if and only if $(\omb^{e_1} \al^{s_1})^g \in M_0$
for all $g \in G_0$. Furthermore, since $G_0 = \brac{\omb^d ,\omb^e \al^s}$, we have that
$M_0$ is normal in $G_0$ if and only if $(\omb^{e_1} \al^{s_1})^{\omb^d} \in M_0$ and
$(\omb^{e_1} \al^{s_1})^{\omb^e \al^s} \in M_0$. Now \be
(\omb^{e_1} \al^{s_1})^{\omb^d} &=& \omb^{-d} \omb^{e_1} \al^{s_1} \omb^d \\
                              &=& \omb^{e_1-d} \al^{s_1} \omb^d \\
                              &=& \omb^{e_1}\al^{s_1} \omb^{d-dp^{s_1}}.
\ee Thus \be
(\omb^{e_1} \al^{s_1})^{\omb^d} \in M_0 & \ifff & \omb^{d-dp^{s_1}} \in M_0 \\
                                      & \ifff & \omb^{d-dp^{s_1}} \in M_0 \cap \brac{\omb} = \brac{\omb^{d_1}} \\
                                      & \ifff & d_1 \mid d(p^{s_1} -1).
\ee Next consider $(\omb^{e_1} \al^{s_1})^{\omb^e \al^s}$. Then
\be (\omb^{e_1} \al^{s_1})^{\omb^e \al^s} &=& (\omb^{e_1 -e} \al^{s_1} \omb^e)^{\al^s} \\
                                        &=& \omb^{(e_1 -e)p^s} (\al^{s_1} \omb^e)^{\al^s} \\
                                        &=& \omb^{(e_1 -e)p^s} \al^{s_1} \omb^{e p^s} \\
                                        &=& \omb^{(e_1 -e)p^s} \al^{s_1} \omb^{e p^s}
                                        (\al^{s_1})^{-1} \al^{s_1} \\
                                        &=& \omb^{(e_1 -e)p^s} (\omb^{e
                                        p^s})^{(\al^{s_1})^{-1}} \al^{s_1} \\
                                        &=& \omb^{(e_1 -e)p^s} \omb^{e p^{R-s_1+s}}
                                        \al^{s_1} \\
                                        &=& \omb^{(e_1 -e)p^s + ep^{R-s_1+s} -e_1}
                                        (\omb^{e_1} \al^{s_1}).
\ee Hence $(\omb^{e_1} \al^{s_1})^{\omb^e \al^s} \in M_0$ if and only if $\omb^{(e_1
-e)p^s + ep^{R-s_1+s} -e_1} \in M_0$, and \be
\omb^{(e_1 -e)p^s + ep^{R-s_1+s} -e_1} \in M_0 & \ifff & \omb^{(e_1 -e)p^s + ep^{R-s_1+s} -e_1} \in M_0 \cap \brac{\omb} \\
                                            & \ifff & d_1 \mid ((e_1 -e)p^s + ep^{R-s_1+s}
                                            -e_1) \\
                                            & \ifff & d_1 \mid (e_1(p^s-1) +
                                            ep^s(p^{R-s_1} -1)).
\ee Thus $M_0$ is normal in $G_0$ if and only if (1) $d_1 \mid d(p^{s_1} -1)$ and (2)
$d_1 \mid (e_1(p^s-1) + ep^s(p^{R-s_1} -1))$. \epf

\bl \label{mintran} {\rm \textbf{(Orbit Length)}} Let $M_0 = \brac{\omb^{d_1},\omb^{e_1}
\al^{s_1}}$ and $G_0 = \brac{\omb^{d},\omb^{e} \al^{s}}$ be subgroups of $\Gamma
L(1,p^R)$ expressed in standard form. Suppose also that $M_0$ is a normal subgroup of
$G_0$ and $G_0$ is transitive on $V^*$. Then $M_0$ has $t_0 =d_1/c$ orbits of equal
length $(p^R-1)/t_0$ in $V^*$, where if $e_1=0$ then $c =1$; and if $e_1 \neq 0$, then
$c$ is defined by$:$ \bnum
\item $d_1 \mid e_1(p^{cs_1} -1)/(p^{s_1}-1)$ and
\item $d_1 \nmid e_1(p^{c's_1} -1)/(p^{s_1}-1)$ for $c' < c $.
\enum \el

Note that conditions (1) and (2) above certaintly define an integer $c \leq
\frac{R}{s_1}$, since by Lemma~\ref{stdlem} (3) applied to $M_0$, $d_1$ divides $e_1
\frac{p^R-1}{p^{s_1}-1}$.

\bpf (\textit{of Lemma~$\ref{mintran}$}) \ Suppose $e_1 \neq 0$. We first determine the
orbits of $H':=\brac{\omb^{d_1}}$ in $V^*$. Now $H'$ acts semiregularly on $V^*$, so $H'$
has $d_1$ orbits of length $(p^R-1)/d_1$. Let $\Omega_1:=\{\brac{\omega^{d_1}},\omega
\brac{\omega^{d_1}},\ldots, \omega^{d_1-1}\brac{\omega^{d_1}}\}$ be the set of
$H'$-orbits in $V^*$. Let $\tau_1 := \omb^{e_1} \al^{s_1}$. Then $\tau_1$ induces a
permutation on the set $\Omega_1$. Let $c$ be the length of the $\brac{\tau_1}$-orbit in
$\Omega_1$ containing $\brac{\om^{d_1}}$. We claim that all the $\brac{\tau_1}$-orbits in
$\Omega_1$ have equal length $c$ and hence $c \mid d_1$.

Now $\brac{\omb^{d_1}}$ is characteristic in $\brac{\omb^d}$ (any subgroup of a cyclic
group is characteristic) and since $\brac{\omb^d} \normlt G_0$, we have
$\brac{\omb^{d_1}} \normlt G_0$. Thus $\Omega_1$ is invariant under $G_0$ and $G_0$ is
transitive on $\Omega_1$. From \cite[Theorem 10.3]{wie} we know that every nontrivial
normal subgroup of a transitive group is half-transitive (that is, orbits have equal
length). Since $M_0 \normlt G_0$ and $H'$ fixes $\Omega_1$ `pointwise', the group
induced by $M_0$ on $\Omega_1$ is equal to the group induced by $\brac{\tau_1}$. Thus the
$\brac{\tau_1}$-orbits in $\Omega_1$ have equal length $c$ and so $c \mid d_1$.

Furthermore, the group $M_0$ has orbits in $V^*$ of equal length $(p^R-1)c/d_1$, and $c$
is such that all of the following conditions hold: \bnum
\item $\tau_1^{c}$ fixes $\brac{\om^{d_1}}$ and
\item if $c' < c$, then $\tau_1^{c'}$ does not fix $\brac{\om^{d_1}}$.
\enum Using similar arguments to those in the proof of Lemma~\ref{gtran}, we have \be
\textrm{$M_0$-orbits in $V^*$ are of equal length $(p^R-1)c/d_1$} & \ifff & \textrm{(1)
$d_1 \mid e_1 \frac{p^{cs_1} -1}{p^{s_1}-1}$, and}\\ & & \textrm{(2) $d_1
 \nmid e_1 \frac{p^{c's_1} -1}{p^{s_1}-1}$ for $c' <c $.} \ee

Now suppose $e_1= 0$. Then since $\al$ fixes each of the orbits of $\brac{\omb^{d_1}}$ in
$V^*$ setwise, $M_0$ has $t_0=d_1$ orbits of equal length $(p^R-1)/d_1$ in $V^*$. \epf

We want to make some observations pertaining to Lemma~\ref{mintran}.

\brk \label{orlength} \rm{Let $M_0$, $G_0$, $t_0$ and $c$ be as in Lemma~\ref{mintran}.
\bnum \item Note that $M_0$ is transitive on $V^*$ when $t_0=d_1=1$, $e_1=0$, and also
when $t_0 = d_1/c =1$, $e_1 \neq 0$.
\item Suppose that $(M,G,V,\mathcal{E})$ (where $M, G \leq A\Gamma L(1,p^R)$ and $V=\F{p^R}$)
is an arc-transitive homogeneous factorisation of $K_{p^R}=(V,E)$ of index $k$ with
factors $\Gamma_i$. Suppose also that $M_0, G_0$ of $M, G \leq A\Gamma L(1,p^R)$ are
expressed in standard form. Then the $\Gamma_i$ are Cayley graphs $Cay(V,S_i)$
(Lemma~\ref{affcayfac}) such that $M_0$ is transitive on each of the $S_i$ (which are
the $M_0$-orbits in $V^*$), with $S_i =-S_i$ and $|S_i| =\frac{p^R-1}{k}$ even when $p$
is odd. Now from Lemma~\ref{mintran}, the normal subgroup $M_0$ of $G_0$ has $t_0$
orbits of equal length $c(p^R-1)/d_1$ in $V^*$ where $t_0 = d_1/c$ and $c=1$ if $e_1=0$.
Thus it follows (by considering the equality $\frac{c(p^R-1)}{d_1} = \frac{p^R-1}{k}$)
that $k =\frac{d_1}{c} =t_0$. \enum } \erk

We shall now use the above results (Lemmas~\ref{stdlem} - \ref{mintran}) on $G_0$ and
$M_0$ to prove the following theorem.

\bt \label{thm1} \ Let $M_0 = \brac{\omb^{d_1},\omb^{e_1} \al^{s_1}}$ and $G_0 =
\brac{\omb^d, \omb^e \al^s}$ be subgroups of $\Gamma L(1,p^R)$ $=$ $\brac{\omb,\al}$
acting on $V^*$ where $ 0 < s,s_1$, $0 \leq e_1 <d_1$ and $0 \leq e < d$. Then $M_0$ and
$G_0$ are in standard form and $M_0$ is a normal subgroup of $G_0$ if and only if the
following conditions are satisfied$:$  \newcounter{kk}
\begin{list} {{\rm \textbf{(\arabic{kk})}}}{\usecounter{kk}}
\item $d_1 \mid (p^R-1)$,
\item $s_1 \mid R$,
\item $e_1(p^R-1)/(p^{s_1}-1) \equiv 0$ $(mod$ $d_1)$.
\item $d \mid d_1$,
\item $s \mid s_1$,
\item $d \mid (e \frac{(p^{s_1}-1)}{(p^s-1)} -e_1)$,
\item $d_1 \mid d(p^{s_1} -1)$ and
\item $d_1 \mid (e_1(p^s-1) + ep^s(p^{R-s_1} -1))$. \end{list}

\vspace{2mm}

\noindent Furthermore, given that $G_0$ is in standard form, $G_0$ is transitive on $V^*$
if and only if $d=1$ and $e=0$ or both of the following conditions hold$:$
\addtocounter{kk}{1}
\begin{list} {{\rm \textbf{(\arabic{kk})}}}
\item $e>0$, $d$ divides $e((p^{ds}-1)/(p^s-1))$, and \end{list}
\addtocounter{kk}{1} \begin{list} {{\rm \textbf{(\arabic{kk})}}}
\item if $1 < d' <d$, then $d$ does not divide $e((p^{d's}-1)/(p^s-1))$. \end{list}

\vspace{3mm}

\noindent In this case, $M_0$ has $d_1/c$ orbits of equal length $c(p^R-1)/d_1$ in $V^*$,
where if $e_1=0$ then $c =1$; and if $e_1 \neq 0$, then $c$ is defined by$:$
\addtocounter{kk}{1}
\begin{list} {{\rm \textbf{(\arabic{kk})}}}
\item $d_1 \mid e_1(p^{cs_1} -1)/(p^{s_1}-1)$ and \end{list}
\addtocounter{kk}{1} \begin{list} {{\rm \textbf{(\arabic{kk})}}}
\item $d_1 \nmid e_1(p^{c's_1} -1)/(p^{s_1}-1)$ for $c' < c $.
\end{list} \et \bpf \ If $M_0$, $G_0$ are in standard form and $M_0$ is a normal subgroup
of $G_0$, then conditions (1) to (8) follow from Lemmas~\ref{stdlem}, \ref{contain} and
\ref{normal}. Further by Lemma~\ref{gtran}, $G_0$ is transitive on $V^*$ if and only if either
$d=1$, $e=0$ or conditions (9) and (10) hold.

Conversely suppose that conditions (1) to (8) hold. Then by Lemma~\ref{stdlem}, conditions (1) to
(3) above imply that $M_0$ is in standard form. We claim that conditions (1) - (6) imply that
$G_0$ is also in standard form.

It is easy to see that conditions (1) and (4) imply $d \mid (p^R-1)$; while conditions
(2) and (5) imply $s \mid R$. Also, from (6), we have  $d \mid (e
\frac{(p^{s_1}-1)}{(p^s-1)} -e_1)$, and hence (multiplying by $(p^R-1)/(p^{s_1}-1)$ (an
integer by (2)))
\[
d \mbox{ divides } e (\frac{p^R-1}{p^s-1}) - e_1 (\frac{p^R-1}{p^{s_1}-1})  \hspace{2cm}
(*).
\]
Now by conditions (3) and (4), we have that $d \mid e_1 (\frac{p^R-1}{p^{s_1}-1})$. Thus
it follows (using $(*)$), that
\[ d \mid e \frac{(p^R-1)}{(p^s-1)}, \]
and so the three conditions of Lemma~\ref{stdlem} are satisfied and $G_0$ is in standard
form.

Next, conditions (4) to (6) and Lemma~\ref{contain} imply that $M_0 \leq G_0$, and conditions (7)
-- (8) and Lemma~\ref{normal} imply that $M_0$ is a normal subgroup of $G_0$.

Finally, the results on the number and length of $M_0$-orbits in $V^*$ follows
immediately from Lemma~\ref{mintran}. \epf

\bct \label{generic} \rm{\textbf{(Generic Construction)} \ Let $K_{p^R}=(V,E)$ be a complete graph
with vertex set $V=\F{p^R}$. Suppose $G=T \rtimes G_0$ and $M =T \rtimes M_0$ are subgroups of
$A\Gamma L(1,p^R)$ acting on $V$ such that $G_0$ and $M_0$ are as in Theorem~$\ref{thm1}$ (that
is, $G_0$ and $M_0$ satisfy all the conditions listed in Theorem~$\ref{thm1}$). Suppose also that
$M_0$ has even order if $p$ is odd. Then let $\M{S} =\{S_1, \ldots, S_k \}$ be the set of all
$M_0$-orbits in $V^*$, and for each $S_i$, we define $\Gamma_i =Cay(V,S_i)$. Finally, let $\M{E}
=\{E\Gamma_i \mid 1 \leq i \leq k\}$ where $E\Gamma_i =\{\{u,v\} \mid v-u \in S_i\}$. } \ect

We shall prove that the generic construction given above for the case where $M$ and $G$
are contained in $A\Gamma L(1,p^R)$ gives rise to arc-transitive homogeneous
factorisations, and that all such factorisations in the one-dimensional case arise in
this way.

\bt \label{thm2} \ The $4$-tuple $(M,G,V,\mathcal{E})$ from Construction~$\ref{generic}$
is a homogeneous factorisation of $K_{p^R}=(V,E)$ of index $k=d_1/c$ $($where $d_1$ and
$c$ are as in Theorem~$\ref{thm1}$$)$ such that each factor $\Gamma_i$ $($$1 \leq i \leq
k$$)$ is undirected and $M$-arc-transitive. Conversely, each arc-transitive homogeneous
factorisation $(M,G,V,\M{E})$ of $K_{p^R}=(V,E)$ of index $k$, with factors $\{\Gamma_i
\mid 1 \leq i \leq k \}$, $G$ a $2$-transitive subgroup of $A\Gamma L(1,p^R)$ and $|M|$
even, can be constructed using Construction~$\ref{generic}$. \et \bpf Let
$(M,G,V,\mathcal{E})$ be the 4-tuple given by Construction~\ref{generic}. By
Theorem~\ref{thm1}, $M_0$ is a normal subgroup of $G_0$ such that $M_0$ is intransitive
on $V^*$ while $G_0$ is transitive on $V^*$. Thus it follows that $G=T \rtimes G_0$ is
2-transitive on $V$ while $M =T \rtimes M_0$ is not (but $M$ is transitive on $V$).
Furthermore, since $|M_0|$ is even if $p$ is odd, it follows that $M$ has even order.

Now observe that the translation group $T$ fixes each $E\Gamma_i$ setwise. Furthermore, as each
$S_i$ is an $M_0$-orbit in $V^*$, it follows that each $\Gamma_i$ admits $M=T \rtimes M_0$ as a
subgroup of automorphisms and is $M$-arc-transitive (see Lemma~\ref{fact1}). Thus each arc set
$A\Gamma_i$ is a nontrivial $M$-orbital. Since $G$ is 2-transitive with normal subgroup $M$ such
that $|M|$ is even, it follows by Proposition~\ref{selfpaired} that all nontrivial $M$-orbitals in
$V$ are self-paired. So elements of each $A\Gamma_i$ come in the form $\{(u,v), (v,u) \}$ and by
definition, $\Gamma_i$ is undirected. Also by Theorem~\ref{thm1}, the number of $M_0$-orbits in
$V^*$ is $d_1/c$ (where $d_1$ and $c$ are as in Theorem~$\ref{thm1}$), and so $k =d_1/c$.

Next we show that $\M{E}$ form a partition of the edge set $E$ of $K_{p^R}$. Note that
each $E\Gamma_i =\{ \{u,v\} \mid v-u \in S_i \}$ corresponds to an $M_0$-orbit $S_i$ in
$V^*$, where $S_i \neq S_j$ for $i \neq j$. It then follows that $E\Gamma_i \neq
E\Gamma_j$ for $i \neq j$. Now take an edge $\{u,v\} \in E$. Then we must have $v-u \in
S_i$ for some $i$, and so $\{u,v\} \in E\Gamma_i$. Thus $\M{E}$ forms a partition of the
edge set $E$.

Finally as $M_0 \normlt G_0$ and $G_0$ is transitive on $V^*$, it follows by
Lemma~\ref{lem:DD} that $G_0$ permutes the the elements of the set $\M{S}$ transitively.
Now suppose $\sigma \in G_0$ maps $S_i$ to $S_j$ for some $1 \leq i,j \leq k$ and $i
\neq j$. Then $\{u,v\} \in E\Gamma_i \imply v-u \in S_i$ and so $(v-u)^{\sigma} =
v^{\sigma} -u^{\sigma} \in S_j$. Thus $\{u^{\sigma}, u^{\sigma}\} \in E\Gamma_j$ and it
follows that $E\Gamma_i^{\sigma} =E\Gamma_j$. Since $T$ fixes each $E\Gamma_i$ setwise,
it follows that $G=T \rtimes G_0$ acts on $\mathcal{E}$. As $G_0$ is transitive on the
sets $\M{S}$, we have $G$ acting transitively on $\mathcal{E}$. Finally, as $M$ fixes
each $E\Gamma_i$ set-wise and $G$ is transitive on $\M{E}$, the result on $(M,G,V,
\M{E})$ being a homogeneous factorisation with $M$-arc-transitive factors $\Gamma_i$
follows.

Conversely, suppose $(M,G,V, \M{E})$ is an arc-transitive homogeneous factorisation of
$K_{p^R}=(V,E)$ of index $k$, with factors $\Gamma_i$ ($1 \leq i \leq k$), $G \leq
A\Gamma L(1,p^R)$ acting $2$-transitively on $V$ and $|M|$ even.

Then $M_0$ and $G_0$ (and hence $M=T \rtimes M_0$ and $G=T \rtimes G_0$) can be
constructed as subgroups of $\Gamma L(1,p^R)$ by checking that all the conditions in
Theorem~\ref{thm1} are satisfied.

Suppose now that the groups $M$ and $G$ are constructed using Theorem~\ref{thm1}. Then
as the factors are $M$-arc-transitive Cayley graphs $\Gamma_i =Cay(V,S_i)$, it follows
that the set $\M{S} =\{S_i \mid 1 \leq i \leq k\}$ are all $M_0$-orbits in $V^*$; now
using Construction~\ref{generic}, both $\M{S}$ and $\Gamma_i$ can be constructed for a
given $M_0$. Lastly, the partition $\M{E}$, which are the edge sets of $\Gamma_i$, are
also given by Construction~\ref{generic}. \epf

Even though in principle (in light of Theorem~\ref{thm2}), Theorem~\ref{thm1} and
Construction~\ref{generic} enable us to construct all possible arc-transitive homogeneous
factorisations $(M,G,V,\M{E})$ where $G$ is an one-dimensional affine 2-transitive
group, explicit classification is still far from complete and seems difficult. Indeed,
we are only able to give (later in this chapter) some families of examples. Thus the
question of complete explicit classification is still open.

\bpro \label{c5:prob} \rm{Give an explicit classification of all arc-transitive homogeneous
factorisations $(M,G,V,\M{E})$ with $G =T \rtimes G_0$ and $M=T \rtimes M_0$, where $G_0$, $M_0$
are as in Theorem~\ref{thm1}. That is to say, find explicitly all solutions to the divisibility
conditions $(1)$ to $(11)$ in Theorem~\ref{thm1}. } \epro

\section{Some explicit examples} \label{c5:const}
We now use the results in Theorem~\ref{thm1} to derive necessary and sufficient
conditions for certain subgroups $G_0$ and $M_0$ of $\Gamma L(1,p^R)$ to satisfy
properties that will enable us to construct interesting examples of arc-transitive
homogeneous factorisations. Note that in the following, $V=\F{p^R}$ and $V^*=V \setminus
\{\0\}$.

\bl \label{gpaley} Let $p$ be any prime number, $R$ a positive integer, and $k \geq 2$
be an integer which divides $p^R-1$ such that if $p$ is odd, then $\frac{p^R-1}{k}$ is
even. Let $G_0=\brac{\omb}= GL(1,p^R)$ and $M_0 =\brac{\omb^k}$ be subgroups of $\Gamma
L(1,p^R)$ acting on $V^*$. Then $M_0$ and $G_0$ are in standard form and $M_0$ is a
normal subgroup of $G_0$. Moreover the following hold. \bnum
\item[$1.$] $G_0$ is transitive on $V^*$.
\item[$2.$] $M_0$ has $k$ orbits of equal length $\frac{p^R-1}{k}$ in $V^*$, and if $p$ is odd,
then $|M_0|$ is even. \enum \el \bpf The assertions about $M_0$ and $G_0$ being in
standard form and $M_0 \normlt G_0$ follow immediately by checking through conditions of
Theorem~\ref{thm1}.

Since $|M_0| =\frac{p^R-1}{k}$, it follows by the assumption of the lemma that if $p$ is
odd, then $|M_0| =\frac{p^R-1}{k}$ is even. Now part (1) and the rest of part (2) of the
above lemma follow since $G_0$ is regular on $V^*$. \epf

\bl \label{newcon} Let $p \equiv 3$ $(mod \ 4)$ be a prime and $h$ be an odd positive
integer. Let $G_0=\brac{\omb^2,\omb \al}$ and $M_0 =\brac{\omb^{4h},\omb^h \al}$ be
subgroups of $\Gamma L(1,p^R)$ acting on $V^*$. Suppose also that $2h \mid (p-1)$ and $R
> 1$ is even. Then the following hold. \bnum
\item[$1.$] $G_0$ and $M_0$ are in standard form and $M_0$ is a normal subgroup of
$G_0$.
\item[$2.$] $G_0$ is transitive on $V^*$, and $M_0$ has $2h$ orbits
of equal length $(p^R-1)/2h$.
\item[$3.$] $|M_0|$ is even. \enum
\el \bpf (1) We shall go through conditions listed in Theorem~\ref{thm1}. For ease of notation, we
will denote by C1 to C10, the conditions (1) to (10) listed in Theorem~\ref{thm1} respectively.
(Following the notations of Theorem~\ref{thm1}, we have $(d,e,s) = (2,1,1)$ and $(d_1,e_1,s_1) =
(4h,h,1)$.)

(\textbf{C1}) Since $R$ is even, condition (C1) follows immediately by the assumption $2h
\mid (p-1)$.

(\textbf{C2, C4, C5 and C7}) Trivially satisfied.

(\textbf{C3}) Here, we need $4 \mid ((p^R-1)/(p-1))$. Now $(p^R-1)/(p-1) = p^{R-1}+
p^{R-2}+ \ldots + p + 1$, is a sum of $R$ odd integers (recall that $p$ is an odd
prime). Since $p \equiv 3$ (mod 4) and $R$ is even, it follows that $4 \mid
((p^R-1)/(p-1))$ holds.

(\textbf{C6}) Here we need $2 \mid (h-1)$. As $h \geq 1$ is an odd integer. The
condition is satisfied.

(\textbf{C8}) For this case we need $4h \mid (h(p-1) + (p^R -p))$. Now
\[ h(p-1) + (p^R -p) = (p-1)(h + p^{R-1} + p^{R-2} +\ldots + p^2 + p).\] Recall that $R$ is
even and $h$ odd, so the expression $h + p^{R-1} + p^{R-2} + \ldots + p^2 + p$ is a sum
of $R$ odd integers and hence is even. It follows that $2h \mid (p-1)$ (as given in the
assumption) and $2 \mid (h + p^{R-1} + p^{R-2} + \ldots + p^2 + p)$, and so $4h \mid
(h(p-1) + (p^R -p))$.

Thus conditions (C1) to (C8) of Theorem~\ref{thm1} are satisfied and part (1) of the lemma holds.

(2) It is easy to check that conditions (C9) and (C10) of Theorem~\ref{thm1} are
satisfied for $G_0=\brac{\omb^2, \omb \al}$. So $G_0$ is transitive on $V^*$.

From Lemma~\ref{mintran} (see also the last part of Theorem~\ref{thm1}), $M_0$ has $k =
4h/c$ orbits of length $\frac{p^R-1}{k}$ in $V^*$ where $c$ is such that:
\newcounter{count}
\begin{list} {{\rm (\alph{count})}}{\usecounter{count}}
\item $4 \mid (p^c-1)/(p-1)$ and
\item $4 \nmid (p^{c'}-1)/(p-1)$ for all $c' < c $.
\end{list} We shall now proceed to determine the possible values for $c$ (and hence $k$).
Clearly $c \geq 2$. If $c'=2$, then $(p^{c'}-1)/(p-1)=p+1$. Since $p \equiv 3$ (mod 4),
we have $4 \mid (p+1)$, and so it follows that $c =2$, and hence $k=2h$ (since $k
=4h/c$).

(3) Now $|M_0| =R(p^R-1)/2h$. Since $2h \mid (p-1)$, we have that $2h \mid (p^R-1)$ and
so $\frac{p^R-1}{2h}$ is an integer. As $R$ is even, it follows that $|M_0| = R
(\frac{p^R-1}{2h})$ is also even. \epf

Recall that in Section~\ref{c3:examples}, we make use of the results in
Lemmas~\ref{gpaley} and \ref{newcon} to construct two examples of arc-transitive
homogeneous factorisations $(M,G,V,\M{E})$ of complete graphs where the group $G \leq
A\Gamma L(1,p^R)$. They are the infinite examples of arc-transitive homogeneous
factorisations (where $G$ is the one-dimensional affine group) that correspond to the
generalised Paley partition $\M{E}_{GP}$ and the twisted generalised Paley partition
$\M{E}_{TGP}$ (see Propositions~\ref{prop:gpaley} and \ref{prop:tgpaley}).\footnote{Note
that the approach in the proofs for Propositions~\ref{prop:gpaley} and \ref{prop:tgpaley}
are essentially the same as that found in the proof of Theorem~\ref{thm1}.} The
arc-transitive factors are respectively called the generalised Paley graphs and twisted
generalised Paley graphs. For completeness, we shall repeat the definitions of these
graphs here.

\vsl \noindent \textbf{(Generalised Paley graph)} \ Let $V=\F{p^R}$ and $V^*=V \setminus \{\0\}$.
Let $k \geq 2$ be an integer which divides $p^R-1$ and is such that, if $p$ is odd, then
$\frac{p^R-1}{k}$ is even. Let $\om$ be a fixed primitive element in $V$. The graph $\GP$ is the
Cayley graph $Cay(V,S)$ with connecting set $S = \brac{\om^k} =\{ 1, \om^k, \om^{2k}, \ldots,
\om^{p^R -1 -k} \} \varsubsetneq V^*$ and is called a \textit{generalised Paley graph} with
respect to $V$.

\vsl \noindent \textbf{(Twisted generalised Paley graph and partition)} \ Again, let $V$ and $V^*$
be as previously defined and $\om$ be a primitive element in $V$. Let $R$ be an even integer and
$p \equiv 3$ (mod 4). Let $h \geq 1$ be an odd integer such that $2h \mid (p-1)$. Let $S_0
=\brac{\om^{4h}} = \{1, \om^{4h}, \om^{8h}, \ldots, \om^{p^R -1 -4h} \} \varsubsetneq V^*$
(observe that $R$ even and $2h \mid (p-1)$ imply that $4h \mid (p^R-1)$). The \textit{twisted
generalised Paley graph}, denoted as $TGPaley(p^R,\frac{p^R-1}{2h})$, with respect to $V$ is the
Cayley graph $Cay(V, S_0 \cup \om^{3h} S_0)$.

\vsl

We also recall that both of the above graphs are undirected, and a generalised Paley
graph $\GP$ admits $T \rtimes \brac{\omb^k,\al}$ as an subgroup of automorphisms; while a
twisted generalised Paley graph $TGPaley(p^R,\frac{p^R-1}{2h})$ admits $T \rtimes
\brac{\omb^{4h}, \omb^h \al}$ as a subgroup of automorphisms.

Furthermore, if $k=2$, then $\GPq$ is the familiar Paley graph with full automorphism
group $T \rtimes \brac{\omb^2, \al}$ (see for example \cite{Peisert2001}). If $h=1$,
then the $TGPaley(p^R,\frac{p^R-1}{2h})$ is the self-complementary, arc-transitive
graph, called the $\M{P}^*$-graph, studied by Peisert in \cite[Section 3]{Peisert2001}.
Also in \cite{Peisert2001}, it was proved that in general, except for one case (see
Lemma~\ref{c5:small} below), $GPaley(p^R, \frac{p^R-1}{2})$ (note that $k=2$) is not
isomorphic to $TGPaley(p^R, \frac{p^R-1}{2})$ (where $h=1$).

\bl \label{c5:small} {\rm \cite[Lemma 6.4]{Peisert2001}} \ Suppose $(R,p,k,h)=(2,3,2,1)$. Then
$GPaley(3^2,4) \cong TGPaley(3^2,4)$. \el

We shall later show that the for $k \geq 3$ and $h > 1$ ($h$ is odd), the twisted generalised
Paley graphs are not isomorphic to the generalised Paley graphs.

\section{Non-isomorphism} \label{c5:noniso}
Let $\Gamma' = TGPaley(p^R,\frac{p^R-1}{2h}) =Cay(V,S')$ be a twisted generalised Paley graph
where $S' = \brac{\om^{4h}} \cup \om^{3h} \brac{\om^{4h}}$ and $p \equiv 3$ (mod 4). Recall that
from the definition of twisted generalised Paley graphs, $h \geq 1$ is an odd integer, $2h \mid
(p-1)$ and $R$ is even (so $R \geq 2$). Also, $\Gamma'$ has valency equal to $|S'|
=\frac{p^R-1}{2h}$ (see Proposition~\ref{prop:tgpaley}). Suppose $\Gamma'$ is isomorphic to a
generalised Paley graph $\Gamma=GPaley(p^R,\frac{p^R-1}{k})$. Then we must have $\frac{p^R-1}{2h}
=\frac{p^R-1}{k}$ and so it follows that $2h=k$. As $2h \mid (p-1)$, it follows that $k \mid
(p-1)$. Now we can invoke Theorem~\ref{thm:no2} which we proved in Chapter~\ref{c6}.


Suppose $\Gamma'=TGPaley(p^R,\frac{p^R-1}{2h}) \cong GPaley(p^R,\frac{p^R-1}{k})=\Gamma$
where $k=2h$ and $R$ is even (since $TGPaley(p^R,\frac{p^R-1}{2h})$ is defined for even
$R$). The case where $k=2$ has already been resolved by Peisert in \cite{Peisert2001}.
Thus in our proof we shall assume that $k \geq 3$. Then since $k$ is even and $k \mid
(p-1)$, we have $k >3$ and $p >3$. We will use Theorem~\ref{thm:no2} to arrive at a
contradiction, and hence show that, for $p>3$, $k>3$, $\Gamma' \ncong \Gamma$.

Note that $Aut(\Gamma') \geq N_{Aut(\Gamma')}(T) \geq T \rtimes \brac{\omb^{4h}, \omb^h
\al}$. Since we assume $\Gamma \cong \Gamma'$, then $T \rtimes \brac{\omb^{4h}, \omb^h
\al}$ is permutationally isomorphic to a subgroup of $Aut(\Gamma) = T \rtimes
\brac{\omb^{2h}, \al}$ (when applying Theorem~\ref{thm:no2} since $k=2h \mid (p-1)$). In
particular, $\brac{\omb^{4h}, \omb^h \al}$ is permutationally isomorphic to a subgroup of
$\brac{\omb^{2h}, \al}$ in its action on $V^*$.

\bt \label{luvsandra} \ For $R$ even, $p \equiv 3$ $($mod $4)$, and $h$ odd such that $2h \mid
(p-1)$, $GPaley(p^R, \frac{p^R-1}{2h}) \ncong TGPaley(p^R, \frac{p^R-1}{2h})$ except for the case
$(R,p,h)=(2,3,1)$, when we have $GPaley(3^2, 4) \cong TGPaley(3^2, 4)$. \et \bpf By
Lemma~\ref{c5:small}, $GPaley(9, 4) \cong TGPaley(9, 4)$. Assume now that $(R,p,h) \neq (2,3,1)$
and that $GPaley(p^R, \frac{p^R-1}{2h}) \cong TGPaley(p^R, \frac{p^R-1}{2h})$. Then as discussed
above, the stabiliser $\brac{\omb^{4h},\omb^h \al}$ is permutation isomorphic to a subgroup $K$ of
$\brac{\omb^{2h}, \al}$ in its action on $V^*$. Moreover we have $p >3$ and hence $\omb$ is an odd
permutation. Since $h$ is odd, $\omb^{4h}$ and $\omb^{2h}$ are even permutations while $\omb^h$ is
an odd permutation.

Suppose first that $\al$ is an odd permutation. Then $\brac{\omb^{4h},\omb^h \al}$
consists entirely of even permutations. It follows that $K$ is contained in $A_{p^R} \cap
\brac{\omb^{2h}, \al} = \brac{\omb^{2h}, \al^2}$. Now $\brac{\omb^{2h}, \al^2}$ has
order $\frac{R(p^R-1)}{4h}$ which equals $|K|$. Hence $K = \brac{\omb^{2h}, \al^2}$.
Consider the case where $R =2$. Then $K= \brac{\omb^{2h}, \al^2}$ becomes
$\brac{\omb^{2h}}$ (since $\al^2=1$ for $R=2$) which is abelian. Thus
$\brac{\omb^{4h},\omb^h \al}$ (permutationally isomorphic to $K$) is also abelian and so
we must have $\omb^{4h} . \omb^h \al = \omb^h \al . \omb^{4h}$. It follows that \be
(\omb^{4h})^{\omb^h \al} & = & \omb^{4h} \\
\omb^{4hp} &=&  \omb^{4h} \\
\omb^{4h(p-1)} &=&  1. \ee This implies that $(p + 1) \mid 4h$. Since $2h \mid (p-1)$, we have $4h
\mid 2(p-1)$. So $(p+1) \mid 2(p-1)$, which is only possible when $p=3$. However, this contradicts
the fact that we have $p >3$.

Next consider the case when $R >2$. Let $r$ be a $\ppd$ of $p^R-1$ (note that $r$ is
odd, and $r$ exists since $p>3$ and $R >2$, see Definition~\ref{def:ppd}). Since
$|\brac{\omb^{2h}}| = \frac{p^R-1}{2h}$ is divisible by $r$ and $r \nmid
|\brac{\al^2}|=R/2$ (because $r >R$, see \cite[Remark 1.1 (a)]{GPPS}), it follows that
$K= \brac{\omb^{2h}, \al^2} = \brac{\omb^{2h}} \rtimes \brac{\al^2}$ has a unique cyclic
Sylow $r$-subgroup $P$ such that $P \leq \brac{\omb^{2h}}$. Thus $\brac{\omb^{4h},\omb^h
\al}$ also has a unique (cyclic) Sylow $r$-subgroup contained in $\brac{\omb^{4h}}$, say
$\brac{\omb^{4hj}}$. Furthermore, $\omb^h \al \in \brac{\omb^{4h},\omb^h \al}$
conjugates a generator of $\brac{\omb^{4hj}}$ to its $p$th power. Hence some element of
$K$ also conjugates an element of $P$ to its $p$th power. Let $P =\brac{\omb^{\ell}}$
where $\frac{p^R-1}{\ell} =r^b$ for some $b \geq 1$. Let $x \in K$ be such that
$(\omb^{\ell})^x = \omb^{\ell p}$. Now $x = \omb^{2hi} \al^{2j}$ for some $i,j$ and
$1\leq j \leq \frac{R}{2}$. So \be
(\omb^{\ell})^x & = & \omb^{\ell p^{2j}} \\
\omb^{\ell p} & = & \omb^{\ell p^{2j}} \\
\omb^{\ell p(p^{2j-1}-1)} & = & 1. \ee It follows that $(p^R-1) \mid \ell p(p^{2j-1}-1)
\imply (p^R-1) \mid \ell (p^{2j-1}-1) \imply \frac{p^R-1}{\ell} =r^b \mid (p^{2j-1}-1)$.
However, this contradicts the fact that $r \nmid (p^{2j-1}-1)$ (by definition of $\ppd$,
since $2j-1 \leq R-1$).

Now suppose that $\al$ is an even permutation. Then $\brac{\omb^{2h}, \al}$ consists
entirely of even permutations. However, $\omb^h \al \in \brac{\omb^{4h},\omb^h \al}$ is
an odd permutation, and as $\brac{\omb^{4h},\omb^h \al}$ is permutationally isomorphic
to $K \leq \brac{\omb^{2h}, \al}$ for which all elements are even permutations, we again
have a contradiction. Thus $\Gamma' \ncong \Gamma$. \epf

To conclude this chapter, we have shown in Theorem~\ref{luvsandra} that the class of generalised
Paley graphs and the class of twisted generalised Paley graphs have exactly one graph in common.
This implies that the corresponding partitions (the generalised Paley partition $\M{E}_{GP}$ and
the twisted generalised Paley partition $\M{E}_{TGP}$) give rise to two infinite and different
families of arc-transitive homogeneous factorisations of complete graphs. It would be interesting
to find if there are any other infinite examples of arc-transitive homogeneous factorisations
$(M,G,V,\M{E})$ where $G \leq A\Gamma L(1,p^R)$ (see Problem~\ref{c5:prob}).

\chapter*{Appendix: Computations using MAGMA} \label{App}

\addcontentsline{toc}{chapter}{\textbf{Appendix: Computations using MAGMA}}
\markboth{Appendix}{Appendix}

This chapter contains some \textsc{Magma} codes that
were used for computations required in the thesis. Although not all the codes used for
performing our computations are written here, we provide examples that sufficiently
illustrate how the rest of results may be obtained (see especially Sections \textbf{A3} and
\textbf{A4}). We also note that the codes presented
here are by no means the most efficient, but merely written to serve their purposes.
Lastly, the version of \textsc{Magma} being used is V2.9-21.

\section*{A1: MAGMA code for constructing an orbital graph}

To construct an $M$-orbital graph given a permutation group $M$
(or \verb+PermGp+ in \textsc{Magma}).

\begin{verbatim}
> Stab_M:=Stabilizer(M,1);
//Construct the point-stabiliser Stab_M of "1" in M.
//(Note that M acts on the set {1,2,....,n} where n is the degree of M.)

> Set:=Orbit(M,2);
//Construct the orbit of Stab_M containing the point "2".

> O:=OrbitalGraph(M,1,Set);
//Construct an M-orbital graph containing the edge {1,2}.

> A:=AutomorphismGroup(O);
//Determine full automorphism group of the orbital graph O.
\end{verbatim}

\newpage

\section*{A2: MAGMA code for Remark~\ref{rk:psl28}}
\begin{verbatim}
> Num:=NumberOfPrimitiveGroups(28);
> for i:= 1 to Num by 1 do;
    > G:=PrimitiveGroup(28,i);
    > if #G eq 504 then;
        > print i;
        > end if;
    > end for;
2
//Thus there is only one primitive group: PrimitiveGroup(28,2),
//of degree 28 with order 504.
> G:=PrimitiveGroup(28,2);
> StabG:=Stabilizer(G,1);
> Orbits(StabG);
[
    GSet{ 1 },
    GSet{ 2, 7, 8, 9, 12, 13, 15, 19, 27 },
    GSet{ 3, 6, 10, 11, 17, 22, 23, 24, 25 },
    GSet{ 4, 5, 14, 16, 18, 20, 21, 26, 28 }
]
> Set:=Orbit(StabG,2);
> O:=OrbitalGraph(G,1,Set);
> A:=AutomorphismGroup(O);
> #A;
504
//Evaluating the order of the group A.
> Num:=#Subgroups(A: IsRegular:=true);
//Compute the number of regular subgroups on 28 points in A.
> Num;
0
//Thus, there is no regular subgroup on 28 points in A.
\end{verbatim}


\newpage

\section*{A3: MAGMA codes used for computing Tables~\ref{tab:TK2} and \ref{tab:TK1}}

\subsection*{\underline{A3-1}}
The \verb+PFind(G,H)+ function below enables us to ``locate" either
$PSL(2,5)$ or $PSL(2,3)$ in $\oline{G_0} = P\Gamma L(2,q)$ where $G_0$ is as in case $2(g)$ of
Theorem~\ref{2_trxclass} (see Section~\ref{c4:SL}). Note that \verb+G+$:=P\Gamma L(2,q)$ and
\verb+H+$:=PSL(2,5)$ or $PSL(2,3)$. Also, if the function returns a warning ``\verb+More than one+
\verb+group+", then our \verb+PFind(G,H)+ function fails. For our investigation in
Section~\ref{c4:SL}, the function always returns the desired group we are seeking.


\begin{verbatim}
> function PFind(G,H);
//G:=PGammaL(2,q) and H:=PSL(2,5) or PSL(2,3).
> j:=1;
> h:=[];
> k:=[];
> Sub:=Subgroups(G: OrderEqual:= #H);
//Locating subgroups of G with order equal |H|.
> for i:= 1 to #Sub by 1 do;
    > k[i]:=Sub[i]`subgroup;
    > ISM:=IsIsomorphic(k[i],H);
//Determine if a subgroup k[i] is isomorphic to H=PSL(2,5) or PSL(2,3).
    > if ISM eq true then;
        > Norm:=Normalizer(G,k[i]);
//Compute the normaliser Norm of H=PSL(2,5) or PSL(2,3) in G=PGammaL(2,q).
        > TN:=IsTransitive(Norm);
//Check if Norm is transitive on the set of 1-spaces of V=V(2,q).
        > if TN eq true then;
            > h[j]:=k[i];
            > j:=j+1;
            > end if;
        > end if;
    > end for;
> if j eq 2 then;
    > return h[j-1];
//If there is only one group (up to conjugacy classes) in G that is
//isomorphic to H, and that the normaliser of H in G is transitive on
//the set of 1-spaces, then the function returns H as a subgroup of G.
    > else
    > return "More than one group";
//If there are 2 or more groups in G isomorphic to H with the desired
//property described earlier, then the function returns a warning.
    > end if;
> end function;
\end{verbatim}

\subsection*{\underline{A3-2}}
We give examples on how some results in Tables~\ref{tab:TK2} and \ref{tab:TK1} are
obtained by using the \verb+PFind+ function.

\vsl
\noindent \textbf{Example A3-2(i)}: We will show that $PSL(2,5)$ in $P\Gamma L(2,9)$ is
transitive on the set of 1-spaces of $V=V(2,9)$ (see Table~\ref{tab:TK2}).
\begin{verbatim}
> g:=PGammaL(2,9);
> h:=PSL(2,5);
> D:=PFind(g,h);
> ISM:=IsIsomorphic(D,PSL(2,5));
//Check to see if D is isomorphic to PSL(2,5).
> ISM;
true
> IsTransitive(D);
true
//PSL(2,5) in PGammaL(2,9) is transitive in the set of 1-spaces.
\end{verbatim}

\vsl
\noindent \textbf{Example A3-2(ii)}: We will show that $\oline{M_0} \cong \Z_2 \times \Z_2$ in
$PGL(2,5)$ has 3 orbits of length 2 in the 1-spaces of $V=V(2,5)$ (see line 1 of
Table~\ref{tab:TK1}). Note that we utilised the fact that $\oline{M_0}$ has orbits of
equal length in the set of 1-spaces (since it is normal in $\oline{G_0}$ which acts
transitively on the 1-spaces).
\begin{verbatim}
> g:=PGL(2,5);
> h:=PSL(2,3);
> D:=PFind(g,h);
> D;
> ISM:=IsIsomorphic(D,PSL(2,3));
//Check to see if D is isomorphic to PSL(2,3).
> ISM;
true
> N:=Normalizer(g,D);
//Construct the normaliser N of PSL(2,3) in PGL(2,5).
> Sub:=Subgroups(N: IsElementaryAbelian:=true, OrderEqual:=4);
> Sub;
Conjugacy classes of subgroups
------------------------------
[1]     Order 4            Length 1
        Permutation group acting on a set of cardinality 6
        Order = 4 = 2^2
            (2, 3)(4, 5)
            (1, 6)(4, 5)
[2]     Order 4            Length 3
        Permutation group acting on a set of cardinality 6
        Order = 4 = 2^2
            (1, 6)(2, 4)(3, 5)
            (2, 3)(4, 5)
> V4_1:=Sub[1]`subgroup;
> V4_2:=Sub[2]`subgroup;
//Setting V4_1 and V4_2 to be the above two elementary abelian
//subgroups of order 4.
> Orbits(V4_1);
[
    GSet{ 1, 6 },
    GSet{ 2, 3 },
    GSet{ 4, 5 }
]
> Orbits(V4_2);
[
    GSet{ 1, 6 },
    GSet{ 2, 3, 4, 5 }
]
//Thus the group V4_1 is the desired group and has 3 orbits of length 2.
\end{verbatim}

\newpage

\section*{A4: MAGMA codes used for computing Tables~\ref{t25} and \ref{t23}}

\subsection*{\underline{A4-1}}
The \verb+Find(G,H,q)+ function below enables us to ``locate"
either $T \rtimes SL(2,5)$ or $T \rtimes SL(2,3)$ in $G= T \rtimes G_0 \leq A\Gamma L(2,q)$
where $G_0$ is as in case $2(g)$ of Theorem~\ref{2_trxclass} (see also Section~\ref{c4:SL}).
Note that if the function returns ``\verb+This function does not work for+
\verb+your case+", then we will have to look at the case separately. However, for our investigation
in Section~\ref{c4:SL}, the \verb+FindGroup+ function always return the group that we are
seeking.
\begin{verbatim}
> function Find(G,H,q);
//To determine T:SL(2,5) or T:SL(2,3) in AGammaL(2,q), where T is
//the translation group acting on V(2,q). Note that the input is
//G:= AGammaL(2,q) and H:=SL(2,5) or SL(2,3). (Note that G and H must
//be permutation groups on V(2,q)).
> StG:=Stabilizer(G,1);
> Set:=Orbit(StG,2);
> j:=0;
> n:=#H * q^2;
//Computing the order of T:H (where H=SL(2,5) or SL(2,3)) in G.
//Also, we let n :=|T:H|.
> Sub:=Subgroups(G: OrderEqual:=n);
//Finding subgroups (up to conjugacy classes) of order n=|T:H| in G.
> for i:= 1 to #Sub by 1 do;
    > A:=Sub[i]`subgroup;
    > StA:=Stabilizer(A,1);
    > T:=IsIsomorphic(StA,H);
//Determine if the subgroup A is isomorphic to the desired group T:H
//in AGammaL(2,q).
    > if T eq true then;
        > N:=Normalizer(StG,StA);
        > TN:=IsTransitive(N,Set);
//Determine if the the normaliser of A in G is 2-transitive.
//If yes, then j will acquire a value of 1 or more.
        > if TN eq true then;
            > j:=j+1;
            > Q:=A;
            > end if;
        > end if;
    > end for;
> if j eq 1 then;
    > return Q;
//If j=1, then there is only 1 group (T:SL(2,5) or T:SL(2,3)) in
//AGammaL(2,q) with the desired property. The function will then
//return the group that we want, that is, Q=T:SL(2,5) or T:SL(2,3).
    > else
//Note that if j is greater than 1, then this means that there
//are more than 1 group with the desired properties. In that case,
//the function will return a warning (see below).
    > return "This function does not work for your case.";
    > end if;
> end function;
\end{verbatim}

\subsection*{\underline{A4-2}} The \verb+Gen(G,H)+ function below enables us to generate
all possible $M$ (and $M_0$) in $G$ (and $G_0$). Note that $G=T \rtimes G_0 \leq A\Gamma
L(2,q)$, where $G_0$ is as in case $2(g)$ of Theorem~\ref{2_trxclass} (see
Section~\ref{c4:SL}). Also, $M=T \rtimes M_0$ is normal in $G$ and $M_0$ contains either
$SL(2,5)$, $SL(2,3)$ or $Q_8$ (see Lemmas~\ref{sl25} and \ref{sl23}). The function also
determines the orders of $M$ and $M_0$, and computes the number of $M_0$-orbits in $V^*$
(where $V=V(2,q)$). See Tables~\ref{t25} and \ref{t23}.
\begin{verbatim}
> function Gen(G,H);
//To determine all possible M_0 and M. Note that the input is H:=SL(2,5)
//or SL(2,3), and G is the normalizer of T:SL(2,5) or T:SL(2,3) in
//AGammaL(2,q) (where T:SL(2,5) or T:SL(2,3) in AGammaL(2,q) is determined
//by the Find(G,H,q) function).
> j:=1;
> h:=[];
> sh:=[];
> k:=[];
> sk:=[];
> w:=[];
> StG:=Stabilizer(G,1);
> Set:=Orbit(StG,2);
> Sub:=Subgroups(G: IsTransitive:=true);
//Find all transitive subgroups of G.
> for i:= 1 to #Sub by 1 do;
    > k[i]:=Sub[i]`subgroup;
    > sk[i]:=Stabilizer(k[i],1);
    > T:=IsTransitive(k[i],2);
//Get rid of the subgroups of G that are 2-transitive.
    > if T eq false then;
        > TC:=IsCyclic(sk[i]);
//Get rid of the subgroups of G whose point-stabilisers are scalars.
        > if TC eq false then;
            > Norm:=Normalizer(G,k[i]);
            > sNorm:=Stabilizer(Norm,1);
            > TN:=IsTransitive(Norm,2);
//Get rid of the subgroups of G whose normalisers (sNorm) in G
//are not 2-transitive.
            > if TN eq true then;
                > SubsNorm:=Subgroups(sNorm: OrderEqual:=#H);
//Determine the subgroups of sNorm with same order as the group H,
//where H=SL(2,5) or SL(2,3).
                > for y:= 1 to #SubsNorm by 1 do;
                    > A:=SubsNorm[y]`subgroup;
                    > TNN:=IsIsomorphic(A,H);
//Determine if sNorm contains subgroups isomorphic to H.
                    > if TNN eq true then;
                        > h[j]:=k[i];
                        > sh[j]:=sk[i];
                        > w[j]:=i;
//Note that w[j]=i tells us the "group number" of k[i] in
//"Sub:=Subgroups(G: IsTransitive:=true);" (see Example 4(i)).
                        > j:=j+1;
                        > end if;
                    > end for;
                > end if;
            > end if;
        > end if;
    > end for;
> for c:= 1 to j-1 by 1 do;
        > print c, "(i)=", w[c], " Order of M=", #h[c],
" Order of M_0=", #sh[c], " No. of orbits in V^*=", #Orbits(sh[c],Set);
    > end for;
> return "================";
> end function;
\end{verbatim}

\subsection*{\underline{A4-3}} The \verb+GPaley(p,R,k)+ function below enables us to
construct a generalised Paley graph $\GP$ with given parameters (integers) $p$, $R$ and
$k$.
\begin{verbatim}
> function GPaley(p,R,k);
> G:=AGL(1,p^R);
> StG:=Stabilizer(G,1);
> Set:=Orbit(StG,2);
> Order:=(p^R-1)/k;
//This the the order of the group <w^k> (see definition of
//generalised Paley graphs).
> j:=0;
> h:=[];
> sh:=[];
> k:=[];
> sk:=[];
> O:=[];
> SubG:=Subgroups(G: IsTransitive:=true);
//Find all subgroups of G=AGammaL(1,p^R) that are transitive.
> for i:= 1 to #SubG by 1 do;
    > k[i]:=SubG[i]`subgroup;
    > sk[i]:=Stabilizer(k[i],1);
    > if #sk[i] eq Order then;
//Check if the point stabiliser of k[i] (a transitive subgroup of G)
//has the same order as that of <w^k>.
        > TS:=IsSemiregular(sk[i],Set);\
        > if TS eq true then;
//Check if sk[i] acts semiregularly on V^* (since <w^k> acts
//semiregularly on V^*). If yes, then j will acquire a value of 1 or more
            > h[j]:=k[i];
            > sh[j]:=sk[i];
            > SetO:=Orbit(sk[i],2);
            > O[j]:=OrbitalGraph(k[i],1,SetO);
//Construct the orbital graph using the group k[i].
            > j:=j+1;
            > end if;
        > end if;
    > end for;
> if j eq 1 then;
    > return O[j-1];
//If j=1 (that is, only one graph is constructed and there is only
//one group in AGammaL(1,p^R) whose point stabiliser has the same
//properties as that of <w^k>) then the function returns the
//desired graph.
    > else
    > return "More than 2 graphs.";
//This happens only when the function fails. However for our
//investigation, we always get the desired graph.
    > end if;
> end function;
\end{verbatim}

\subsection*{\underline{A4-4}} We will show by giving examples on how the functions \verb+Find(G,H,q)+,
\verb+Gen(G,H)+ and \verb+GPaley(p,R,k)+ are used to give us the results in
Tables~\ref{t25} and \ref{t23}.

\vs
\noindent \textbf{Example A4-4(i)}: Line $(1)$ of Table~\ref{t25} ($\Gamma_i \cong G(9^2,2)$).
\begin{verbatim}
> g:=AGammaL(2,9);
> h:=Stabilizer(ASL(2,5),1);
> D:=Find(g,h,9);
> #D;
9720
> N:=Normalizer(g,D);
//Compute the normaliser of T:SL(2,5) in AGammaL(2,9).
> Gen(N,h);
1 (i)= 100  Order of M= 9720  Order of M_0= 120  No. of orbits in V^*= 2
//Note that the "(i) = 100" tells us that M is the #100 group
//in "Subgroups(N: IsTransitive:=true)".
2 (i)= 107  Order of M= 19440  Order of M_0= 240  No. of orbits in V^*= 2
3 (i)= 111  Order of M= 38880  Order of M_0= 480  No. of orbits in V^*= 2
================
> Sub:=Subgroups(N: IsTransitive:=true);
> M111:=Sub[111]`subgroup;
//We are letting M111 be the #111 group in
//"Subgroups(N: IsTransitive:=true)".
> M100:=Sub[100]`subgroup;
> M107:=Sub[107]`subgroup;
> sM100:=Stabilizer(M111,1);
> sM100:=Stabilizer(M100,1);
> sM107:=Stabilizer(M107,1);
> Set111:=Orbit(sM111,2);
> Set100:=Orbit(sM100,2);
> Set107:=Orbit(sM107,2);
> OG111:=OrbitalGraph(M111,1,Set111);
> OG100:=OrbitalGraph(M100,1,Set100);
> OG107:=OrbitalGraph(M107,1,Set107);
//Constructing the orbital graphs of M111, M100 and M107,
//that is, the M-orbital graphs where M=M111, M100 or M107.
> IsIsomorphic(OG100,OG107);
true
> IsIsomorphic(OG111,OG107);
true
//So all 3 orbital graphs are isomorphic to one another.
//The next step is to show that the orbital graphs just constructed
//are isomorphic to a corresponding TGPaley graph. We note that this
//TGPaley graph is self-complementary (since index k=2) and admits
//T:<w^4,wa> as an arc-transitive subgroup of automorphisms (see relevant
//section on the TGPaley graphs). Also <w^4,wa> has a cyclic normal
//subgroup <w^4> of order 20 such that the quotient <w^4,wa>/<w^4> is a
//cyclic subgroup of order 4. Moreover, <w^4,wa> has 2 orbits of equal
//length in V^*. We will make use of all these facts in our construction
//of the required TGPaley graph.
> G:=AGammaL(1,81);
> subG:=Subgroups(G: OrderEqual:=6480);
//Note that |T:<w^4,wa>|=81*80=6480.
> for i:=1 to #subG by 1 do;
    > A:=subG[i]`subgroup;
    > sA:=Stabilizer(A,1);
    > TC:=IsCyclic(sA);
    > if TC eq false then;
//Since <w^4,wa> is not cyclic.
        > num:=#Orbits(sA);
        > if num eq 3 then;
//Since <w^4,wa> has 2 orbits in V^* (or 3 orbits in V).
            > subA:=CyclicSubgroups(A: OrderEqual:=20);
//Since <w^4,wa> has a cyclic subgroup of order 20.
            > k:=#subA;
            > if k ge 1 then;
                > print i;
                > end if;
            > end if;
        > end if;
    > end for;
1
5
//Thus groups #1 and #5 are possible candidates for T:<w^4,wa>.
> nsA1:=NormalSubgroups(sA1: IsCyclic:=true, OrderEqual:=20);
> nsA5:=NormalSubgroups(sA5: IsCyclic:=true, OrderEqual:=20);
//We make use of the fact that <w^4> is cyclic and normal in <w^4,wa>.
> #nsA1;
1
> #nsA5;
1
> N1:=nsA1[1]`subgroup;
> N5:=nsA5[1]`subgroup;
> F1:=sA1/N1;
> F5:=sA5/N1;
//Taking quotients. Note that <w^4,wa>/<w^4> is a cyclic
//subgroup of order 4.
> #F1;
4
> #F5;
4
> IsCyclic(F1);
false
> IsCyclic(F5);
true
//Thus A5 is our desired group T:<w^4,wa>.
> SetO:=Orbit(sA5,2);
> OT:=OrbitalGraph(A5,1,SetO);
//Constructing the corresponding TGPaley graph.
> IsIsomorphic(OG111,OT);
true
//This shows that the 3 orbital graphs (OG111, OG100, OG107) are
//isomorphic to a TGPaley graph.
\end{verbatim}

\vs \vs
\noindent \textbf{Example A4-4(ii)}: Lines $(2)-(3)$ of Table~\ref{t23}
($\Gamma_i \cong G(7^2,6)$ and $\Gamma_i \cong G(7^2,2)$).
\begin{verbatim}
> g:=AGL(2,7);
> h:=Stabilizer(ASL(2,3),1);
> D:=Find(g,h,7);
> #D;
1176
> N:=Normalizer(g,D);
> Gen(N,h);
1 (i)= 12  Order of M= 392  Order of M_0 = 8  No. of orbits in V^*= 6
2 (i)= 20  Order of M= 1176  Order of M_0 = 24  No. of orbits in V^*= 2
3 (i)= 21  Order of M= 1176  Order of M_0 = 24  No. of orbits in V^*= 2
4 (i)= 28  Order of M= 3528  Order of M_0 = 72  No. of orbits in V^*= 2
=====================
> Sub:=Subgroups(N: IsTransitive:=true);
> s12:=Sub[12]`subgroup;
> s20:=Sub[20]`subgroup;
> s21:=Sub[21]`subgroup;
> s28:=Sub[28]`subgroup;
> Set12:=Orbit(Stabilizer(s12,1),2);
> Set20:=Orbit(Stabilizer(s20,1),2);
> Set21:=Orbit(Stabilizer(s21,1),2);
> Set28:=Orbit(Stabilizer(s28,1),2);
> O12:=OrbitalGraph(s12,1,Set12);
> O20:=OrbitalGraph(s20,1,Set20);
> O21:=OrbitalGraph(s21,1,Set21);
> O28:=OrbitalGraph(s28,1,Set28);
> IsIsomorphic(O20,O21);
true
> IsIsomorphic(O20,O28);
true
//Thus the graphs O20, O21 and O28 are pair-wise isomorphic.
//Note that by using the same approach as in Example 4(i), it can
//also be shown that they are all TGPaley graphs.
\end{verbatim}

\vs \vs \noindent \textbf{Example A4-4(iii)}: Test for isomorphism in
Line $(6)$ of Table~\ref{t23} ($G(11^2,5)$).
\begin{verbatim}
> O:=GPaley(11,2,5);
//Construct generalised Paley graph with p=11, R=2, k=5.
> A:=AutomorphismGroup(O);
> A;
Permutation group A acting on a set of cardinality 121
Order = 5808 = 2^4 * 3 * 11^2
> g:=AGL(2,11);
> h:=Stabilizer(ASL(2,3),1);
> D:=Find(g,h,11);
> #D;
2904
> N:=Normalizer(g,D);
> #N;
29040
> Gen(N,h);
1 (i)= 14   Order of M= 968   Order of M_0 = 8   No. of orbits in V^*= 15
2 (i)= 24   Order of M= 2904   Order of M_0 = 24   No. of orbits in V^*= 5
3 (i)= 25   Order of M= 4840   Order of M_0 = 40   No. of orbits in V^*= 3
4 (i)= 29   Order of M= 5808   Order of M_0 = 48   No. of orbits in V^*= 5
=====================
> SUB:=Subgroups(N: IsTransitive:=true);
> M24:=SUB[24]`subgroup;
//Let M24 be the #24 group in
//"Subgroups(N: IsTransitive:=true)".
> M29:=SUB[29]`subgroup;
> Set24:=Orbit(Stabilizer(M24,1),2);
> Set29:=Orbit(Stabilizer(M29,1),2);
> O24:=OrbitalGraph(M24,1,Set24);
> O29:=OrbitalGraph(M29,1,Set29);
> ISM:=IsIsomorphic(O24,O29);
> ISM;
true
//This shows that the orbital graphs for M24 and M29 are isomorphic.
> ISM:=IsIsomorphic(O,O29);
> ISM;
false
//This shows that the orbital graphs for M24 and M29 are not isomorphic
//to a GPaley graph with corresponding parameters. We do not need to
//check if they are isomorphic to a TGPaley graphs since the index k of
//the homogeneous factorisation is odd (homogeneous factorisation into
//TGPaley graphs always has even index).
\end{verbatim}

\vs \vs \noindent \textbf{Example A4-4(iv)}: To show that the orbital graph ($G(23^2,6)$) in Line
$(8)$ of Table~\ref{t23} is isomorphic to a $K$-orbital graph where $K < A\Gamma L(1,23^2)$. (We
can employ a similar approach for Line $(7)$ of Table~\ref{t23} ($G(23^2,66)$)).
\begin{verbatim}
> g:=AGL(2,23);
> h:=ASL(2,3);
> SL:=Stabilizer(h,1);
> D:=Find(g,h,23);
> #D;
12696
> N:=Normalizer(g,D);
> Gen(N,SL);
1 (i)= 10  Order of M= 4232  Order of M_0= 8  No. of orbits in V^*= 66
2 (i)= 18  Order of M= 12696  Order of M_0= 24  No. of orbits in V^*= 22
3 (i)= 19  Order of M= 46552  Order of M_0= 88  No. of orbits in V^*= 6
4 (i)= 23  Order of M= 25392  Order of M_0= 48  No. of orbits in V^*= 11
5 (i)= 25  Order of M= 139656  Order of M_0= 264  No. of orbits in V^*= 2
================
>
> Sub:=Subgroups(N: IsTransitive:=true);
> M19:=Sub[19]`subgroup;
//Note that our desired group has 6 orbits in V^*.
> sM19:=Stabilizer(M19,1);
> Set19:=Orbit(sM19,2);
> O19:=OrbitalGraph(M19,1,Set19);
> A19:=AutomorphismGroup(O19);
> #A19;
46552
//Computing the full automorphism group of O19 (or the graph G(23^2,6)).
> OGP:=GPaley(23,2,6);
//Constructing the corresponding GPaley graph (note that we do not need
//to consider the TGPaley case; see explanation in Section 5.2.1).
> IsIsomorphic(O19,OGP);
false
//Thus our orbital graph is not a generalised Paley graph.
> G:=AGammaL(1,23^2);
> SubG:=Subgroups(G: OrderEqual:=46552);
//We want to find a subgroup of AGammaL(1,23^2) that has the same
//order as the full automorphism group of O19, whose point stabiliser
//is non-cyclic (since the graph is non-GPaley) and has 6 orbits in V^*.
> for i:= 1 to #SubG by 1 do;
    > A:=SubG[i]`subgroup;
    > sA:=Stabilizer(A,1);
    > TC:=IsCyclic(sA);
    > if TC eq false then;
        > num:=#Orbits(sA);
        > if num eq 7 then;
            > print i;
            > end if;
        > end if;
    > end for;
3
> A:=SubG[3]`subgroup;
> sA:=Stabilizer(A,1);
> #Orbit(sA,2);
88
> Set:=Orbit(sA,2);
> OA:=OrbitalGraph(A,1,Set);
> IsIsomorphic(O19,OA);
true
//Thus O19 (or the graph G(23^2,6)) is isomorphic to a K-orbital graph
//where K < AGammaL(1,23^2).
\end{verbatim}

\begin{verbatim}

\end{verbatim}

\begin{verbatim}

\end{verbatim}

\addcontentsline{toc}{chapter}{\textbf{Bibliography}}
\bibstyle{plain}




\end{document}